\documentclass[12pt, reqno]{amsart}
\usepackage{amsmath}
\usepackage{amssymb}
\usepackage{epsfig}
\usepackage{graphicx}
\usepackage{color}
\usepackage{fullpage}
\definecolor{shadecolor}{gray}{0.875}
\usepackage{amscd}
\usepackage{comment}
\usepackage{enumitem}

\numberwithin{equation}{section}

\input xy
\xyoption{all}

\calclayout
\allowdisplaybreaks[3]

\theoremstyle{plain}
\newtheorem{prop}{Proposition}[section]

\newtheorem{theo}[prop]{Theorem}
\newtheorem{coro}[prop]{Corollary}

\newtheorem{lemm}[prop]{Lemma}

\theoremstyle{definition}
\newtheorem{defi}[prop]{Definition}

\newtheorem{conj}[prop]{Conjecture}

\newtheorem{rema}[prop]{Remark}

\newtheorem{prin}[prop]{Principle}
\newtheorem{exam}[prop]{Example}

\newtheorem{clai}[prop]{Claim}
\newtheorem{obse}[prop]{Observation}
\newtheorem{cons}[prop]{Construction}

\def\bP{{\mathbb P}}

\def\bZ{{\mathbb Z}}

\def\Br{\mathrm{Br}}

\def\Eff{\overline{\mathrm{Eff}}}
\def\Pic{\mathrm{Pic}}

\def\Ext{\mathrm{Ext}}

\def\Mor{\mathrm{Mor}}
\def\Aut{\mathrm{Aut}}
\def\Nef{\mathrm{Nef}}
\def\Hilb{\mathrm{Hilb}}

\def\Pic{\mathrm{Pic}}
\def\neg{\mathrm{neg}}
\def\Sec{\mathrm{Sec}}
\def\Manin{\mathrm{Manin}}

\def\Jac{\mathrm{Jac}}
\def\maxdef{\mathrm{maxdef}}
\def\MBB{\mathrm{MBBbound}}

\def\CH{\mathrm{CH}}
\def\Sing{\mathrm{Sing}}

\makeatother
\makeatletter

\author{Brian Lehmann}
\address{Department of Mathematics \\
Boston College  \\
Chestnut Hill, MA \, \, 02467}
\email{lehmannb@bc.edu}

\author{Sho Tanimoto}
\address{Graduate School of Mathematics, Nagoya University, Furocho Chikusa-ku, Nagoya, 464-8602, Japan}
\email{sho.tanimoto@math.nagoya-u.ac.jp}

\title[Classifying sections of del Pezzo fibrations]{Classifying sections of del Pezzo fibrations, II}

\begin{document}
\date{\today}

\begin{abstract}
Let $X$ be a del Pezzo surface over the function field of a complex curve.  We study the behavior of rational points on $X$ leading to bounds on the counting function in Geometric Manin's Conjecture.  A key tool is the Movable Bend and Break Lemma which yields an inductive approach to classifying relatively free sections for a del Pezzo fibration over a curve. 
Using this lemma we prove Geometric Manin's Conjecture for certain split del Pezzo surfaces of degree $\geq 2$ admitting a birational morphism to $\mathbb P^2$ over the ground field.
\end{abstract}

\maketitle

\section{Introduction}

In an unpublished note (\cite{Bat88}), Batyrev developed a heuristic argument for Manin's Conjecture for a trivial family of Fano varieties over an $\mathbb{F}_{q}$-curve.
His heuristic relies on several assumptions about the geometry of sections of a Fano fibration over an $\mathbb{F}_{q}$-curve as well as some point counting estimates.  If we focus on the geometry of sections, the collection of analogous conjectures over the function field of a complex curve is known as Geometric Manin's Conjecture.  In this paper we study Geometric Manin's Conjecture for del Pezzo surfaces over the function field of a complex curve.

Let $B$ be a smooth projective curve over an algebraically closed field of characteristic $0$.  A del Pezzo fibration over $B$ is an algebraic fiber space $\pi: \mathcal{X} \to B$ such that $\mathcal{X}$ is projective, $\mathcal X$ has only Gorenstein terminal singularities, and the general fiber of $\pi$ is a smooth del Pezzo surface.  Let $\Sec(\mathcal{X}/B)$ denote the parameter space of sections of $\pi$.  Our main results describe the relationship between irreducible components of $\Sec(\mathcal{X}/B)$
and the geometric invariants used in Geometric Manin's Conjecture.

In the companion paper \cite{LTdPI}, we studied the analogous problem for del Pezzo fibrations over $\mathbb{P}^{1}$.  The conceptual approach in this paper is similar (particularly in Theorem \ref{theo:toomuchdeforming2} and in Section \ref{sect:genpoints}): we systematically use Bend-and-Break to reduce questions about sections to properties of rational curves in the fibers.  \cite{LTdPI} explains how our main results can be applied to study the Abel-Jacobi map for components of $\Sec(\mathcal{X}/B)$ and the enumerativity of certain Gromov-Witten invariants; we will not explain these applications in this paper but refer the reader to \cite{LTdPI} for details.

\subsection{Main results}
The following definition identifies the ``well-behaved'' sections of a del Pezzo fibration.

\begin{defi}
Let $\pi: \mathcal{X} \to B$ be a del Pezzo fibration.  We say that a section $C$ is relatively free if $C$ is contained in the smooth locus of $\mathcal X$, $H^{1}(C,N_{C/\mathcal{X}}) = 0$, and $N_{C/\mathcal{X}}$ is generically globally generated.
\end{defi}

\begin{rema}
\cite{GHS03} shows that every del Pezzo fibration $\pi: \mathcal{X} \to B$ admits a section.  Since a Gorenstein terminal threefold has only lci singularities, starting from any given section we can construct a relatively free section by gluing on sufficiently many free rational curves contained in the fibers of $\pi$ and smoothing (see Remark \ref{rema:enoughsections}). 
Thus every del Pezzo fibration admits many families of relatively free sections.
\end{rema}

It is natural to separate irreducible components of $\Sec(\mathcal{X}/B)$ into three types:
\begin{enumerate}
\item Non-dominant families of sections.
\item Dominant families for which a general section is not relatively free.
\item Dominant families for which a general section is relatively free.
\end{enumerate}
We prove structural theorems for all three types of irreducible components.  We first analyze the components of $\Sec(\mathcal{X}/B)$ which parametrize a non-dominant family of sections.  %The following theorem shows that such sections are usually contained in the exceptional set for Geometric Manin's Conjecture (and if they are not, then they give a negligible contribution to the counting function).

\begin{theo} \label{theo:maintheorem1}
Let $\pi: \mathcal{X} \to B$ be a del Pezzo fibration such that $-K_{\mathcal X/B}$ is relatively nef.   There is a finite union of surfaces $Y$ whose intersections with a general fiber of $\pi$ are rational curves of anticanonical degree $\leq2$ and a constant $C(\mathcal{X})$ such that the following holds.  Suppose that $M \subset \Sec(\mathcal{X}/B)$ is a component parametrizing a non-dominant family of sections of height $\geq C(\mathcal{X})$.  Then the sections parametrized by $M$ sweep out one of these surfaces $Y$.
\end{theo}

Here $C(\mathcal{X})$ is an explicit constant determined by the behavior of low degree sections.

\begin{rema}
\cite[3.3 Theorem]{Corti} shows that any del Pezzo surface over the function field of $B$ admits an integral model with Gorenstein terminal singularities such that $-K_{\mathcal{X}/B}$ is relatively nef.  Thus every del Pezzo fibration admits a birational model where Theorem \ref{theo:maintheorem1} applies.
\end{rema}

\begin{rema}
Theorem \ref{theo:maintheorem1} can be used to show that there is a proper closed set $V \subsetneq \mathcal{X}$ which contains all sections which do not deform to dominate $\mathcal{X}$.
\end{rema}

Our second main result addresses dominant families of sections which are not relatively free.  It shows that the existence of sections of this type is controlled by a bounded family of surfaces.  %The following theorem shows that such sections usually give a negligible contribution to the counting function in Geometric Manin's Conjecture (and if they do not, then they are contained in the exceptional set).

\begin{theo} \label{theo:maintheorem2}
Let $\pi: \mathcal{X} \to B$ be a del Pezzo fibration such that $-K_{\mathcal X/B}$ is relatively nef.   There is a bounded family of surfaces $Y$ whose intersections with a general fiber of $\pi$ are rational curves of anticanonical degree $2$ and a constant $D(\mathcal{X})$  such that the following holds.  Suppose that $M \subset \Sec(\mathcal{X}/B)$ is a component parametrizing a dominant family of sections of height $\geq D(\mathcal{X})$ which are not relatively free.  Then the general section parametrized by $M$ is contained in some surface  in our family.
\end{theo}

Here $D(\mathcal{X})$ is an explicit constant determined by the behavior of low degree sections.

Our final main result shows that any component of $\Sec(\mathcal{X}/B)$ that parametrizes relatively free sections of large height will contain in its closure a union of a $\pi$-vertical rational free curve and a relatively free section of smaller height.  As in \cite{HRS04}, this allows us to study families of relatively free sections of large height via induction by appealing to the gluing-and-smoothing structure of stable maps.

\begin{theo}[Movable Bend-and-Break for relatively free sections] \label{theo:maintheorem3}
Let $\pi: \mathcal{X} \to B$ be a del Pezzo fibration such that $-K_{\mathcal{X}/B}$ is relatively ample.  There is a constant $Q(\mathcal{X})$ satisfying the following property.  Suppose that $M \subset \Sec(\mathcal{X}/B)$ is a component that parametrizes a dominant family of relatively free sections $C$ satisfying $-K_{\mathcal{X}/B} \cdot C \geq Q(\mathcal{X})$.  Then the closure of $M$ in $\overline{M}_{g(B),0}(\mathcal{X})$ contains a point representing a stable map whose domain has two components; one maps birationally to a relatively free section and the other maps birationally to a free $\pi$-vertical rational curve. 
\end{theo}

Here $Q(\mathcal{X})$ is an explicit constant determined by the behavior of low degree sections.  

\begin{rema}
\cite[1.10 Theorem]{Corti} shows that any del Pezzo fibration whose generic fiber has degree $\geq 3$ will admit a birational model which has Gorenstein terminal singularities and a relatively ample anticanonical divisor.
\end{rema}

\subsection{Geometric Manin's Conjecture}

Inspired by the thin set version of Manin's conjecture and the conjectural description of the exceptional set in \cite{LST18}, the authors proposed the first version of Geometric Manin's Conjecture in \cite{LT17}.  The statement relies on the following invariant from the Minimal Model Program.

\begin{defi}
Let $X$ be a smooth projective variety over a field of characteristic $0$.  Let $L$ be a big and nef $\mathbb{Q}$-Cartier divisor on $X$.  We define the Fujita invariant, or the $a$-invariant, to be
\begin{equation} \label{eq: ainv}
a(X, L) = \min \{ t \in \mathbb R \mid K_{X} + tL \in \overline{\mathrm{Eff}}^1(X)\}.
\end{equation}
When $L$ is nef but not big, we formally set $a(X, L) = +\infty$.

When $X$ is singular, we define the Fujita invariant as the Fujita invariant of the pullback of $L$ to any smooth model. This is well-defined because of \cite[Proposition 2.7]{HTT15}.
\end{defi}

Roughly speaking Geometric Manin's Conjecture predicts two things.  First, the conjecture predicts that there should be a ``thin exceptional set'' which can be described using the Fujita invariant as in \cite{LST18}.

\begin{prin}
\label{prin:exceptionalset}
Let $\pi: \mathcal{X} \to B$ be a Fano fibration over a smooth projective curve.  All ``pathological'' irreducible components of $\Sec(\mathcal{X}/B)$ will parametrize sections which are contained in a bounded family of subvarieties $\mathcal{Y} \subset \mathcal{X}$ such that the Fujita invariant of $\mathcal{Y}_{\eta}$ with respect to $-K_{\mathcal{X}/B}$ is at least as large as the Fujita invariant of $\mathcal{X}_{\eta}$.
\end{prin}

Here the notion of ``pathological'' is flexible.  If we interpret patholgical to mean ``not relatively free'', then Principle \ref{prin:exceptionalset} for del Pezzo fibrations is established by Theorem \ref{theo:maintheorem1} and Theorem \ref{theo:maintheorem2}.  We can also interpret pathological more narrowly to mean ``components which grow too quickly'' and we give a precise accounting of Principle \ref{prin:exceptionalset} for del Pezzo fibrations in this setting in Section \ref{sect:maninconj}.

Second, Geometric Manin's Conjecture predicts that the number of relatively free families representing a given numerical class is bounded above.  

\begin{prin}
\label{prin:GMC}
Let $\pi: \mathcal{X} \to B$ be a Fano fibration over a smooth projective curve.  After removing all components of $\Sec(\mathcal{X}/B)$ which factor through the exceptional set constructed by \cite{LST18}, there should be exactly $|\mathrm{Br}(X)|$ irreducible components of $\Sec(\mathcal{X}/B)$ representing each sufficiently positive nef curve class which admits a relatively free section.
 \end{prin}
 
We do not know how to establish Principle \ref{prin:GMC} for del Pezzo fibrations in general. However, the inductive structure given by Theorem \ref{theo:maintheorem3} allows us to prove a weaker bound on the number of components.  The following result shows that the number of components of $\Sec(\mathcal{X}/B)$ only grows polynomially in the degree as predicted by Batyrev.

\begin{theo} \label{theo:maintheorembatyrev}
Let $\pi: \mathcal{X} \to B$ be a del Pezzo fibration such that $-K_{\mathcal{X}/B}$ is relatively ample.  Then there is a polynomial $P(d)$ such that the number of components of $\Sec(\mathcal{X}/B)$ parametrizing sections of height $\leq d$ is bounded above by $P(d)$.
\end{theo}

One can rephrase Geometric Manin's Conjecture to emphasize the analogy with the number theoretic version. Let $\pi : \mathcal X \to B$ be a del Pezzo fibration. For simplicity, we assume that $\rho(\mathcal X_\eta)\geq 2$ where $\rho(\mathcal X_\eta)$ is the Picard rank of $\mathcal X_\eta$. We say that a component $M\subset \Sec(\mathcal X/B)$ is an accumulating component if $M$ parametrizes sections contained in a surface whose generic fiber is a union of $-K_{\mathcal X_{\eta}}$-lines. We say that a component $M\subset \Sec(\mathcal X/B)$ is a Manin component if it is not accumulating.

Let $\mathrm{Manin}_i$ be the set of Manin components parametrizing sections $C$ with $-K_{\mathcal X/B}.C = i$. For $q > 1$ and any positive integer $d$ we define the counting function by
\begin{equation*}
N(\mathcal{X},-K_{\mathcal{X}/B},q,d) := \sum_{i = 1}^{d} \sum_{M \in \Manin_{i}} q^{\dim M}.
\end{equation*}
This counting function is inspired by Batyrev's heuristic for Manin's Conjecture for $\mathbb F_q$-curves on a smooth Fano variety; the term $q^{\dim M}$ represents the ``expected'' number of $\mathbb{F}_{q}$-points on $M$.  The ``expected'' asymptotic growth rate of the counting function is
\begin{equation*}
N(\mathcal{X},-K_{\mathcal{X}/B},q,d) \mathrel{\mathop{\sim}_{d \to \infty}} c q^{d} d^{\rho(\mathcal{X}_{\eta})-1}.
\end{equation*}
Theorem \ref{theo:maintheorembatyrev} implies a weaker upper bound on the counting function.  For simplicity we assume that the general fiber of $\pi$ contains a $(-1)$-curve.

\begin{theo} \label{theo:maintheoremgrowthrate}
Let $\pi: \mathcal{X} \to B$ be a del Pezzo fibration such that $\mathcal X$ is smooth, $-K_{\mathcal{X}/B}$ is relatively ample, and the general fiber is a del Pezzo surface that is not $\mathbb{P}^{2}$ or $\mathbb{P}^{1} \times \mathbb{P}^{1}$.  Then there is some non-negative integer $r$ such that
\begin{equation*}
N(\mathcal{X},-K_{\mathcal{X}/B},q,d) = O(q^{d} d^{r}).
\end{equation*}
\end{theo}

\begin{rema} 
Suppose we fix a del Pezzo surface $\mathcal{X}_{\eta}$ over the function field of $B$.  Geometric Manin's Conjecture for $\mathcal{X}_{\eta}$ implicitly relies on the choice of an integral model of $\mathcal{X}_{\eta}$ over $B$.  (Analogously, in the number field setting Manin's Conjecture depends upon a choice of metrization.)

Just as in the number field setting, if $\mathcal{X}_{\eta}$ has an integral model where Theorem \ref{theo:maintheoremgrowthrate} applies then one should be able to deduce the analogous bound for every integral model of $\mathcal{X}_{\eta}$.  (In contrast, although we expect Theorem \ref{theo:maintheorem1} and Theorem \ref{theo:maintheorem3} to hold for any integral model of $\mathcal{X}_{\eta}$ we do not see an easy way to pass the statements between different integral models.)
\end{rema}

\subsection{Classifying components of $\Sec(\mathcal{X}/B)$}

Suppose that $\pi: \mathcal{X} \to B$ is a del Pezzo fibration such that $-K_{\mathcal{X}/B}$ is relatively ample.  Furthermore, suppose that we can classify all sections of low height by hand.  Then the components of $\Sec(\mathcal{X}/B)$ of large height can be described inductively using Theorem \ref{theo:maintheorem1}, Theorem \ref{theo:maintheorem2}, and Theorem \ref{theo:maintheorem3}.  This technique is illustrated in several examples in \cite[Section 8]{LTdPI}.  In particular, this type of argument should allow us to fully prove Geometric Manin's Conjecture for specific examples.
However it is challenging to prove a general statement.  We are able to prove this strong version in a somewhat restrictive situation.

\begin{theo} \label{theo:maintheorem4}
Let $\pi: \mathcal{X} \to B$ be a del Pezzo fibration such that every fiber of $\pi$ is a del Pezzo surface with canonical singularities of degree $\geq 2$.  Suppose that $\mathcal{X}$ is $\mathbb{Q}$-factorial and that the generic fiber $\mathcal{X}_{\eta}$ is the blow-up of $\mathbb{P}^{2}_{k(B)}$ at a finite number of points defined over $k(B)$.  Then there is some section $C$ such that for every integral numerical class $\alpha$ in $C + \Nef_{1}(\mathcal{X}_{\eta})$ there is exactly one component of $\Sec(\mathcal{X}/B)$ which parametrizes relatively free sections and represents $\alpha$.
\end{theo}

\begin{exam}
Suppose we fix a smooth del Pezzo surface $S$ of degree $\geq 2$ and a curve $B$.  Applying Theorem \ref{theo:maintheorem4} to $\pi: S \times B \to B$ we see there is a translate of $\Nef_{1}(S)$ in $N_{1}(S)$ such that every curve class in this translate is represented by a unique component of $\Mor(B,S)$. 
\end{exam}

Theorem \ref{theo:maintheorem4} verifies the strong version of Geometric Manin's Conjecture for this type of del Pezzo fibration.  Note that in this situation $\mathcal X_\eta$ is rational so that $\mathrm{Br}(\mathcal X)$ is trivial.  For the del Pezzo fibrations addressed by Theorem \ref{theo:maintheorem4} we obtain the expected growth rate
\begin{equation*}
N(\mathcal{X},-K_{\mathcal{X}/B},q,d) \mathrel{\mathop{\sim}_{d \to \infty}} c q^{d} d^{\rho(\mathcal{X}_{\eta})-1}.
\end{equation*}
where the leading constant $c$ is described by Theorem~\ref{theo:asymptoticformula}.

A key new feature of Theorem \ref{theo:maintheorem4} is that the proof does not rely on induction.  In particular, we do not need any ad hoc arguments to analyze base cases. Instead, we define a monoid structure on the set $\mathcal R$ of components of nef curve classes on a del Pezzo surface given by gluing and smoothing. This monoid acts on the set of components of relatively free sections of sufficiently large height by gluing and smoothing and Movable Bend and Break tells us that this set is finitely generated by $\mathcal R$. By exhibiting many relations among components of relatively free sections, we show that every nef class in a translate of the cone is represented by one irreducible relatively free component.

\subsection{Comparison to previous works}
\label{subsec:history}
Batyrev developed a heuristic for Manin's Conjecture over global function fields in the unpublished notes \cite{Bat88}.  (This heuristic is explained in \cite[Section 4.7]{Tschinkel09} and \cite[Section 1.2]{Bou}.)  This perspective motivated the formulation of the Batyrev-Manin Conjecture developed in \cite{BM}.  This theme was also revisited in \cite{Manin95} which proves an exponential bound on the number of components of $\Mor(\bP^{1},X)$ and discusses the analogy with rational point counts.

Suppose that $X$ is a Fano variety and $B$ is a smooth projective curve over a finite field. 
Batyrev's heuristic for sections of $\pi: X \times B \to B$ is based on the following three assumptions:
\begin{enumerate}
\item For each nef class $\alpha \in \mathrm{Nef}_1(X)_{\mathbb Z}$, every irreducible component of the moduli space $\mathrm{Mor}(B, X, \alpha)$ has the expected dimension $-K_{X} \cdot B + \dim(X)(1-g(B))$;
\item For each nef class $\alpha \in \mathrm{Nef}_1(X)_{\mathbb Z}$, $\mathrm{Mor}(B, X, \alpha)$ is irreducible;
\item If we use the naive estimate $|\mathrm{Mor}(B, X, \alpha)(\mathbb F_q)| \approx q^{\dim \mathrm{Mor}(B, X, \alpha)}$ then we do not affect the asymptotic growth rate of the number of rational points.
\end{enumerate}
As stated the assumptions (1) and (2) are not valid even for large degree curves.  This failure is related to the properties of the exceptional set in Manin's Conjecture for rational points.  
(Ellenberg and Venkatesh suggest that (3) may be approached using homological stability -- see \cite{EV05} -- but as of now we do not have many examples. \cite{BS20} proves a statement in this direction for low degree Fano hypersurfaces.)

The classification of components of $\Mor(\mathbb{P}^{1},X)$ for a Fano variety $X$ has a long and rich history.  For $\mathbb{P}^{1}$ the predictions (1),(2) implicit in Batyrev's heuristics have been verified in the following cases (and in many other special cases):
\begin{itemize}
\item generalized flag varieties (\cite{Thomsen98} and \cite{KP01});
\item most Fano hypersurfaces (\cite{HRS04}, \cite{RY16}, and \cite{BV16});
\item toric varieties (\cite{Bou09} and \cite{Bou16});
\item moduli spaces of vector bundles on curves (\cite{Cas04} and \cite{BM20}); 
\item smooth del Pezzo surfaces (\cite{Testa09}), and;
\item Fano threefolds (\cite{CS09}, \cite{Cas04}, \cite{LT17}, \cite{LT18}, and \cite{BLRT20}).
\end{itemize}
In contrast, there has not been much progress toward the classification of irreducible components of $\Mor(B,X)$ when $B$ has genus $\geq 1$.  Aside from Bourqui's pioneering results for toric varieties, to the best of our knowledge there are only a few examples of homogeneous varieties $X$ where the irreducible components of $\Mor(B,X)$ have been fully classified (and primarily when $B$ is an elliptic curve  -- \cite{Bruguieres87}, \cite{Ballico89}, \cite{Perrin12}, \cite{PP13}).

More generally, one would like to classify the irreducible components of the space of sections of a Fano fibration. While there are many powerful theorems describing the qualitative nature of sections of Fano fibrations, there are fewer quantitative results.  \cite{LTdPI} develops Geometric Manin's Conjecture for sections of Fano fibrations over $\mathbb{P}^{1}$ and proves some results for del Pezzo fibrations.  In this paper we address del Pezzo fibrations over curves of higher genus.  The key contributions of this paper are:

\begin{itemize}
\item This is the first paper which uses the inductive strategy of \cite{HRS04} to classify curves of high genus on a Fano variety.  The framework and perspective we develop should be applicable in other situations.  When working with higher genus curves, there are many new technical obstacles -- in particular, the most difficult part of the argument concerns dominant families of sections which are not relatively free.
\item For a special class of del Pezzo fibrations, we give a conceptual proof of Geometric Manin's Conjecture using a gluing-and-smoothing monoid structure (Theorem~\ref{theo:maintheorem4}).  Previous techniques have tended to be somewhat ad hoc.
\end{itemize}

\bigskip

\noindent
{\bf Acknowledgements:}
The authors would like to thank Yuri Tschinkel for answering our question regarding del Pezzo surfaces over non-closed fields.
The authors would also like to thank the anonymous referees for carefully reading our manuscript and significantly improving the exposition of the paper.

Brian Lehmann was supported by NSF grant 1600875.  Sho Tanimoto was partially supported by Inamori Foundation, by JSPS KAKENHI Early-Career Scientists Grant number 19K14512,  by JSPS Bilateral Joint Research Projects Grant number JPJSBP120219935, and by MEXT Japan, Leading Initiative for Excellent Young Researchers (LEADER).

\section{Preliminaries}
\label{sec:preliminaries}

Let $k$ be an algebraically closed field of characteristic $0$. Let $B$ be a smooth projective curve defined over $k$.  In this paper our ground fields will be $k$ and $k(B)$. A variety is a reduced irreducible separated scheme which has finite type over the ground field. In this paper a component of a scheme means an irreducible component unless otherwise specified.  When we take a component of a scheme, we always endow it with its reduced structure.

\subsection{Positive cycles}
We will use $\sim_{rat}$ to denote rational equivalence of cycles, $\sim_{alg}$ to denote algebraic equivalence of cycles,
and $\equiv$ to denote numerical equivalence of cycles. 

Let $X$ be a projective variety.  We will let $N^{1}(X)$ denote the space of $\mathbb{R}$-Cartier divisors up to numerical equivalence and let $\Eff^{1}(X)$ and $\Nef^{1}(X)$ denote respectively the pseudo-effective and nef cones of divisors.  Dually, $N_{1}(X)$ denotes the space of real $1$-cycles up to numerical equivalence and $\Eff_{1}(X)$ and $\Nef_{1}(X)$ denote respectively the pseudo-effective and nef cones of curves.  We will denote by $N_{1}(X)_{\mathbb{Z}}$ the lattice of integral curve classes inside of $N_{1}(X)$ and by $N^{1}(X)_{\mathbb{Z}}$ the lattice of integral divisors inside of $N^{1}(X)$. 

We say that a reduced irreducible curve $C$ on $X$ is movable if $C$ is a member of a family of curves which dominates $X$.

\subsection{Height functions}

\begin{defi}
A Fano fibration $\pi: \mathcal{X} \to B$ is a surjective morphism with connected fibers from a Gorenstein terminal projective variety $\mathcal{X}$ to a smooth projective curve $B$ such that a general fiber is a smooth Fano variety.

We will always denote the generic point of $B$ by $\eta$ and the generic fiber of $\pi$ by $\mathcal{X}_{\eta}$.
\end{defi}

Given a section $C$ of $\pi$ and a divisor $L$ on $\mathcal{X}$, the height of $C$ with respect to $L$ is defined to be $L \cdot C$.  The following statement is the Northcott property for Fano fibrations over $B$.

\begin{lemm} \label{lemm:northcott}
Let $\pi: \mathcal{X} \to B$ be a Fano fibration.  Fix a $\mathbb{Q}$-Cartier divisor $L$ on $\mathcal{X}$ whose restriction to the generic fiber of $\pi$ is ample.  For any constant $\gamma$, the sections whose height with respect to $L$ is $\leq \gamma$ form a bounded family.
\end{lemm}

One consequence is that there is a lower bound on the possible values of $L \cdot C$ as we vary $C$ over all sections.

\begin{defi} \label{defi:neg}
Let $\pi: \mathcal{X} \to B$ be a Fano fibration and let $L$ be a $\mathbb{Q}$-Cartier divisor on $\mathcal{X}$ such that the restriction of $L$ to the generic fiber is  ample.  We define $\neg(\mathcal{X},L)$ to be the smallest value of $L \cdot C$ as we vary $C$ over all sections of $\pi: \mathcal{X} \to B$.
\end{defi}

Another consequence is captured by the following observation:

\begin{obse} \label{obse:goingthroughgenpoints}
Let $\pi: \mathcal{X} \to B$ be a Fano fibration.  Fix an irreducible component $M$ of $\Sec(\mathcal{X}/B)$.  We let $M^{[m]}$ denote the family of sections parametrized by $M$ with $m$ marked points and we denote the evaluation map by $\mathrm{ev}_m: M^{[m]} \to \mathcal{X}^m$.

Fix an integer $d$ and a relatively ample divisor $L$.  Suppose that we vary $M$ over all components of $\Sec(\mathcal{X}/B)$ of $L$-degree $\leq d$ for which the evaluation map $\mathrm{ev}_m$ does not map dominantly to $\mathcal{X}^{\times m}$.  By Lemma \ref{lemm:northcott} the union of closures of the images of these maps will be a proper closed subset of $\mathcal{X}^{\times m}$.  In particular, if we fix $m$ general points of $\mathcal{X}$, then any irreducible component $M$ of $\Sec(\mathcal{X}/B)$ of degree $\leq d$ which parametrizes a section through these points will yield a dominant map onto $\mathcal{X}^{\times m}$.
\end{obse}

\begin{rema}
From now on, when we say that an irreducible component $M$ of $\Sec(\mathcal{X}/B)$ parametrizes sections through $m$ general points, we will mean ``general'' in the sense of Observation \ref{obse:goingthroughgenpoints}.  Note that implicitly the meaning of ``general'' will depend upon the relative anticanonical degree of the curves parametrized by $M$ even when this degree is not explicitly mentioned.
\end{rema}

\subsection{Generic Fujita invariants}

Let us recall the definition of the Fujita invariant from the introduction.

\begin{defi}
Let $X$ be a smooth projective variety over a field of characteristic $0$.  Let $L$ be a big and nef $\mathbb{Q}$-Cartier divisor on $X$.  We define the Fujita invariant, or the $a$-invariant, to be
\begin{equation} \label{eq: ainv}
a(X, L) = \min \{ t \in \mathbb R \mid K_{X} + tL \in \overline{\mathrm{Eff}}^1(X)\}.
\end{equation}
When $L$ is nef but not big, we formally set $a(X, L) = +\infty$.

When $X$ is singular, we define the Fujita invariant as the Fujita invariant of the pullback of $L$ to any smooth model. This is well-defined because of \cite[Proposition 2.7]{HTT15}.
\end{defi}

Note that the $a$-invariant is geometric: it does not change under field extension.  We will be interested in how the Fujita invariant behaves over $k(B)$. 

\begin{lemm}[\cite{LTdPI} Lemma 3.3]
Let $\pi: \mathcal{X} \to B$ be a Fano fibration and let $L$ be a $\mathbb{Q}$-Cartier divisor on $\mathcal X$ such that the restriction of $L$ to the generic fiber $\mathcal X_\eta$ is big and nef.  Then for any smooth Fano fiber $F$ of $\pi$ we have
\begin{equation*}
a(X_{\eta},L|_{X_{\eta}}) = a(F,L|_{F}).
\end{equation*}
\end{lemm}

As in \cite{LTdPI} we will call this quantity the generic $a$-invariant of $\mathcal X$ with respect to $L$.  For del Pezzo surfaces, it is easy to work out the behavior of the $a$-invariant of the anticanonical divisor when restricted to subvarieties.  This leads to the following description:

\begin{lemm}[\cite{LTdPI} Lemma 3.4] \label{lemm:genericainvfordp}
Let $\pi: \mathcal{X} \to B$ be a del Pezzo fibration.  Then:
\begin{itemize}
\item A subvariety $Y$ will have $a(Y_{\eta},-K_{\mathcal{X}/B}) > 1$ if and only if its intersection with a general fiber $F$ is a union of curves of the following types: $(-1)$-curves, or rational curves in $|-K_F|$ when $F$ has degree $1$.
\item A subvariety $Y$ will have $a(Y_{\eta},-K_{\mathcal{X}/B}) = 1$ if and only if its intersection with a general fiber $F$ is a union of curves of the following types: irreducible fibers of a conic fibration on $F$, the rational curves in $|-K_{F}|$ if $F$ has degree $2$, and the rational curves which lie in $|-2K_{F}|$ or the pullback of the anticanonical linear series on a degree $2$ del Pezzo surface if $F$ has degree $1$.
\end{itemize}
\end{lemm}

\begin{coro} \label{coro:closedainv}
Let $\pi: \mathcal{X} \to B$ be a del Pezzo fibration.  The union of all subvarieties $Y$ with  $a(Y_{\eta},-K_{\mathcal{X}/B}) > 1$ is a closed subset of $\mathcal{X}$.
\end{coro}

The subvarieties $Y$ with generic $a$-invariant equal to $1$ are a little more complicated; note that they need not form a bounded family on $\mathcal{X}$ (even though the corresponding subvarieties of $\mathcal{X}_{\eta}$ do form a bounded family).  However, we do have a weaker boundedness statement.

\begin{prop} \label{prop:boundednessequalainv}
Let $\pi: \mathcal{X} \to B$ be a del Pezzo fibration.  Fix a bounded family of sections of $\pi$.  The family of surfaces $Y \subset \mathcal{X}$ which contain a section in our family and satisfy $a(Y_{\eta},-K_{\mathcal{X}/B}) = 1$ is bounded.
\end{prop}

\begin{proof}
Lemma \ref{lemm:genericainvfordp} (2) describes the possible types of a general fiber of $\pi|_{Y}$: they are rational curves of anticanonical degree $2$.  Note that if we fix a point in a del Pezzo surface there are only finitely many curves of these types through that point.

Suppose we fix a section $C$ of $\pi$.  The observation above shows that there are only finitely many surfaces which contain this section and satisfy $a(Y_{\eta},-K_{\mathcal{X}/B}) = 1$.  Constructing an incidence correspondence we obtain the boundedness of the surfaces $Y$ as in the statement of the proposition.
\end{proof}

\subsection{Bend-and-Break}

We will need the following result controlling the behavior of Bend-and-Break for sections.

\begin{lemm} \label{lemm:strongerbandb}
Let $\pi: \mathcal{X} \to B$ be a Fano fibration.  Fix a set of $n$ general points $\{x_{i}\}_{i=1}^{n}$ of $\mathcal{X}$ where $n \geq 2$.  Suppose there is a one-dimensional family of sections which contain all $n$ points.  Then the closure of this family in $\overline{M}_{g(B), 0}(\mathcal{X})$ contains a stable map $f: C \to \mathcal{X}$ such that $C$ has at least two components $C_{1},C_{2}$ which are not contracted by $f$ and such that $f(C_{1})$ and $f(C_{2})$ each contain one of the $x_{i}$.
\end{lemm}

\begin{proof}
When $B \cong \mathbb{P}^{1}$ this statement is proved in \cite[Lemma 4.1]{LTdPI}.  Suppose that $g(B) \geq 1$.  The statement of Bend-and-Break shows that the family of sections deforms to a stable map $f: C \to \mathcal{X}$ such that
\begin{itemize}
\item $f(C)$ contains all $n$ general points, and
\item there is some component $C_{1}$ of $C$ which is rational, is not contracted by $f$, and whose $f$-image contains one of the general points.
\end{itemize}
Since $C_{1}$ is rational and $g(B) \geq 1$ the image of $C_{1}$ must be $\pi$-vertical, so $f(C_{1})$ can contain at most one general point.  Thus the other general points will be contained in other components of $f(C)$, proving the statement.
\end{proof}

\subsection{Vector bundles on curves}

In this section we quickly review some facts about rank $2$ vector bundles on curves.

\begin{lemm}[\cite{Hartshorne77} V.2.12.(b), \cite{Nagata70} Theorem 1] \label{lemm:indecomrk2}
Let $\mathcal{E}$ be a rank $2$ vector bundle on the smooth curve $B$.
\begin{enumerate}
\item Suppose that $\mathcal{E}$ is unstable.  Consider the destabilizing exact sequence
\begin{equation*}
0 \to \mathcal{L}_{2} \to \mathcal{E} \to \mathcal{L}_{1} \to 0
\end{equation*}
given by the maximal slope subbundle $\mathcal{L}_{2}$.  If $\mathcal{E}$ is indecomposable then
\begin{equation*}
0 < \deg(\mathcal{L}_{2}) - \deg(\mathcal{L}_{1}) \leq 2g(B) - 2.
\end{equation*}
\item Suppose that $\mathcal{E}$ is semistable.  Consider the exact sequence
\begin{equation*}
0 \to \mathcal{L}_{1} \to \mathcal{E} \to \mathcal{L}_{2} \to 0
\end{equation*}
given by a maximal slope rank $1$ subbundle $\mathcal{L}_{1}$.  Then we have
\begin{equation*}
0 \leq \deg(\mathcal{L}_{2}) - \deg(\mathcal{L}_{1}) \leq g(B).
\end{equation*}
\end{enumerate}
\end{lemm}

Using Kodaira vanishing and the exact sequences above, one obtains:

\begin{coro} \label{coro:indecomposablevanishing}
Let $\mathcal{E}$ be a rank $2$ vector bundle on the smooth curve $B$.
\begin{enumerate}
\item If $\mathcal{E}$ is unstable, indecomposable, and $\deg(\mathcal{E}) \geq 6g(B) - 4$ then $h^{1}(B,\mathcal{E}) = 0$.
\item If $\mathcal{E}$ is semistable and $\deg(\mathcal{E}) \geq 5g(B) - 2$ then $h^{1}(B,\mathcal{E}) = 0$.
\end{enumerate}
\end{coro}

We will also need the following fact about the variation of the Harder-Narasimhan filtration in families. 

\begin{theo}[\cite{HL97} Theorem 2.3.2] \label{theo:hnfiltrationinfamilies}
Let $T$ be a variety and let $\mathcal{E}$ be a vector bundle of rank $2$ on $T \times B$.  We can think of $T$ as the parameter space for a family of rank $2$ vector bundles on $B$.  
Then there exists a non-empty Zariski open subset $T' \subset T$ such that either
\begin{itemize}
\item for each $t \in T'$, $\mathcal E|_{\{t\} \times B}$ is semistable, or;
\item there exists a line subbundle $\mathcal L \subset \mathcal E|_{T'\times B}$ on $T'\times B$ such that for each $t \in T'$, $\mathcal L|_{\{t\} \times B}$ is the maximal destabilizing subsheaf of $\mathcal E|_{\{t\} \times B}$. 
\end{itemize}
\end{theo}

\subsection{Conics on del Pezzo surfaces}
Finally we record the following lemma for later applications:

\begin{lemm}
\label{lemm:conicsondP}
Let $S$ be a smooth del Pezzo surface of degree $d$ over an algebraically closed field $k$ of characteristic $0$.
Let $C$ be an anticanonical conic on $S$, i.e.~$C$ is rational and $-K_S \cdot C = 2$.
Then $C$ satisfies one of the following conditions:
\begin{enumerate}
\item $C^2 = 0$ and $C$ is a member of a conic fibration $\rho : S \to \mathbb P^1$;
\item $d = 2$ and $C$ is a rational member of $|-K_S|$;
\item $d = 1$ and $C$ is a rational member of $|L|$ where $L$ is the pullback of the anticanonical divisor via a blow down $S \to S'$ to a degree $2$ del Pezzo surface, or;
\item $d = 1$ and $C$ is a rational member of $|-2K_S|$.
\end{enumerate}
Moreover a general member of each family has at worst nodal singularities.
\end{lemm}

\begin{proof}
First of all it follows from adjunction that $C^2 \geq 0$ and $C^2$ is even.
On the other hand, by the Hodge index theorem we have $dC^2 -4 \leq 0$. 
Thus when $d \geq 3$, we conclude that $C^2 = 0$ and the linear series $|C|$ defines a conic fibration.
When $d = 2$, we have two possibilities: $C^2 = 0$ or $2$. When $C^2 = 0$, $|C|$ defines a conic fibration. When $C^2 = 2$, $C$ is linearly equivalent to $-K_S$. Thus $C$ is a rational member of $|-K_S|$.
When $d = 1$, there are three possibilities: $C^2 = 0, 2$, or $4$.
Again when $C^2 = 0$ the linear series $|C|$ defines a conic fibration. When $C^2 = 2$, we claim that $C + K_S$ is linearly equivalent to a $(-1)$-curve. Indeed, using the exact sequence
\begin{equation*}
0 \to H^{0}(S,\mathcal{O}_{S}(C + K_{S})) \to H^{0}(S,\mathcal{O}_{S}(C)) \to H^{0}(C,\mathcal{O}_{C})
\end{equation*}
one sees that $H^{0}(S,C + K_{S})>0$ and then the intersection-theoretic properties of $C+K_{S}$ show that it is linearly equivalent to a $(-1)$-curve.  Thus $C$ is a member of a family in (3). When $C^2 = 4$, $C$ is linearly equivalent to $-2K_S$. Thus our assertion follows.

Finally for the last claim, any general member of a conic fibration is smooth. For a degree $2$ del Pezzo surface $S$, $|-K_S|$ defines a double cover $S \to \mathbb P^2$ ramified along a smooth quartic curve $D$ and any rational member of $|-K_S|$ is the pullback of a tangent line to $D$. A generic tangent line has only one tangency point so that a general rational curve in $|-K_{S}|$ has only one node.
This proves the claim for families in (2) and (3). Finally for a degree $1$ del Pezzo surface $S$, $|-2K_S|$ defines a double cover $S \to Q$ to a singular quadric cone ramified along a smooth complete intersection of $Q$ and a cubic surface. A general rational member of $|-2K_S|$ corresponds to a hyperplane section bitangent to the branch divisor, thus a general rational curve in $|-2K_{S}|$ has two nodes. Thus our assertion follows.
\end{proof}

\section{Families of sections}

\subsection{Moduli spaces of sections}

Let $\pi : \mathcal X \rightarrow B$ be a Fano fibration with $\dim \mathcal X \leq 3$. By definition $\mathcal X$ has only Gorenstein terminal singularities. When $\dim(\mathcal X) = 3$, this is equivalent to saying that $\mathcal X$ has cDV singularities (\cite[Corollary 5.38]{KM98}).  Such singularities are analytically isomorphic to hypersurface singularities, and in particular $\mathcal X$ is locally complete intersection in the Zariski topology. 

We let $\Sec(\mathcal{X}/B)$ denote the open subset of $\Hilb(\mathcal{X})$ parametrizing sections of $\pi$ and let $\Sec(\mathcal{X}/B)_{d}$ denote the sublocus parametrizing sections of height $d$.  Suppose that $C$ is a general curve parametrized by a component $M \subset \Sec(\mathcal{X}/B)$.  The expected dimension of $M$ is
\begin{equation*}
-K_{\mathcal{X}/B} \cdot C + (\dim \mathcal X-1)(1-g(B)).
\end{equation*}
The expected dimension is a lower bound for the dimension of $M$ by \cite[Chapter I.2 Theorem I.2.15.2]{Kollar}. Indeed, let $\mathcal I_C \subset \mathcal O_{\mathcal X}$ be the ideal sheaf of $C$. Using the fact that $\mathcal X$ is lci, one can conclude that we have the exact sequence
\[
0 \to \mathcal I_C/\mathcal I_C^2 \to \Omega^1_{\mathcal X}|_C \to \Omega^1_C \to 0.
\]   
Our claim follows by combining this with \cite[Chapter I.2 Theorem 2.15.2]{Kollar}.

When $C$ avoids the singular locus of $\mathcal X$ and $H^{1}(C,N_{C/\mathcal{X}}) = 0$, then the expected dimension coincides with the actual dimension and $C$ represents a smooth point of $M$.  The quantity $H^{0}(C,N_{C/\mathcal{X}})$ is an upper bound for the dimension of $M$ assuming that $C$ avoids the singularities of $\mathcal{X}$.

\subsubsection{Compactifications} Using the functor defining the Hilbert scheme, we see that each component $M \subset \Sec(\mathcal{X}/B)$ admits a natural embedding into the stack $\mathcal{M}_{g(B),0}(\mathcal{X})$.  Since every stable map parametrized by $M$ is birational onto its image these curves have trivial automorphism group over $\mathcal{X}$, so that $M$ also admits a natural embedding into a component $\overline{M}$ of the coarse moduli space $\overline{M}_{g(B),0}(\mathcal{X})$. 
Note that the domain of any stable map parametrized by $\overline{M}$ will be the union of a section with a finite collection of trees of $\pi$-vertical rational curves.

\subsubsection{Fixed points} \label{subsubsec:fixedpoints} We will also need parameter spaces for sections containing a set of fixed points.  When the base has genus $0$ we can appeal to \cite{Shen12} which constructs a moduli space for rational curves through fixed points.  We will briefly discuss the situation for curves of higher genus.

Let $q_{1},\ldots,q_{r}$ be $r$ points on $\mathcal{X}$ whose images $p_{i} := \pi(q_{i})$ are distinct.  Recall that \cite[Section I.1]{Kollar} constructs the moduli space $\Mor(B,\mathcal{X},p_{i} \mapsto q_{i})$ of morphisms $f: B \to \mathcal{X}$ sending $p_{i} \mapsto q_{i}$.  This space admits a map to $\Aut(B,p_{i})$.  Since the genus of $B$ is at least $1$ and $r \geq 1$, this automorphism group is finite.  Using the universal properties, we see that the component of $\Mor(B,\mathcal{X},p_{i} \mapsto q_{i})$ lying above the trivial automorphism is isomorphic to the sublocus $\Sec(\mathcal{X}/B,q_{1},\ldots,q_{r})$ in $\Sec(\mathcal{X}/B)$ parametrizing sections containing the points $\{ q_{i} \}$. 
It is clear from the construction that we have natural inclusions $\Sec(\mathcal{X}/B,q_{1},\ldots,q_{r}) \to \Sec(\mathcal{X}/B,q_{1},\ldots,q_{r-1})$. 
The expected dimension of $\Sec(\mathcal{X}/B,q_{1},\ldots,q_{r})$ is 
\begin{equation*}
-K_{\mathcal{X}/B} \cdot C + (\dim \mathcal X-1)(1-g(B)-r)
\end{equation*}
and every component has at least the expected dimension.  If $C$ avoids the singular locus of $\mathcal X$ and $H^{1}(C,N_{C/\mathcal{X}}(-q_{1}-\ldots-q_{r})) = 0$ then the expected dimension coincides with the actual dimension and $C$ represents a smooth point of $M$.  If $C$ is in the smooth locus of $\mathcal X$ the quantity $H^{0}(C,N_{C/\mathcal{X}}(-q_{1}-\ldots-q_{r}))$ is an upper bound for the dimension of $M$.

\subsection{Dominant families}

Suppose that $M \subset \Sec(\mathcal{X}/B)$ parametrizes a dominant family of sections.  Let $C$ be a general section parametrized by $M$; we will study the properties of the normal bundle.

\begin{defi}
We say that a section $C$ is relatively free if $C$ avoids the singular locus of $\mathcal X$, $H^{1}(C,N_{C/\mathcal{X}}) = 0$, and $N_{C/\mathcal{X}}$ is generically globally generated, in the sense that the evaluation map
\begin{equation*}
H^{0}(C,N_{C/\mathcal{X}}) \otimes \mathcal{O}_{C} \to N_{C/\mathcal{X}}
\end{equation*}
is surjective at the generic point of $C$.
\end{defi}

\begin{rema}
According to \cite[II.3 Definition 3.1]{Kollar} a curve $f: C \to X$ is free if $C$ avoids the singular locus of $X$, $H^{1}(C,f^{*}T_{X}) = 0$, and $f^{*}T_{X}$ is globally generated.  Note that our notion of relatively free is not quite a ``relative version'' of the notion of freeness, since we only require the normal bundle to be generically globally generated.  Our decision to use a weaker notion is motivated by Corollary \ref{coro:freecurvesmove}.
\end{rema}

The following proposition connects the existence of deformations of $C$ through general points with the generically globally generated condition of $N_{C/\mathcal{X}}$.

\begin{prop} \label{prop:deffixpoints} 
Let $\pi : \mathcal X \rightarrow B$ be a Fano fibration.  
Fix different points $q_{1},\ldots,q_{m}$ of $\mathcal X$ and let $M$ denote a component of $\Sec(\mathcal{X}/B,q_{1},\ldots,q_{m})$.  Suppose that the sections parametrized by $M$ dominate $\mathcal{X}$ and the general such section avoids the singular locus of $\mathcal X$.  Then for a general section $C$ parametrized by $M$ and for a general point $p \in B$ we have that $H^{0}(C, N_{C/\mathcal X}(-q_{1}-\ldots-q_{m})) \to N_{C/\mathcal X}|_{p}$ is surjective.

Conversely, suppose we fix a section $f: B \rightarrow \mathcal X$ whose image $C$ avoids the singularities of $\mathcal X$. 
%Let $p_{1},\ldots,p_{m}$ denote the points of $C$ lying over distinct points $q_{1},\ldots,q_{m}$ of $B$.   
Suppose that $q_{1},\ldots,q_{m}$ are distinct points of $C$ such that $H^1(C, N_{C/\mathcal X}(-q_{1}-\ldots-q_{m}))= 0$.  Let $M \subset \Sec(\mathcal{X}/B,q_{1},\ldots,q_{m})$ denote the unique component containing $f$.  If for a general point $p \in C$ we have that $H^{0}(C, N_{C/\mathcal X}(-q_{1}-\ldots-q_{m})) \to N_{C/\mathcal X}|_{p}$ is surjective, then $M$ parametrizes a dominant family of curves on $\mathcal{X}$.
\end{prop}

\begin{proof}
This is proved by \cite[Section 2]{Shen12} when the genus is $0$.  Let $p_1, \cdots, p_m$ be the images of $q_1, \cdots, q_m$ on $B$. In the situations where $\Sec(\mathcal{X}/B,q_{1},\ldots,q_{m})$ can be identified with a component of $\Mor(B,\mathcal{X},p_{i} \mapsto q_{i})$, this follows from the deformation theory for morphisms described by \cite[Section II.3]{Kollar}.  The only remaining case is when $g(B) = 1$ and $m = 0$, and it follows from similar arguments.
\end{proof}

In particular, we obtain an alternative description of relative free sections.

\begin{coro} \label{coro:freecurvesmove}
Let $C$ be a section that is general in its family.  Then $C$ is relatively free if and only if $C$ avoids the singularities of $\mathcal X$, $H^{1}(C,N_{C/\mathcal{X}}) = 0$, and $C$ is movable.
\end{coro}

We next give two statements relating the space of sections of $N_{C/\mathcal{X}}$ with the existence of deformations of $C$ through general points.

\begin{prop} \label{prop:pointnormalbundledim1} 
Let $\pi: \mathcal{X} \to B$ be a Fano fibration of relative dimension $1$.  Let $C$ be a section of $\pi$.   If $H^{1}(C,N_{C/\mathcal{X}}) = 0$ and for some positive integer $m$ we have $h^{0}(C,N_{C/\mathcal{X}}) \geq m$ then deformations of $C$ go through $m$ general points of $\mathcal{X}$.  

Conversely, if deformations of $C$ go through $m$ general points of $\mathcal{X}$ then
\begin{equation*}
h^{0}(C,N_{C/\mathcal{X}}) \geq m.
\end{equation*}
\end{prop}

\begin{prop} \label{prop:pointsandnormalbundle} 
Let $\pi: \mathcal{X} \to B$ be a del Pezzo fibration.
\begin{enumerate}
\item 
Let $C$ be a section of $\pi$ contained in the smooth locus of $\mathcal X$ whose normal bundle is unstable with the maximal destabilizing exact sequence 
\[
0 \to \mathcal L_2 \to N_{C/\mathcal{X}} \to \mathcal L_1 \to 0.  
\]

Suppose that $-K_{\mathcal X/B} \cdot C \geq4g(B)-4$ and $H^{1}(C,N_{C/\mathcal{X}}) = 0$.  Suppose $m$ is a positive integer such that
\begin{equation*}
h^{0}(C,\mathcal{L}_{1}) \geq m.
\end{equation*}
Then $C$ is relatively free and deformations of $C$ go through at least $m$ general points of $\mathcal{X}$.  

Conversely, suppose that $C$ is general in moduli.  Fix a positive integer $m$.  If $C$ goes through at least $m$ general points of $\mathcal{X}$ then we have
\begin{equation*}
h^{0}(C,\mathcal{L}_{1}) \geq m.
\end{equation*}

\item Let $C$ be a section of $\pi$ contained in the smooth locus of $\mathcal X$ whose normal bundle is semistable.  

Suppose that $-K_{\mathcal X/B} \cdot C \geq 6g(B)-2$ and $H^{1}(C,N_{C/\mathcal{X}}) = 0$.  Suppose $m$ is a positive integer such that
\begin{equation*}
h^{0}(C,N_{C/\mathcal{X}}) \geq 2m.
\end{equation*}
Then $C$ is relatively free and deformations of $C$ go through at least $m$ general points of $\mathcal{X}$.  

Conversely, fix a positive integer $m$.  If deformations of $C$ go through at least $m$ general points of $\mathcal{X}$ then we have
\begin{equation*}
h^{0}(C,N_{C/\mathcal{X}}) \geq 2m.
\end{equation*}

\end{enumerate}
\end{prop}

We will only prove Proposition \ref{prop:pointsandnormalbundle}, since the proof of Proposition \ref{prop:pointnormalbundledim1} is similar but easier.

\begin{proof}
Case (1): 
We start by proving the reverse implication.  Suppose that $C$ is general in moduli and goes through $m$ general points.  By Proposition \ref{prop:deffixpoints}, if $q_1, \cdots, q_m$ are general points on $C$ then for any $j=1,2,\ldots,m$ we have a surjection
\begin{equation*}
H^{0}(C,N_{C/\mathcal{X}}(-q_{1}-\ldots-q_{j-1})) \to N_{C/\mathcal{X}}(-q_{1}-\ldots-q_{j-1})|_{q_{j}}.
\end{equation*}
Using the exact sequence
\begin{equation*}
0 \to \mathcal L_{2}(-q_{1}-\ldots-q_{m-1})  \to N_{C/\mathcal{X}}(-q_{1}-\ldots-q_{m-1}) \to  \mathcal L_{1}(-q_{1}-\ldots-q_{m-1})\to 0
\end{equation*}
we see that we must have $h^{0}(C,\mathcal{L}_{1}(-q_{1}-\ldots-q_{m-1})) > 0$.  
Indeed, there are at least two sections of $H^{0}(C,N_{C/\mathcal{X}}(-q_{1}-\ldots-q_{m-1}))$ which generate $N_{C/\mathcal{X}}(-q_{1}-\ldots-q_{m-1})|_{q_{m}}$. Then one of these sections must generate $\mathcal L_1(-q_1 - \ldots -q_{m-1})|_{q_m}$. Thus our claim follows. 
Since the points $q_{i}$ were generic, 
we deduce that $h^{0}(C,\mathcal{L}_{1}) \geq m$.

Next we prove the forward implication.  Suppose that $-K_{\mathcal X/B} \cdot C \geq4g(B)-4$, $H^{1}(C,N_{C/\mathcal{X}}) = 0$, and for some positive integer $m$ we have
\begin{equation*}
h^{0}(C,\mathcal{L}_{1}) \geq m.
\end{equation*}
Since $h^{1}(C,N_{C/\mathcal{X}}) = 0$ we have $h^1(C, \mathcal L_1) = 0$.
Moreover since $\deg(\mathcal L_2) > -\frac{1}{2}K_{\mathcal X/B} \cdot C \geq 2g(B)-2$ we have $h^1(C, \mathcal L_2) = 0$.
Since we have $\deg(\mathcal L_2) > \deg(\mathcal L_1)$, we must have
\[
h^{0}(C,\mathcal{L}_{2}) \geq m.
\]
If we fix $m-1$ general points $q_{1},\ldots,q_{m-1}$ in $C$ then for $i=1,2$
\begin{equation*}
h^{0}(C,\mathcal{L}_{i}(-q_{1}-\ldots-q_{m-1})) = h^{0}(C,\mathcal{L}_{i}) - (m-1).
\end{equation*}
Thus $h^{1}(C,\mathcal{L}_{i}(-q_{1}-\ldots-q_{m-1})) = 0$ and $h^{0}(C,\mathcal{L}_{i}(-q_{1}-\ldots-q_{m-1})) \geq 1$.  Since for $i=1,2$ we have a surjection
\begin{equation*}
\mathcal L_i(-q_{1}-\ldots-q_{m-1}) \to \mathcal{L}_{i}(-q_{1}-\ldots-q_{m-1})|_{p}
\end{equation*}
for a general point $p$, the Snake Lemma shows that $N_{C/\mathcal{X}}(-q_{1}-\ldots-q_{m-1})$ surjects onto $N_{C/\mathcal{X}}(-q_{1}-\ldots-q_{m-1})|_{p}$ for a general point $p$.  Applying Proposition \ref{prop:deffixpoints} we obtain the desired statement.

Case (2): 
We first prove the reverse implication.  Let $C'$ be a general section in our family. Fix $m$ general points $q_{1},\ldots,q_{m}$ on $C'$.
By Proposition \ref{prop:deffixpoints}
\begin{equation*}
H^{0}(C',N_{C'/\mathcal{X}}(-q_{1}-\ldots-q_{j-1})) \to N_{C'/\mathcal{X}}(-q_{1}-\ldots-q_{j-1})|_{q_{j}}
\end{equation*}
is surjective for any $j=1,2,\ldots,m$.  Thus for $j$ in this range
\begin{equation*}
h^{0}(C',N_{C'/\mathcal{X}}(-q_{1}-\ldots-q_{j})) = h^{0}(C',N_{C'/\mathcal{X}}(-q_{1}-\ldots-q_{j-1})) - 2
\end{equation*}
and we conclude $h^{0}(C',N_{C'/\mathcal{X}}) \geq 2m$. Then our assertion follows from upper semicontinuity of $h^0$.

We next prove the forward implication.  Suppose that $h = -K_{\mathcal X/B} \cdot C \geq 6g(B)-2$
and $H^{1}(C,N_{C/\mathcal{X}}) = 0$.
We let $k \geq 0$ denote the maximal number of general points contained in deformations of $C$.  Assume for a contradiction that $h + 2- 2g(B) > 2k + 1$.  
We fix $k$ general points $q_1, \cdots, q_k$ on $C$ and set $D = q_{1} + \ldots + q_{k}$.  Using Proposition \ref{prop:deffixpoints} inductively as we twist down by general points, we see that $h^1(C, N_{C/\mathcal X}(-D)) = 0$.
Since $k$ is the maximal number of general points, Proposition \ref{prop:deffixpoints} shows that $N_{C/\mathcal X}(-D)$ is not generically globally generated.  Thus the image of the evaluation map
\begin{equation*}
H^0(C, N_{C/\mathcal X}(-D)) \otimes \mathcal{O}_{C} \to N_{C/\mathcal X}(-D)
\end{equation*}
is a rank $1$ subsheaf.  We let $\mathcal{L}_{D}$ denote the saturation of this subsheaf.  Note that by our assumption on $k$ we have
\[
h^0(C, \mathcal {L}_{D}) = h^0(C, N_{C/\mathcal X}(-D)) > 1.
\]
We next study how $\mathcal{L}_{D}$ varies as we change $D$.  Choose a different general point $q_{k+1}$ and set $D' = q_{1} + \ldots + q_{k-1} + q_{k+1}$.  Just as above we obtain a saturated subsheaf $\mathcal{L}_{D'}$ of $N_{C/\mathcal{X}}(-D')$.
Note that we have
\[
h^0(C, \mathcal{L}_{D}(-q_{k+1})) = h^0(C, N_{C/\mathcal X}(-q_1 - \cdots -q_{k +1})) = h^0(C, \mathcal{L}_{D'}(-q_{k}))
\]
and that all of these spaces of sections have dimension at least $1$.  This implies that both $\mathcal{L}_{D}(-q_{k+1})$ and $\mathcal{L}_{D'}(-q_{k})$ agree with the saturation of the subsheaf of $N_{C/\mathcal{X}}(-q_{1} - \ldots - q_{k+1})$ generated by global sections.  In particular $\mathcal{L}_{D'} = \mathcal{L}_{D}(q_{k}-q_{k+1})$.

Define $\mathcal{L} = \mathcal{L}_{D}(D)$.  Note that by our argument above $\mathcal{L}$ does not depend on the choice of the sum of $k$ general points.  Consider the exact sequence
\[
0 \to \mathcal L \to N_{C/\mathcal X} \to \mathcal K \to 0.
\]
Then $\mathcal K$ is invertible because $\mathcal L$ is saturated.
Since we have
\[
h^0(C, N_{C/\mathcal X}) = h^0(C, N_{C/\mathcal X}(-q_1 - \cdots - q_k)) + 2k,
\]
we must have
\begin{align*}
h^0(C, \mathcal L) = h^0(C, \mathcal L(-q_1 - \cdots - q_k)) + \ell\\
h^0(C, \mathcal K) = h^0(C, \mathcal K(-q_1 - \cdots - q_k)) + k,
\end{align*}
for some $\ell \leq k$.  Since the sections of $N_{C/\mathcal{X}}(-q_{1} - \cdots - q_{k})$ generate $\mathcal{L}(-q_{1} - \cdots - q_{k})$ we also must have 
\[
h^0(C, \mathcal K(-q_1 - \cdots - q_k)) = h^1(C, \mathcal L(-q_1 - \cdots - q_k)) = h^1(C, \mathcal L) + k-\ell.
\]
Since $h^0(C, \mathcal L(-q_1 - \cdots - q_k)) > 1$ we see that $g(B) \geq h^1(C, \mathcal L(-q_1 - \cdots - q_k)) = k - \ell +h^1(C, \mathcal L)$.  Thus
\begin{align*}
\deg(\mathcal{K}) & \leq h^{0}(C,\mathcal{K}) + g(B) - 1 \\
& \leq (g(B) + k) + g(B) - 1 \\
& \leq \ell + 3g(B) - 1 \\
& < h^{0}(C,\mathcal{L}) + 3g(B) - 1.
\end{align*}
If $h^0(C, \mathcal L) \leq g(B)$, then we conclude that $\deg(\mathcal K)< 4g(B)-1$.  This would also imply that $\deg(\mathcal{L}) \leq 2g(B) - 1$.
Combining we get $h < 6g(B)-2$, contradicting our height bound. 
Hence we must have $h^0(C, \mathcal L) \geq g(B)+1$ so that $h^1(C, \mathcal L) = 0$.  

Since $N_{C/\mathcal{X}}$ is semistable we see that $\deg(\mathcal L) \leq \deg(\mathcal K)$.  Since $h^1(C, \mathcal L) = h^1(C, \mathcal K) = 0$ we deduce that deformations of $C$ go through at least $\deg(\mathcal L) + 1-g(B)$ general points.  Set $s = \deg(\mathcal{L}) + 1 - g(B)$.  Note in particular that $k \geq s$.
For a set of $s$ general points $\{ q_i \}$ we have $h^0(C, \mathcal L(-q_1 - \cdots -q_s)) = 0$ which contradicts with $h^0(C, \mathcal L(-q_1 - \cdots - q_k)) > 0$.   
Thus our assertion follows.
\end{proof}

We will frequently use the following two useful corollaries which allow us to easily show that a curve must be relatively free.

\begin{coro} \label{coro:generalpointsandfreeness}
Let $\pi: \mathcal{X} \to B$ be a del Pezzo fibration and let $C$ be a section.
\begin{enumerate}
\item 
Suppose that $C$ avoids the singular locus of $\mathcal X$ and the normal bundle of $C$ is unstable with the maximal destabilizing exact sequence
\[
0 \to \mathcal L_2 \to N_{C/\mathcal X} \to \mathcal L_1 \to 0.
\]  
If we have
\begin{equation*}
h^{0}(C,\mathcal{L}_{1})  \geq g(B)+1
\end{equation*}
then $C$ is relatively free.
\item Suppose that $C$ avoids the singular locus of $\mathcal X$ and the normal bundle of $C$ is semistable.  If 
\begin{equation*}
h^{0}(C,N_{C/\mathcal{X}}) \geq 4g(B)+1
\end{equation*}
then $C$ is relatively free.
\item
In particular, if $C$ is a section containing $2g(B)+1$ general points of $\mathcal{X}$ then a general deformation of $C$ avoids the singular locus of $\mathcal{X}$ and is relatively free.
\end{enumerate}
\end{coro}

\begin{proof}
In case (1), since $h^{0}(C,\mathcal{L}_{1})>g(B)$ we must have $\deg(\mathcal L_1) \geq 2g(B)-1$.  Since $\deg(\mathcal{L}_{2}) > \deg (\mathcal{L}_{1})$ this implies that $-K_{\mathcal{X}/B} \cdot C \geq 4g(B) - 1$ and $H^{1}(C,N_{C/\mathcal{X}}) = 0$. 
Then the statement follows from Proposition \ref{prop:pointsandnormalbundle}.

In case (2), write $0 \to \mathcal{L}_{1} \to N_{C/\mathcal{X}} \to \mathcal{L}_{2} \to 0$ as in Lemma \ref{lemm:indecomrk2}.  We have
\begin{equation*}
h^{0}(C,\mathcal{L}_{1}) + h^{0}(C,\mathcal{L}_{2}) \geq h^{0}(C,N_{C/\mathcal{X}}) \geq 4g(B)+1.
\end{equation*}
If $h^{0}(C,\mathcal{L}_{1}) > g(B)$ then $h^{1}(C,\mathcal{L}_{1}) = 0$.  In this case $\deg(\mathcal{L}_{2}) \geq \deg(\mathcal{L}_{1}) \geq 2g(B) - 1$ and $H^{1}(C,N_{C/\mathcal{X}}) = 0$.  If $h^{0}(C,\mathcal{L}_{1}) \leq g(B)$ then $h^{0}(C,\mathcal{L}_{2}) \geq 3g(B)+1$.  However, in this case we have $\deg(\mathcal{L}_{2}) - \deg(\mathcal{L}_{1}) > g(B)$, a contradiction to Lemma \ref{lemm:indecomrk2}.

Thus we conclude that $H^{1}(C,N_{C/\mathcal{X}}) = 0$.  By Riemann-Roch this implies that $\deg(N_{C/\mathcal{X}}) \geq 6g(B) - 1$.  We conclude the relative freeness condition using Proposition \ref{prop:pointsandnormalbundle}.

Finally we prove (3). If a general deformation of $C$ avoids the singularities of $\mathcal X$ then relative freeness follows from Proposition~\ref{prop:pointsandnormalbundle} and (1), (2) above, so we only need to prove that a general deformation of $C$ is contained in the smooth locus.  Let $M$ denote the component of $\Sec(\mathcal X/B)$ containing $C$ and suppose that all sections parametrized by $M$ meet with the singular locus of $\mathcal X$.  Let $C'$ be a general member of $M$.
We choose a resolution $\beta: \mathcal Y \to \mathcal X$ and let $\widetilde{C}'$ be the strict transform of $C'$.  Since $\mathcal X$ has terminal singularities
\[
-K_{\mathcal Y/B} \cdot \widetilde{C}' < -K_{\mathcal X/B} \cdot C'.
\]
On the other hand since $\widetilde{C}'$ contains $2g(B)+1$ general points of $\mathcal Y$ it must be relatively free. This means that the dimension of the component containing $\widetilde{C}'$ is equal to $-K_{\mathcal Y/B} \cdot C + 2(1-g(B))$.  This contradicts with
\[
\dim M \geq -K_{\mathcal X/B} \cdot C + 2(1-g(B))
\]
and our assertion follows.
\end{proof}

The analogue of Corollary \ref{coro:generalpointsandfreeness} in the relative dimension $1$ case is:

\begin{coro} \label{coro:freenessonsurface}
Let $\pi: \mathcal{X} \to B$ be a Fano fibration of relative dimension $1$.  Let $C$ be a section of $\pi$.  If either
\begin{enumerate}
\item deformations of $C$ go through $g(B)+1$ general points of $\mathcal{X}$, or
\item $-K_{\mathcal{X}/B} \cdot C \geq 2g(B)$
\end{enumerate}
then a general deformation of $C$ is relatively free.
\end{coro}

Finally, we will need the following observation.

\begin{lemm} \label{lemm:expdimatpoints}
Let $\pi: \mathcal{X} \to B$ be a Fano fibration of relative dimension $1$.  Let $C$ be a section of $\pi$.  If $-K_{\mathcal{X}/B} \cdot C \geq 2g(B)$ then for any codimension $2$ subset $Z$ of $\mathcal{X}$ there is a deformation of $C$ that avoids $Z$.
\end{lemm}

\begin{lemm} \label{lemm:expdimatpoints2}
Let $\pi: \mathcal{X} \to B$ be a del Pezzo fibration.  Let $C$ be a relatively free section of $\pi$ that is general in its family (and thus avoids the singular locus of $\mathcal{X}$).  Suppose that either:
\begin{enumerate}
\item the normal bundle of $C$ is semistable and $-K_{\mathcal{X}/B} \cdot C \geq 5g(B)$, or
\item the normal bundle of $C$ is unstable and indecomposable and $-K_{\mathcal{X}/B} \cdot C \geq 6g(B)-2$, or
\item the normal bundle of $C$ is split and the two summands $\mathcal{L}_{1},\mathcal{L}_{2}$ satisfy $\deg(\mathcal{L}_{i}) \geq 2g(B)$ for $i=1,2$.
\end{enumerate}
Then for any codimension $2$ subset $Z$ of $\mathcal{X}$ there is a deformation of $C$ that avoids $Z$.  In particular, this statement holds for any section $C$ whose deformations go through $\geq 2g(B)+1$ general points of $\mathcal{X}$.
\end{lemm}

We will only prove Lemma \ref{lemm:expdimatpoints2}, since the proof of Lemma \ref{lemm:expdimatpoints} is similar but easier.

\begin{proof}
In each of the three circumstances above we have $h^{1}(C,N_{C/\mathcal{X}}(-p)) = 0$ for every point $p$ in $C$.  Thus the space $\Sec(\mathcal{X}/B,p)$ has the expected dimension, and the first statement follows easily.  The second statement follows from Corollary \ref{coro:generalpointsandfreeness} and Riemann-Roch.
\end{proof}

\subsection{Gluing criterion}

We will need the following lemma allowing us to glue relatively free sections to free vertical curves.
\begin{lemm}
\label{lemm: gluingandsmoothing}
Let $\pi: \mathcal{X} \to B$ be a Fano fibration.  Suppose that $f: C \to \mathcal{X}$ is a stable map satisfying the following conditions:
\begin{enumerate}
\item The domain of $C$ is a genus $g(B)$ curve which consists of a single curve $C_{0}$ isomorphic to $B$ attached to several trees of rational curves, and these trees of rational curves are contracted by $\pi \circ f$.
\item The image of $f$ is contained in the smooth locus of $\mathcal{X}$.
\item  The restriction of $f$ to $C_{0}$ is an isomorphism from $C_{0}$ to a relatively free section of $\pi$.
\item For each rational component $T_{i}$ of $C$, $f^*T_{\mathcal X}|_{T_i}$ is globally generated. 
\end{enumerate}
Then $f$ is a smooth point of $\overline{M}_{g(B),0}(\mathcal{X})$.
Furthermore assume that $f$ is an immersion on a neighborhood of every node on $C_0$. 
Then a general point of the corresponding component of $\overline{M}_{g(B),0}(\mathcal{X})$ parametrizes a relatively free section of $\pi$.
\end{lemm}

The proof is a combination of well-known techniques in deformation theory; see for example \cite[Section 1]{Testa09} and \cite[Lemma 2.6]{GHS03}.

\subsection{Shen's work on curves with unstable normal bundles}
\label{subsec:Shen}

Suppose that a section $C$ has unbalanced normal bundle.  If we look at all deformations of $C$ which contain a suitable number of general points then the resulting sections sweep out a surface $\Sigma$ in $\mathcal{X}$.  In \cite{Shen12} Shen studied the geometry of this surface $\Sigma$ in the setting of rational curves on $3$-folds with unbalanced normal bundles.  In this section, we will make a few modifications to develop a similar story for sections of higher genus.  The goal is Proposition \ref{prop:normalbundleinsurface} which we later use to produce rational curves in $\Sigma$ by applying Bend-and-Break to a large family of sections.

Suppose that we have a del Pezzo fibration $\pi : \mathcal X \rightarrow B$. Let $M\subset \Sec (\mathcal X/B)$ be a component generically parametrizing relatively free sections. 
Let $C$ be a general section parametrized by $M$.
We assume that $N_{C/\mathcal X}$ is unstable so that it fits into the maximal destabilizing exact sequence
\[
0 \to \mathcal L_2 \to N_{C/\mathcal X} \to \mathcal L_1 \to 0.
\]
We define $h = -K_{\mathcal X/B} \cdot C$ and assume that $h \geq 5g(B)-2$. 
As the general section parametrized by $M$ is relatively free, by definition we have $h^{1}(C,N_{C/\mathcal{X}}) = 0$.  This implies that $h^{1}(C,\mathcal{L}_{1}) = 0$, and since $\deg(\mathcal{L}_{2}) > \deg(N_{C/\mathcal{X}})/2 = h/2$ we also have $h^{1}(C, \mathcal L_2) = 0$.  Define $k = \deg (\mathcal L_1) + 1 -g(B)$.  Since $h^{1}(C,\mathcal{L}_{1}) = 0$ we must have have $k \geq 1$.

Assume that $\deg(\mathcal L_2) > \deg (\mathcal L_1) + 1$.   Fix $k$ general points $q_1, \ldots, q_k$ on $C$ so that $H^1(C, \mathcal L_i(-q_1- \ldots -q_k)) = 0$.  Then the deformation space $\Sec(\mathcal X/B; q_1, \cdots, q_k)$ is smooth at $C$ and at least $2$ dimensional.  Denote the unique component of $\Sec(\mathcal X/B; q_1, \cdots, q_k)$ that contains $C$ by $M(C; q_1, \cdots, q_k)$.  As discussed in Section \ref{subsubsec:fixedpoints} we have an inclusion $M(C; q_1, \cdots, q_k) \subset M$. Let
\begin{equation*}
\mathcal U(C; q_1, \cdots, q_k) \to M(C; q_1, \cdots, q_k)
\end{equation*}
be the universal family with the evaluation map $v : \mathcal U(C; q_1, \cdots, q_k) \to \mathcal X$. Since by Proposition \ref{prop:pointsandnormalbundle} deformations of $C$ can not go through more than $k$ general points of $\mathcal{X}$, the closure $\Sigma$ of the image $v(\mathcal U(C; q_1, \cdots, q_k))$ is a surface. Let $\Sigma' \to \Sigma$ be the normalization.  By the universal property of normalizations there is some neighborhood of the curve $C$ in $\mathcal U(C; q_1, \cdots, q_k)$ such that the evaluation map $v$ restricted to this neighborhood factors through $\Sigma'$.  Then we have

\begin{lemm}
Suppose $M$ parametrizes a relatively free section $C$ as above.  Then $\Sigma'$ is smooth along $C \setminus \{q_1, \cdots, q_k\}$.
\end{lemm}
\begin{proof}
We claim that $\mathcal L_2(-q_1 - \cdots - q_k)$ is basepoint free.  
 If $\deg(\mathcal L_2) - \deg (\mathcal L_1) + g(B) -1 \geq 2g(B)$, then $\deg (\mathcal L_2(-q_1 - \cdots - q_k)) \geq 2g(B)$ and thus this line bundle is basepoint free. 
 Otherwise $\deg (\mathcal L_2) - \deg (\mathcal L_1) \leq g(B)$. Since we have $2\deg (\mathcal L_1) +g(B) \geq \deg (\mathcal{L}_{1}) + \deg (\mathcal{L}_{2}) \geq 5g(B) - 2$, in this situation
\begin{equation*}
k = \deg (\mathcal L_1) -g(B) + 1   \geq g(B).
\end{equation*}
Thus the class $[\mathcal L_2(-q_1 - \cdots - q_k)]$ in $\Pic(B)$ is generic and has degree $\geq g(B)+1$ so that it is basepoint free.
Now our assertion follows from the argument of \cite[Lemma 2.4]{Shen12}.   
Indeed, the destabilizing exact sequence shows that for any $q \in C \setminus\{q_1, \cdots, q_k\}$ the map $H^0(C, N_{C/\mathcal X}) \to N_{C/\mathcal X}|_q$ has rank $1$.  We are thus in a situation where we can apply \cite[Lemma 2.3]{Shen12}.
\end{proof}

Then arguing as in \cite[Proposition 2.5]{Shen12}, we have
\begin{lemm}
\label{lemm:Shen_independent}
Suppose $M$ parametrizes a relatively free section $C$ as above.  Let $k = \deg (\mathcal{L}_{1}) - g(B)+1$ and choose general points $q_1, \cdots, q_k$ on $C$.
The surface $\Sigma = \Sigma(C; q_1, \cdots, q_k)$ is independent of $q_1, \cdots, q_k$.
\end{lemm}

\begin{proof}
We already showed that $\mathcal L_2(-q_1 - \cdots - q_k)$ is basepoint free. Hence for any $q \in C \setminus\{q_1, \cdots, q_k\}$ we have
\[
h^0(C, N_{C/\mathcal X}(-q_1 - \cdots - q_k-q)) = h^0(C, N_{C/\mathcal X}(-q_1 - \cdots - q_k))-1.
\]
Thus the arguments in \cite[Proposition 2.5]{Shen12} yield our assertion.
\end{proof}

Finally arguing as in \cite[Corollary 2.7]{Shen12}, we obtain 
\begin{prop}
\label{prop:normalbundleinsurface}
Suppose $M$ parametrizes a relatively free section $C$ as above.  Denote the maximal destabilizing exact sequence for the unstable normal bundle as
\[
0 \to \mathcal L_2 \to N_{C/\mathcal X} \to \mathcal L_1 \to 0.
\]
Let $k = \deg (\mathcal L_1) - g(B)+1$ and choose general points $q_1, \cdots, q_k$ on $C$.
Let $\Sigma = \Sigma(C; q_1, \cdots, q_k)$ and let $\Sigma' \to \Sigma$ be the normalization. Then $\Sigma'$ is smooth along the strict transform of $C$ and we have $N_{\Sigma'/\mathcal X}|_C = \mathcal L_1$ and $N_{C/\Sigma'} = \mathcal L_2$.
\end{prop}

\section{Breaking curves on surfaces}

The following conjecture is essential for understanding sections of Fano fibrations.

\begin{conj}[Movable Bend-and-Break for sections] \label{conj:movablebandb}
Let $\pi: \mathcal{X} \to B$ be a Fano fibration.  There is a constant $Q = Q(\mathcal{X})$ such that the following holds.  Suppose that $C$ is a relatively free section of $\pi$ satisfying $-K_{\mathcal{X}/B} \cdot C > Q(\mathcal{X})$.  Then $C$ deforms (as a stable map) to a union of a relatively free section with a $\pi$-vertical free curve.
\end{conj}

In this section we prove Conjecture \ref{conj:movablebandb} for sections of surfaces over $B$.  We first need a couple lemmas about the intersection theory of reducible fibers of a map from a surface to a curve.

\begin{lemm}[\cite{LTdPI} Lemma 4.3] \label{lemm:blowingup}
Let $Y$ be a smooth projective surface with a morphism $\pi : Y \rightarrow B$ such that a general fiber of $\pi$ is isomorphic to $\mathbb P^1$.  Let $F$ be a singular fiber of $\pi$ with components $\{ E_{i} \}_{i=1}^{r}$.  Suppose that $E_{1}$ is a $(-1)$-curve that has multiplicity $1$ in the fiber $F$.  Then there is another $(-1)$-curve in the fiber $F$.
\end{lemm}

\begin{coro} \label{coro:isolatingsection}
Let $Y$ be a smooth projective surface with a morphism $\pi : Y \rightarrow B$ such that a general fiber of $\pi$ is isomorphic to $\mathbb P^1$.  Let $C$ be a movable section of $\pi$.  Then there is a birational morphism $\phi: Y \to \mathbb{F}$ to a ruled surface $\mathbb{F}$ which is an isomorphism on a neighborhood of $C$.  This ruled surface satisfies
\begin{equation*}
K_{Y/B} \cdot C \leq \neg(\mathbb{F},-K_{\mathbb{F}/B}) \leq g(B).
\end{equation*}
(where $\neg$ is defined as in Definition \ref{defi:neg}.)
\end{coro}

\begin{proof}
Every reducible fiber of $\pi$ will carry a $(-1)$-curve.  Lemma \ref{lemm:blowingup} guarantees that if this $(-1)$-curve intersects $C$ then there is another $\pi$-vertical $(-1)$-curve in the same fiber that is disjoint from $C$.  Thus we may inductively contract $(-1)$-curves that are disjoint from $C$ to obtain $\mathbb{F}$.

To see the final statement, write $\rho_{*}C \equiv C'_{0} + kF$ where $C'_{0}$ is a section of minimal height and $F$ is a general fiber of the projective bundle.  Since $\rho_{*}C$ is movable we have $\rho_{*}C \cdot C'_{0} \geq 0$.  Thus by applying adjunction to $\rho_{*}C$ we see
\begin{equation*}
-K_{Y/B} \cdot C = -K_{\mathbb{F}/B} \cdot \rho_{*}C = (\rho_{*}C)^2 = \rho_{*}C \cdot (C'_{0} + kF) \geq k.
\end{equation*}
On the other hand
\begin{equation*}
0 \leq \rho_{*}C \cdot C'_{0} = C^{\prime 2}_{0} + k.
\end{equation*}
Since $-K_{\mathbb{F}/B} \cdot C'_{0} = C^{\prime 2}_{0}$ this yields the first inequality.  The second inequality holds true for every ruled surface and is a consequence of Lemma \ref{lemm:indecomrk2}.
\end{proof}

It will be helpful to have a numerical version of Corollary \ref{coro:isolatingsection}. 

\begin{defi} \label{defi:intersectionprofile}
Let $\pi: \mathcal{X} \to B$ be a Fano fibration such that $\mathcal X$ is smooth.  An intersection profile $\lambda$ for $\pi$ is a choice of a component $F_{0}$ in each fiber $F$ of $\pi$ such that $F_{0}$ has multiplicity $1$ in $F$.  We will denote the finite set of intersection profiles for $\pi$ by $\Lambda$.
\end{defi}

Note that any section $C$ naturally identifies an intersection profile for $\pi$ by selecting the components of fibers which meet $C$. 
 By repeating the arguments of Corollary \ref{coro:isolatingsection} we obtain:

\begin{lemm} \label{lemm:intersectionprofilebirationalmap}
Let $\pi: \mathcal{X} \to B$ be a Fano fibration of relative dimension $1$.  Fix an intersection profile $\lambda$ for $\pi$.  Then there is a birational morphism $\phi: \mathcal{X} \to \mathbb{F}$ to a ruled surface $\mathbb{F}$ over $B$ such that $\phi$ contracts every component of every fiber not identified by $\lambda$.
\end{lemm}

\begin{lemm}
\label{lemm:forBB}
Let $Y$ be a smooth projective surface with a morphism $\pi : Y \rightarrow B$ such that a general fiber of $\pi$ is isomorphic to $\mathbb P^1$.
Let $F$ be a reducible fiber of $\pi$ with components $\{ E_i \}_{i=1}^{r}$.
We assume that $E_1$ and $E_2$ have multiplicity $1$ in $F$. Suppose $Q = \sum_{i=1}^{r} a_i E_i$ is an effective $\mathbb{Q}$-divisor such that
\[
Q \cdot E_j =
\begin{cases}
1 & \text{ if } j = 1\\
-1 & \text{ if } j = 2\\
0 & \text{ otherwise}
\end{cases}
\]
Then we have
\[
-K_Y \cdot Q > 0.
\]
\end{lemm}

\begin{proof}
We prove this by induction on the number of components $r$ of $F$.
When $F$ consists of two components our assertion is trivial.

Suppose that $r$ is greater than $2$. If some component $E_j$ with $j \neq 1, 2$ is a $(-1)$-curve, then we can contract $E_j$ via $\phi: Y \to Y'$ and apply the induction hypothesis to deduce that
\begin{equation*}
-K_{Y} \cdot Q = -K_{Y'} \cdot \phi_{*}Q > 0.
\end{equation*}
So we may assume that there is no $(-1)$-curve other than $E_1$ and $E_2$.  Suppose that $E_1$ is a $(-1)$-curve.
By Lemma \ref{lemm:blowingup} $F$ must contain a different $(-1)$-curve, and by the argument above we may assume this is $E_{2}$.  Similarly, if $E_{2}$ is a $(-1)$-curve then we may assume that $E_{1}$ is also.  Thus we have reduced to the case when $E_{1}$ and $E_{2}$ are both $(-1)$-curves.

Since the multiplicity of $E_{2}$ in $F$ is $1$, there is a unique curve $E_3$ adjacent to $E_2$. If we contract $E_2$ and denote the blow down by $\phi : Y \rightarrow Y'$, then we have
\[
\phi_*Q \cdot \phi_*(E_j) = 
\begin{cases}
1 & \text{ if } j = 1\\
-1 & \text{ if } j = 3\\
0 & \text{ otherwise}
\end{cases}
\]
Then we have
\[
-K_Y \cdot Q = -K_{Y'} \cdot \phi_*Q + 1 >0
\]
by the induction hypothesis.
\end{proof}

The following theorem is the statement of Movable Bend-and-Break for sections of surface fibrations.

\begin{theo} \label{theo:bandbforsurfaces}
Let $\pi: Y \to B$ be a Fano fibration of relative dimension $1$.  Suppose that $C$ is a section satisfying
\[
-K_{Y/B} \cdot C \geq \max\{ 2, 2g(B)+1, 4g(B) + 1 -\neg(Y,-K_{Y/B})\}.
\]
Then there is a relatively free section $C_{0}$ and a general fiber $F$ such that
\[
C \sim_{alg} C_{0} + F.
\]
Furthermore $C$ deforms as a stable map to a morphism whose domain has two components where each component is mapped birationally onto either $C_{0}$ or $F$.
\end{theo}

\begin{proof}
Let $n = -K_{Y/B} \cdot C + 1 - g(B)$.  By our degree assumption we know that $H^{1}(C,N_{C/Y}) = 0$.  Since general points impose codimension $1$ conditions on the moduli space of sections, we see that deformations of $C$ contain $n$ general points of $Y$.  Thus there is a $1$-parameter family of curves through $n-1$ general points of $Y$.  By Lemma \ref{lemm:strongerbandb} $C$ deforms into the union of a section with some $\pi$-vertical curves, at least one of which goes through a general point.  Thus we can write
\[
C \sim_{alg} C_0 + mF + T
\]
where $F$ denotes a general fiber of $\pi$, $m \geq 1$, and $T$ is an effective $\pi$-vertical curve which does not deform.

Fix a fiber $F_{0}$ and let $T_{0}$ denote the sum of the components of $T$ contained in $F_{0}$.  By \cite[Lemma 1-2-10]{Matsuki02} there are two possibilities: either $T_{0}$ is proportional to a multiple of $F_{0}$ or $T_{0}$ has non-vanishing intersection against some components of $F_{0}$.  In the latter case $T_{0}$ must have exactly the intersection pattern as in the statement of Lemma \ref{lemm:forBB}.  By applying Lemma~\ref{lemm:forBB} and summing over all fibers, we conclude that $-K_{Y/B} \cdot T > 0$ unless $T = 0$.  We let $b$ denote this non-negative constant.

We divide the argument into two cases.  First suppose that $H^{1}(C_{0},N_{C_{0}/Y}) = 0$.
Since $-K_{Y/B} \cdot C_0 + 1 - g(B) = n-2m-b$, Proposition \ref{prop:pointnormalbundledim1} shows that $C_0$ can only go through $n-2m-b$ general points.  A fiber $F$ can go through only $1$ general point and $T$ cannot go through any general points.  But the broken curve $C_{0} + mF + T$ should pass through $n-1$ general points of $Y$.  Thus $b = 0$ and $m=1$ yielding the desired expression with $C' = C_{0}$.  We see that $C_{0}$ is relatively free using Corollary \ref{coro:freenessonsurface}.

Next suppose that $H^{1}(C_{0},N_{C_{0}/Y}) \neq 0$.
We will show that this case is impossible using our degree assumption.
The maximum number of general points that $C_{0}$ can contain is bounded above by $H^{0}(C_{0},N_{C_{0}/Y}) \leq g(B)$.  Thus we must have $m \geq n - g(B) - 1$ so that the entire curve can contain $n-1$ points.  But we also know that 
\begin{align*}
\neg(Y,-K_{Y/B}) & \leq -K_{Y/B} \cdot C_{0} \\
& = n - 1 + g(B) - 2m - b \\
& \leq 3g(B) + 1 - n \\
& = 4g(B) + K_{Y/B} \cdot C
\end{align*}
which contradicts our degree assumption.

To see the final statement, choose a deformation from $C$ to $C_{0} + F$.  The generic curve in this deformation family will be a section, yielding a map $\sigma: B \to Y$.  Taking a closure in the space of stable maps, we obtain a limit stable map whose image in $Y$ is $C_{0} + F$ which satisfies the desired properties.
\end{proof}

The following example illustrates that the intersection bound in Theorem \ref{theo:bandbforsurfaces} must depend on the choice of surface $Y$ and not just the genus of $B$.  In particular, this justifies the presence of $\neg(Y,-K_{Y/B})$ in the statement of Theorem \ref{theo:bandbforsurfaces}.

\begin{exam}
Let $\mathbb{F}_{e}$ denote the Hirzebruch surface whose rigid section has self-intersection $-e$ and let $\pi: \mathbb{F}_{e} \to \mathbb{P}^{1}$ denote the projective bundle structure.  Every movable section $C$ on $\mathbb{F}_{e}$ satisfies $-K_{\mathbb{F}_{e}/\mathbb{P}^{1}} \cdot C \geq e$.  Thus Theorem \ref{theo:bandbforsurfaces} can only apply to a section $C$ when the height of $C$ is $\geq e+1$.
\end{exam}

We also prove a breaking statement which deals with numerical equivalence.

\begin{prop} \label{prop:totalbreaking}
Let $\pi: Y \to B$ be a Fano fibration of relative dimension $1$. Let $C$ be a general member of a dominant family of sections on $Y$.  Let $\rho: Y \to \mathbb{F}$ be the birational map to a ruled surface $\mathbb{F}$ obtained by applying Corollary \ref{coro:isolatingsection} to $Y$ and $C$.  Then:
\begin{enumerate}
\item We have $C \equiv C_{0} + T + kF$ where $T$ is an effective $\pi$-vertical curve, $C_{0}$ is a section on $Y$ satisfying $-K_{Y/B} \cdot C_{0} \leq \neg(\mathbb{F},-K_{\mathbb{F}/B})$, and $k = \frac{1}{2}(-K_{Y/B} \cdot C) - \frac{1}{2}\neg(\mathbb{F},-K_{\mathbb{F}/B})$.
\item Define
\begin{equation*}
s = \frac{1}{2}(K_{Y/B} \cdot C) +\frac{1}{2}\neg(\mathbb{F},-K_{\mathbb{F}/B})+  \max \{4g(B) - 2,2g(B)  -\neg(\mathbb{F},-K_{\mathbb{F}/B})\} 
\end{equation*}
Then for any $r \geq 0$, there is a unique dominant family of relatively free sections $C_{1}$ on $Y$ which satisfy
\begin{equation*}
C_{1} \equiv C + \left(s + r \right)F
\end{equation*}
where $F$ denotes a general fiber of $\pi$.  For any $2r+1$ general points of $Y$ we can find a deformation of $C_{1}$ containing these $2r+1$ points.
\end{enumerate}
\end{prop}

\begin{proof}
Applying Corollary \ref{coro:isolatingsection} to $Y$ and $C$ we find a ruled surface $\mathbb{F}$ and a birational map $\rho: Y \to \mathbb{F}$ which is an isomorphism on a neighborhood of $C$.

Let $C'_{0}$ denote a section on $\mathbb{F}$ with minimal self intersection, so that $-K_{\mathbb F/B} \cdot C'_{0} = \neg(\mathbb{F},-K_{\mathbb{F}/B})$.  Then we have $\rho_{*}C \equiv C'_{0} + kF$ where $F$ denotes a general fiber of the projective bundle map and $k = \frac{1}{2}(-K_{Y/B} \cdot C) - \frac{1}{2}\neg(\mathbb{F},-K_{\mathbb{F}/B})$.  Write $\rho^{*}C'_{0} = C_{0} + T$.  Lemma \ref{lemm:forBB} implies that $K_{Y/B} \cdot T \leq 0$, proving (1).

To prove (2), we need to appeal to the geometry of ruled surfaces.  Define
\begin{equation*}
m = \max \{4g(B) - 2,2g(B)  -\neg(\mathbb{F},-K_{\mathbb{F}/B})\} 
\end{equation*}
so that $s = -k+m$. Suppose that $\mathcal{E}$ is a rank $2$ bundle such that $\mathbb{F} = \mathbb{P}_{B}(\mathcal{E})$.  After twisting $\mathcal{E}$, we may suppose that the section $C'_{0}$ constructed above is a section of $\mathcal{O}_{\mathbb{F}/B}(1)$.

First suppose that $\mathcal{E}$ is indecomposable.  Since $\mathcal{E}$ has a section, by Lemma \ref{lemm:indecomrk2} we see that $\mathcal{E}(D)$ is globally generated for any divisor $D$ on $B$ of degree at least $4g(B)-2$.  In particular, $\mathcal{E}(D)$ is globally generated for any $D$ of degree $m+r$ with $r \geq 0$.  Then sections of $\mathcal{E}(D)$ yield a basepoint free family of sections $C'_{1}$ on $\mathbb{F}$ with numerical class
\begin{equation*}
\rho_{*}C + (-k + m+r)F.
\end{equation*}
Furthermore, the additional twist by $rF$ guarantees that we can find deformations of $C'_{1}$ through any $2r+1$ general points of $\mathbb{F}$.  Since the same logic applies if we replace $D$ by any numerically equivalent divisor, we see that the component of $\Sec(\mathbb{F}/B)$ containing $C'_{1}$ is a projective bundle over $\Jac(B)$, and in particular, is irreducible.  To construct the desired family of sections on $Y$, note that when $C'_{1}$ and $F$ are general in their respective families they avoid all $\rho$-exceptional centers by Lemma~\ref{lemm:expdimatpoints}.  Thus by pulling back we obtain sections of the desired class on $Y$.  Note that pushforward and pullback by $\rho$ induce birational maps on the moduli space of sections of this class, so that sections on $Y$ of the given numerical class also must form an irreducible family.

Second suppose that $\mathcal{E}$ is decomposable.  Since $C'_{0}$ is a minimal section, we can write $\mathcal{E} = \mathcal{L}_{1} \oplus \mathcal{O}$.  Note that we have $\deg(\mathcal{L}_{1}) = \neg(\mathbb{F},-K_{\mathbb{F}/B})$.  Thus, $\mathcal{E}(D)$ is globally generated for any divisor $D$ on $B$ of degree $\geq 2g(B) - \neg(\mathbb{F},-K_{\mathbb{F}/B})$.  
We conclude by the same argument as before.
\end{proof}

\section{Non-dominant families of sections}

Let $\pi: \mathcal{X} \to B$ be a del Pezzo fibration such that $-K_{\mathcal{X}/B}$ is relatively nef.  Suppose that $M \subset \Sec(\mathcal{X}/B)$ is a component which defines a non-dominant family of sections.  The main goal of this section is to show that if the height of the sections parametrized by $M$ is sufficiently high then these sections sweep out a surface $Y$ with $a(Y_\eta, -K_{\mathcal X/B}|_{Y}) \geq a(X_\eta, -K_{\mathcal X/B})$.  Using this result we show there is a finite set of surfaces which contain all such families of sections.

We start with a couple results concerning surfaces swept out by sections.  The first lemma shows that if we have a family $N$ of surfaces $Y \subset \mathcal{X}$ satisfying $a(Y_{\eta},-K_{\mathcal{X}/B}|_{Y}) = a(\mathcal{X}_{\eta},-K_{\mathcal{X}/B})$ then we can stratify the parameter space $N$ according to the set of minimal models of $Y$. %\brian{I added here.}

\begin{lemm} \label{lemm:ruledsurfacedeform}
Let $\pi: \mathcal{X} \to B$ be a del Pezzo fibration. Suppose that we have a bounded family $\mathcal{S} \to N$ where each fiber is a surface $Y \subset \mathcal{X}$ such that $Y_{\eta}$ is geometrically irreducible and $a(Y_{\eta},-K_{\mathcal{X}/B}|_{Y}) = a(\mathcal{X}_{\eta},-K_{\mathcal{X}/B})$. Then there is a stratification $\{ N_{i}\}$ of $N$ into locally closed subsets with the following properties:
\begin{enumerate}
\item There is a family $\mathcal{S} \to N_{i}$ such that every fiber is a resolution of the corresponding surface $Y \subset \mathcal{X}$.
\item There is a base change $W_{i} \to N_{i}$ such that if we fix a resolution $\widetilde{Y}$ parametrized by $N_{i}$ and fix a birational map $\phi: \widetilde{Y} \to \mathbb{F}$ as in Lemma \ref{lemm:intersectionprofilebirationalmap} then the base change $\mathcal{S}_{W_{i}}$ admits a birational morphism over $W_{i}$ to a family where every fiber is isomorphic to a ruled surface $\mathbb{F}_{w} \to B$ and there is a point in $W_{i}$ representing $\widetilde{Y}$ such that the restriction of the map to this fiber is $\phi$.
\end{enumerate}
Furthermore, we may assume that the degrees of the terms in the Harder-Narasimhan filtration for the bundle $\mathcal{E}$ that defines $\mathbb{F}_{w}$ is constant for the entire family.
\end{lemm}

\begin{proof}
By repeatedly taking resolutions and restricting to the smooth locus we can stratify $N$ into a finite collection of irreducible locally closed subsets $N_{i}$ such that over each $N_{i}$ there is a smooth family of surfaces $\mathcal{S}_{i} \to N_{i}$ where the fiber over a point in $N_{i}$ is a  resolution of the corresponding point in $N$. 
The main step is to show that after repeatedly replacing $N_{i}$ by a non-empty open subset we may ensure the desired properties hold for $N_{i}$.  For notational clarity we will continue to call this open set $N_{i}$ even while making these changes.

Since every surface $Y$ satisfies $a(Y_{\eta},-K_{\mathcal{X}/B}|_{Y}) = a(\mathcal{X}_{\eta},-K_{\mathcal{X}/B})$, every resolution $\widetilde{Y}$ will admit a morphism to $B$ whose general fiber is birational to an irreducible $-K_{\mathcal{X}/B}$-conic. In particular, the Picard rank of $\widetilde{Y}$ is the same as the rank of $H^{2}(\widetilde{Y},\mathbb{Z})$.  Thus it is constant in each smooth family.  Moreover the N\'eron-Severi groups of $\widetilde{Y}$ form a local system over $N_i$.

Take a base change $W_{i} \to N_{i}$ that kills the monodromy action on the N\'eron-Severi groups of the fibers over $N_{i}$.  We let $\mathcal{S}_{W_{i}}$ denote the smooth family obtained by base change.  Note that the restriction map $N^{1}(\mathcal{S}_{W_{i}}) \to N^{1}(\widetilde{Y}_{w})$ is surjective for every surface $\widetilde{Y}_{w}$ in our family. 

For any fiber $\widetilde{Y}_{w}$ in our family consider the components of the reducible fibers of $\widetilde{Y}_{w} \to B$.  Using a relative Hilbert scheme argument, after perhaps shrinking $W_{i}$ (and $N_{i}$) we obtain a finite collection of irreducible divisors $D_{i}$ on $\mathcal{S}_{W_{i}}$ such that the intersection of the $D_{i}$ with the surfaces $\widetilde{Y}_{i}$ are exactly the components of the reducible fibers.  Due to the monodromy condition, after perhaps shrinking $W_{i}$ (and $N_{i}$) again we see that in fact the restriction of $D_{i}$ to each surface is irreducible and that when $i \neq j$ the restrictions of $D_{i}$ and $D_{j}$ to every fiber in our family are different components.

Using the constancy of intersection numbers, if the restriction of $D_{i}$ to some $\widetilde{Y}_{w}$ is a $(-1)$-curve then it is a $(-1)$-curve in every fiber.  By running the relative MMP for $\mathcal{S}_{W_{i}}$ over $W_{i} \times B$ we can contract the divisor $D_{i}$, and this contraction will result in the contraction of the corresponding $(-1)$-curve on each fiber.  This shows that any morphism $\phi: \widetilde{Y}_{w} \to \mathbb{F}$ obtained by contracting $(-1)$-curves in fibers can be spread out to the entire family over $W_{i}$.

Finally, by Theorem \ref{theo:hnfiltrationinfamilies} we can shrink $N_{i}$ so that for any of our families of ruled surfaces obtained as above the bundles $\mathcal{E}$ defining the ruled surfaces have the property that the terms in the Harder-Narasimhan filtration are of constant degree.  This guarantees that all the desired properties hold for our new $N_{i}$.  Repeating the argument on the complement and appealing to Noetherian induction we deduce the desired statement. 
\end{proof}

Next we show that if we have a family $N$ of surfaces $Y \subset \mathcal{X}$ satisfying $a(Y_{\eta},-K_{\mathcal{X}/B}|_{Y}) = a(\mathcal{X}_{\eta},-K_{\mathcal{X}/B})$ then we can deform a section in one surface $Y$ into other surfaces which lie in the same stratum as $Y$.

\begin{coro} \label{coro:curvescoveringfamily}
Let $\pi: \mathcal{X} \to B$ be a del Pezzo fibration. Suppose that we have a bounded family $\mathcal{S} \to N$ where each fiber is a surface $Y \subset \mathcal{X}$ such that $Y_{\eta}$ is geometrically irreducible and $a(Y_{\eta},-K_{\mathcal{X}/B}|_{Y}) = a(\mathcal{X}_{\eta},-K_{\mathcal{X}/B})$.  Let $\{N_{i}\}$ be the stratification of $N$ as in Lemma \ref{lemm:ruledsurfacedeform}.

Let $M \subset \Sec(\mathcal{X}/B)$ be a component.  Suppose that some sublocus of $M$ parametrizes a family of sections which sweep out a surface $Y$ parametrized by $N_{i}$.  Take the strict transform of these sections on the resolution $\widetilde{Y}$ identified by $N_{i}$ and consider the corresponding birational map $\rho: \widetilde{Y} \to \mathbb{F}$ as in Corollary \ref{coro:isolatingsection}.  Suppose that the sections $C_{\mathbb{F}}$ on $\mathbb{F}$ obtained in this way satisfy
\begin{equation*}
-K_{\mathbb{F}/B} \cdot C_{\mathbb{F}} \geq \neg(\mathbb{F},-K_{\mathbb{F}/B}) + 2 \max \{4g(B) - 2,2g(B)  -\neg(\mathbb{F},-K_{\mathbb{F}/B})\}.
\end{equation*}
Then for every surface $S$ parametrized by $N_{i}$ there is a unique component of $\Sec(\widetilde{S}/B)$ parametrizing sections of anticanonical degree $-K_{\mathbb{F}/B} \cdot C_{\mathbb{F}}$ which have intersection profile compatible with the birational map to $\mathbb{F}$ and the image of these sections in $\mathcal{X}$ is parametrized by $M$.  In particular every surface parametrized by $N_{i}$ is swept out by sections parametrized by $M$.
\end{coro}

\begin{proof}
Lemma \ref{lemm:ruledsurfacedeform} yields a base change $W_{i} \to N_{i}$ and a universal family $\mathcal{S}_{W_{i}} \to W_{i}$ such that $\mathcal{S}_{W_{i}}$ admits a birational map to a family of ruled surfaces $\mathcal{T}_{W_{i}}$ over $W_{i}$ whose restriction to $\widetilde{Y}$ coincides with $\rho$.  Due to our height restriction on $C_{\mathbb{F}}$, Proposition \ref{prop:totalbreaking} shows that every ruled surface in $\mathcal{T}_{W_{i}}$ admits a unique dominant family of sections which has the same anticanonical height as $C_{\mathbb{F}}$.  Using the description of these curves as sections of a twist of the bundle defining the ruled surface, we see that in fact all these sections form a single family on $\mathcal{T}_{W_{i}}$.  By taking the strict transform of this family on $\mathcal{S}_{W_{i}}$ and pushing forward to $\mathcal{X}$ we obtain the desired statement.
\end{proof}

The following theorem is the main result in this section.

\begin{theo} \label{theo:toomuchdeforming2}
Let $\pi: \mathcal{X} \to B$ be a del Pezzo fibration such that $-K_{\mathcal{X}/B}$ is relatively nef.  Let $M$ denote a component of $\Sec(\mathcal{X}/B)$ parametrizing sections $C$ satisfying
\begin{align*}
-K_{\mathcal{X}/B} \cdot C \geq  C(\mathcal X) := \max\{ &3g(B)+1, -2\neg(\mathcal{X},-K_{\mathcal{X}/B}) + 6g(B) - 2,  \\
 & 2g(B)-1 + 2 \max\{4g(B)-2, 3g(B)-1 - \neg(\mathcal X, -K_{\mathcal X/B})\}\}.
\end{align*}

Suppose that the closure of the locus swept out by the corresponding sections is a surface $Y$.  Then $a(Y_\eta, -K_{\mathcal X/B}|_{Y}) \geq a(X_\eta, -K_{\mathcal X/B})$.

If equality of $a$-invariants is achieved, then these sections satisfy the height bound of Corollary \ref{coro:curvescoveringfamily} with respect to the surface $Y$.   Moreover, $Y$ is swept out by a (possibly different) family of sections of height at most $2g(B) - 1 + 2\max \{4g(B) - 2,3g(B) - 1 - \neg(\mathcal{X},-K_{\mathcal{X}/B})\}$.
\end{theo}

\begin{proof}
Let $\phi: \widetilde{Y} \to Y$ denote a resolution of singularities.  Let $\widetilde{C}$ denote the strict transform of a general deformation of the section $C$.  By assumption the deformations of $\widetilde{C}$ are Zariski dense on $\widetilde{Y}$; thus the natural map $\psi: \widetilde{Y} \to B$ is an algebraic fiber space. Moreover our height bound guarantees that there is at least a $3$-parameter family of deformations of $C$ on $\mathcal{X}$, and hence also of $\widetilde{C}$ on $\widetilde{Y}$.  By Bend-and-Break we see that $\widetilde{Y}$ is generically a $\mathbb P^1$-bundle over the base.

Since by assumption $-K_{\mathcal{X}/B} \cdot C \geq 3g(B)+1$, the dimension of the space of deformations of $C$ is bounded below by $3g(B)+1 + 2 - 2g(B) = g(B)+3$.  Since the deformations of $C$ are contained in $Y$, we can find a deformation of $\widetilde{C}$ through $g(B)+3$ general points of $\widetilde{Y}$.  This implies that $H^{0}(\widetilde{C},N_{\widetilde{C}/\widetilde{Y}}) \geq g(B)+3$, and in particular by Corollary \ref{coro:freenessonsurface} $H^{1}(\widetilde{C},N_{\widetilde{C}/\widetilde{Y}}) = 0$.  Thus if $M$ denotes the the family of deformations of $\widetilde{C}$ on $\widetilde{Y}$, then
\begin{equation*}
\dim(M) = -K_{\widetilde{Y}/B} \cdot \widetilde{C} + (1-g(B)).
\end{equation*}
Since $M$ has at least the expected dimension on $\mathcal{X}$, we obtain
\begin{equation*}
(K_{\widetilde{Y}/B} - \phi^{*}K_{\mathcal{X}/B}|_{Y}) \cdot \widetilde{C}  \leq g(B)-1.
\end{equation*}
Apply Corollary \ref{coro:isolatingsection} to find a ruled surface $\mathbb{F}$ and a birational map $\rho: \widetilde{Y} \to \mathbb{F}$ which is an isomorphism on a neighborhood of $\widetilde{C}$.  As in Proposition \ref{prop:totalbreaking} (1) we can write $\widetilde{C} \equiv \widetilde{C}_{0} + kF + T$ where $F$ is a general fiber of the map to $B$ and
\begin{equation*}
k =  \frac{1}{2}(-K_{\widetilde{Y}/B} \cdot \widetilde{C}) - \frac{1}{2}\neg(\mathbb{F},-K_{\mathbb{F}/B}).
\end{equation*}
Then
\begin{align}
g(B) -1 & \geq K_{\widetilde{Y}/B} \cdot \widetilde{C} - \phi^{*}K_{\mathcal{X}/B}|_{Y} \cdot \widetilde{C} \nonumber \\
& =  K_{\widetilde{Y}/B} \cdot \widetilde{C} - \phi^{*}K_{\mathcal{X}/B}|_{Y} \cdot \widetilde{C}_{0} + k (-\phi^{*}K_{\mathcal{X}/B}|_{Y} \cdot F) -\phi^{*}K_{\mathcal{X}/B}|_{Y} \cdot T \nonumber \\
& \label{eqn:toomuchdeformingequation} \geq  K_{\widetilde{Y}/B} \cdot \widetilde{C} + \neg(\mathcal{X},-K_{\mathcal{X}/B}) + k (-\phi^{*}K_{\mathcal{X}/B}|_{Y} \cdot F)
\end{align}
where we have used the fact that $-K_{\mathcal{X}/B}$ is relatively nef at the last step.
Suppose for a contradiction that there is an inequality $a(Y_\eta, -K_{\mathcal X/B}|_{Y}) < a(X_\eta, -K_{\mathcal X/B})$.  This is equivalent to saying that $-\phi^{*}K_{\mathcal{X}/B}|_{Y} \cdot F \geq 3$.  Rearranging and substituting in our value for $k$, we find
\begin{align*}
\frac{1}{2}(-K_{\widetilde{Y}/B} \cdot \widetilde{C}) \leq  \frac{3}{2}\neg(\mathbb{F},-K_{\mathbb{F}/B}) -\neg(\mathcal{X},-K_{\mathcal{X}/B}) + g(B) - 1.
\end{align*}
Using the inequality $\neg(\mathbb{F},-K_{\mathbb{F}/B}) \leq g(B)$ from Corollary \ref{coro:isolatingsection}, we see that
\begin{align*}
-K_{\mathcal{X}/B} \cdot C & \leq -K_{\widetilde{Y}/B} \cdot \widetilde{C} + g(B) - 1 \\
& \leq -2\neg(\mathcal{X},-K_{\mathcal{X}/B}) + 6g(B) - 3.
\end{align*}
contradicting our height bound.

Suppose now that the equality of $a$-invariants is achieved so that $-\phi^{*}K_{\mathcal{X}/B}|_{Y} \cdot F = 2$.  Using this to simplify Equation \eqref{eqn:toomuchdeformingequation} we see that $-\neg(\mathcal{X},-K_{\mathcal{X}/B}) +g(B) - 1 \geq -\neg(\mathbb{F},-K_{\mathbb{F}/B})$.  Furthermore we know $\neg(\mathbb{F},-K_{\mathbb{F}/B}) \leq g(B)$ by Corollary \ref{coro:isolatingsection}.  Thus
\begin{align}
-K_{\widetilde{Y}/B} \cdot \widetilde{C} & \geq -\phi^{*}K_{\mathcal{X}/B}|_{Y} \cdot \widetilde{C} + 1 - g(B)  \nonumber \\
& \geq g(B) + 2 \max\{4g(B)-2, 3g(B)-1 - \neg(\mathcal X, -K_{\mathcal X/B}) \} \nonumber \\
& \label{eqn:toomuchdeformingequation5} \geq \neg(\mathbb{F},-K_{\mathbb{F}/B}) + 2 \max \{4g(B) - 2,2g(B)  -\neg(\mathbb{F},-K_{\mathbb{F}/B})\}
\end{align}
proving the first claim.  To see the last claim we will appeal to Proposition \ref{prop:totalbreaking} (2) with $r=0$.  This guarantees the existence of a dominant family of sections $\widetilde{C}_{1}$ on $\widetilde{Y}$ which satisfy $\widetilde{C} \equiv \widetilde{C}_{1} - sF$ where
\begin{equation*}
s = \frac{1}{2}(K_{Y/B} \cdot C) +\frac{1}{2}\neg(\mathbb{F},-K_{\mathbb{F}/B}) + \max \{4g(B) - 2,2g(B)  -\neg(\mathbb{F},-K_{\mathbb{F}/B})\} 
\end{equation*}
and $F$ denotes a general fiber of $\pi$.  Since $\widetilde{C}$ and $\widetilde{C}_{1}$ differ by a multiple of $F$, the equality of $a$-invariants shows that
\begin{align*}
g(B) -1 & \geq K_{\widetilde{Y}/B} \cdot \widetilde{C} - \phi^{*}K_{\mathcal{X}/B}|_{Y} \cdot \widetilde{C} \\
& = K_{\widetilde{Y}/B} \cdot \widetilde{C}_{1} - \phi^{*}K_{\mathcal{X}/B}|_{Y} \cdot \widetilde{C}_{1}.
\end{align*}
Thus
\begin{align*}
- \phi^{*}K_{\mathcal{X}/B}|_{Y} \cdot \widetilde{C}_{1} & \leq g(B) - 1  -K_{\widetilde{Y}/B} \cdot \widetilde{C}_{1} \\
& =  g(B) -1 + \neg(\mathbb{F},-K_{\mathbb{F}/B}) + 2\max \{4g(B) - 2,2g(B)  -\neg(\mathbb{F},-K_{\mathbb{F}/B})\} \\
& \leq 2g(B) - 1 + 2 \max\{4g(B)-2, 3g(B)-1 - \neg(\mathcal X, -K_{\mathcal X/B}) \} 
\end{align*}
where the last step is given by the bound $\neg(\mathbb{F},-K_{\mathbb{F}/B}) \leq g(B)$ from Corollary \ref{coro:isolatingsection}.  
Altogether this shows that $Y$ is swept out by sections satisfying the desired height bound.
\end{proof}

\begin{coro} \label{coro:nondomfambounded}
Let $\pi: \mathcal{X} \to B$ be a del Pezzo fibration such that $-K_{\mathcal X/B}$ is relatively nef. Then there is a proper closed subset $V \subsetneq \mathcal{X}$ such that any component $M \subset \Sec(\mathcal{X}/B)$ parametrizing a non-dominant family of sections will parametrize sections contained in $V$.
\end{coro}

\begin{proof}
Rigid sections are bounded in height and thus are contained in a closed subset.  Thus it suffices to consider sections which sweep out a surface.
By Theorem \ref{theo:toomuchdeforming2}, any component $M \subset \Sec(\mathcal{X}/B)$ parametrizing such a non-dominant family will satisfy one of the following properties: 
\begin{enumerate}
\item $M$ will parametrize sections $C$ satisfying $-K_{\mathcal{X}/B} \cdot C <  C(\mathcal X)$, or
\item $M$ will parametrize sections $C$ satisfying $-K_{\mathcal{X}/B} \cdot C \geq  C(\mathcal X)$ which sweep out a $2$-dimensional subvariety $Y$ satisfying $$a(Y_{\eta},-K_{\mathcal{X}/B}|_{Y}) > a(\mathcal{X}_{\eta},-K_{\mathcal{X}/B}),$$ or
\item $M$ will parametrize sections $C$ satisfying $-K_{\mathcal{X}/B} \cdot C \geq  C(\mathcal X)$  which sweep out a $2$-dimensional subvariety $Y$ satisfying $$a(Y_{\eta},-K_{\mathcal{X}/B}|_{Y}) = a(\mathcal{X}_{\eta},-K_{\mathcal{X}/B})$$ and which contains a section of height at most $2g(B) - 1 + 2\max \{4g(B) - 2,3g(B) - 1 - \neg(\mathcal{X},-K_{\mathcal{X}/B})\}$.
\end{enumerate}
Lemma \ref{lemm:northcott} shows that curves of the first type lie in a bounded family.  Thus, the union of the subvarieties swept out by the non-dominant families satisfying this height bound will be a proper closed subset of $\mathcal X$.  Corollary \ref{coro:closedainv} shows that the surfaces $Y$ defined by the components of the second type will lie in a proper closed subset of $\mathcal X$.  Proposition \ref{prop:boundednessequalainv} shows that the surfaces defined by the components of the third type form a bounded family $\mathcal{S} \to N$. We stratify the parameter space $N$ using Lemma \ref{lemm:ruledsurfacedeform}.  Let $N_{i}$ denote a component of this stratification.  Combining Corollary \ref{coro:curvescoveringfamily} and Theorem~\ref{theo:toomuchdeforming2}, we see that for any $d \geq C(\mathcal{X})$ a family of non-dominant sections of height $d$ which sweeps out a surface parametrized by $N_{i}$ will in fact sweep out every surface parametrized by $N_{i}$.  Thus if $\dim(N_{i}) \geq 1$ then we would have a dominant family, a contradiction.  This shows that there are only finitely many surfaces parametrized by our bounded family which can contain sections as in (3), proving our claim.
\end{proof}

\section{Dominant families which are not relatively free}

Let $\pi: \mathcal{X} \to B$ be a del Pezzo fibration.  If $B$ has genus $\geq 1$ then a dominant family of sections need not be relatively free.  In fact, such families can have arbitrarily large height and can have larger than the expected dimension.

\begin{exam}
Let $B$ be a curve of genus $\geq 1$.  For $d$ sufficiently large, there is a unique component $M \subset \Mor(B,\mathbb{P}^{1})$ parametrizing degree $d$-maps whose dimension is the expected value $2d + (1-g(B))$.

Give $\mathcal{X} = \mathbb{P}^{1} \times \mathbb{P}^{1} \times B$ the structure of a del Pezzo fibration using the projection $\pi: \mathcal{X} \to B$.  Consider the component $M \subset \Sec(\mathcal{X}/B)$ parametrizing sections such that the induced map $f: B \to \mathbb{P}^{1} \times \mathbb{P}^{1}$ is a degree $d$ cover of a fiber of the first projection.  The dimension of $M$ is $2d + (1-g(B)) + 1$.  This is always greater than the expected dimension $2d + 2(1-g(B))$.  In this situation we have $N_{C/\mathcal{X}} = \mathcal{O}_{C} \oplus f^{*}\mathcal{O}(2)$ and the discrepancy between the actual and expected dimension is accounted for by $H^{1}(C,N_{C/\mathcal{X}})$.
\end{exam}

In this section we will show that the existence and behavior of such families is controlled by the generic Fujita invariant.  We have two main goals.  First, we show that there is a bounded family of surfaces $Y$ with $a(Y_{\eta},-K_{\mathcal{X}/B}) = a(\mathcal{X}_{\eta},-K_{\mathcal{X}/B})$ which contains any section of this type of sufficiently large degree (Theorem \ref{theo: non-relativelyfreecurve}).  Second, we bound the dimension of such families so that we can control their contributions to the counting function (Lemma \ref{lemm:dominantexpecteddim}).

We start with a proposition that identifies the geometric properties of these dominant families under the assumption that the general curve does not intersect the singular locus of $\mathcal{X}$.

\begin{prop}
\label{prop:dominantnonrelativelyfree}
Let $\pi: \mathcal{X} \to B$ be a del Pezzo fibration.  Let $M$ denote a component of $\Sec(\mathcal{X}/B)$ parametrizing a dominant family of sections such that a general member $C$ satisfies
\begin{align*}
-K_{\mathcal{X}/B} \cdot C \geq 6g(B)-2.
\end{align*}
If the general section parametrized by $M$ is not relatively free and is contained in the smooth locus of $\mathcal X$, then:
\begin{enumerate}
\item The normal bundle of $C$ is unstable and split.
\item Let $q$ denote the maximal number of general points of $\mathcal{X}$ contained in a curve $C$ parametrized by $M$. Then $q \leq g(B)$.
\end{enumerate}
\end{prop}

\begin{proof}
Note that this situation can only happen if $g(B) > 0$.  Let $C$ be a general section in our family.  We split into two cases.

\textbf{Case 1:} $N_{C/\mathcal{X}}$ is semistable.  Using our height bound, we see the line bundles $\mathcal{L}_{1}$, $\mathcal{L}_{2}$ as in Lemma \ref{lemm:indecomrk2} have degree $\geq 2g(B)-1$.  Thus $h^{1}(C,N_{C/\mathcal{X}}) = 0$.  Since we know that $g(B) > 0$ we also see that $N_{C/\mathcal{X}}$ is generically globally generated.  Thus $C$ is relatively free, a contradiction.

\textbf{Case 2:} $N_{C/\mathcal{X}}$ is unstable.  Consider the maximal destabilizing exact sequence
\[
0 \to \mathcal L_2 \to N_{C/\mathcal X} \to \mathcal L_1 \to 0
\]
If $N_{C/\mathcal{X}}$ is indecomposable, then Lemma \ref{lemm:indecomrk2} shows that $\mathcal{L}_{1}$ and $\mathcal{L}_{2}$ have degree $\geq 2g(B)$.  In this case $C$ must be relatively free, a contradiction, finishing the proof of (1).

Since we know that the normal bundle of $C$ is split and unstable, we can write $N_{C/\mathcal{X}} = \mathcal{L}_{1} \oplus \mathcal{L}_{2}$ with $\deg(\mathcal{L}_{1}) < \deg(\mathcal{L}_{2})$.  Suppose that deformations of $C$ go through the maximum number $q$ of general points of $\mathcal{X}$.  If $q > g(B)$, then Proposition \ref{prop:pointsandnormalbundle} shows $h^{0}(C,\mathcal{L}_{1}) > g(B)$.  This implies that $\deg(\mathcal{L}_{1}) > 2g(B)-2$, and thus also $\deg(\mathcal{L}_{2}) > 2g(B)-2$.  We deduce that $h^{1}(C,N_{C/\mathcal{X}}) = 0$ and that the family is relatively free, a contradiction.  This proves (2). 
\end{proof}

Next we would like to extend the above proposition to sections which meet the singular locus of $\mathcal X$. For this purpose we introduce the following invariant:

\begin{defi}
Let $\pi : \mathcal X \to B$ be a del Pezzo fibration. Let $\beta: \mathcal Y \to \mathcal X$ be a resolution of singularities such that $\beta : Y_\eta \to X_\eta$ is an isomorphism.  We define
\[
m(\mathcal Y / \mathcal X) : = \max \{K_{\mathcal Y/\mathcal X} \cdot C \, | \, C \in \Sec(\mathcal Y/B)\}.
\]
This is a non-negative integer as $K_{\mathcal Y/\mathcal X}$ is a $\pi$-vertical effective divisor and any section $C$ meets with a vertical component with multiplicity at most $1$. Then we define
\[
m(\mathcal X/B) : = \min \{m(\mathcal Y / \mathcal X) | \beta : \mathcal Y \to \mathcal X \textrm{ a resolution}\}.
\]
\end{defi}

\begin{coro}
\label{coro:themaximumnumber}
Let $\pi: \mathcal{X} \to B$ be a del Pezzo fibration.  Let $M$ denote a component of $\Sec(\mathcal{X}/B)$ parametrizing a dominant family of sections $C$ such that
\begin{align*}
-K_{\mathcal{X}/B} \cdot C \geq 6g(B)-2 + m(\mathcal X/B).
\end{align*}
Suppose that the generic section parametrized by $M$ is not relatively free.  If we let $q$ denote the maximal number of general points of $\mathcal{X}$ contained in a curve $C$ parametrized by $M$ then $q \leq g(B)$.
\end{coro}

\begin{proof}
Let $\beta : \mathcal Y \to \mathcal X$ be a resolution achieving $m(\mathcal X/B)$.
Let $\widetilde{C}$ be the strict transform of a general $C \in M$.
If a general member of $M$ avoids the singular locus, then our assertion follows from Proposition~\ref{prop:dominantnonrelativelyfree}.
If all members of $M$ pass through the singular locus of $\mathcal X$, then $\widetilde{C}$ cannot be relatively free on $\mathcal Y$. Indeed if it is relatively free, then the dimension of $M$ is equal to $-K_{\mathcal Y/B} + 2(1-g(B))$. However since $\mathcal{X}$ has terminal singularities $-K_{\mathcal Y/B} \cdot \widetilde{C} < -K_{\mathcal X/B} \cdot C$ which is a contradiction with $\dim M \geq -K_{\mathcal X/B} + 2(1-g(B))$.

Now our height bound implies that
\[
-K_{\mathcal Y/B} \cdot \widetilde{C} \geq 6g(B)-2.
\]
Thus our assertion follows from Proposition~\ref{prop:dominantnonrelativelyfree}.
\end{proof}

This yields an estimate on the dimension of the components of $\Sec(\mathcal{X}/B)$ parametrizing curves of this type. 

\begin{lemm} \label{lemm:dominantexpecteddim}
Let $\pi: \mathcal{X} \to B$ be a del Pezzo fibration such that $-K_{\mathcal{X}/B}$ is relatively nef.  Suppose that $M \subset \Sec(\mathcal{X}/B)$ is a component that parametrizes a dominant family of non-relatively free sections $C$ with
\begin{align*}
-K_{\mathcal{X}/B} \cdot C \geq 6g(B) - 2 + m(\mathcal{X}/B).
\end{align*}
Then $\dim(M) \leq -K_{\mathcal{X}/B} \cdot C + 2 - g(B) + m(\mathcal{X}/B)$.
\end{lemm}

\begin{proof}
Choose a resolution $\beta: \mathcal{Y} \to \mathcal{X}$ achieving the minimum in the definition of $m(\mathcal{X}/B)$.  Let $C'$ be the strict transform of a general member of our family of sections on $\mathcal{Y}$.  Arguing as in the proof of Corollary \ref{coro:themaximumnumber} we see that the parameter space $M'$ of deformations of $C'$ has the same dimension as $M$ and that the sections $C'$ have height $\geq 6g(B) - 2$.

Since $\mathcal{Y}$ is smooth we have $\dim(M') \leq H^{0}(C',N_{C'/\mathcal{Y}})$.  By Proposition \ref{prop:dominantnonrelativelyfree} $N_{C'/\mathcal{Y}}$ splits.  Proposition \ref{prop:deffixpoints} shows that $N_{C'/\mathcal{Y}}$ is generically globally generated so that both summands must have non-negative degree.  Our height bound implies that at least one of the summands will have degree $\geq 2g(B) - 1$.  Thus $H^{1}(C',N_{C'/\mathcal{Y}}) \leq g(B)$.  By Riemann-Roch we obtain $H^{0}(C',N_{C'/\mathcal{Y}}) \leq -K_{\mathcal{Y}/B} \cdot C' + 2 - g(B)$, which proves the statement.
\end{proof}

The next theorem is our main result describing dominant families of sections that fail to be relatively free.

\begin{theo} \label{theo: non-relativelyfreecurve}
Let $\pi: \mathcal{X} \to B$ be a del Pezzo fibration such that $-K_{\mathcal{X}/B}$ is relatively nef.  Let $M$ denote a component of $\Sec(\mathcal{X}/B)$ parametrizing a dominant family of sections such that for a general $C \in M$ we have 
\begin{align*}
-K_{\mathcal{X}/B} \cdot C \geq \max \{& 6g(B)-2 + m(\mathcal X/B), -2\neg(\mathcal{X},-K_{\mathcal{X}/B}) + 12g(B) - 2, \\
&   8g(B) + 2\max \{4g(B) - 2,5g(B) - 1 - \neg(\mathcal{X},-K_{\mathcal{X}/B}) \}  \},
\end{align*}
and the general section $C$ is not relatively free. 

 Let $q$ denote the maximal number of general points contained in deformations of $C$.  Fix $q$ general points of $\mathcal{X}$ and fix a component $N$ of the sublocus of $M$ parametrizing sections through these $q$ points.  Then the sections parametrized by $N$ will sweep out a surface $Y$ that satisfies $a(Y_\eta, -K_{\mathcal X/B}|_{Y}) = a(\mathcal{X}_\eta, -K_{\mathcal X/B})$.  These sections will satisfy the height bound of Corollary \ref{coro:curvescoveringfamily} with respect to $Y$.

Furthermore, there is a component $M_{1}$ of $\Sec(\mathcal{X}/B)$ which satisfies the following conditions:
\begin{enumerate}
\item The sections $C_{1}$ parametrized by $M_{1}$ have height
\begin{equation*}
-K_{\mathcal{X}/B} \cdot C_{1} \leq 8g(B) - 1 +  2\max \{4g(B) - 2,5g(B) - 1 - \neg(\mathcal{X},-K_{\mathcal{X}/B}) \} .
\end{equation*}
\item The general surface $Y$ constructed above is swept out by sections parametrized by a sublocus of $M_{1}$.
\item There is an integer $r \geq 1$ such that the closure of $M$ in $\overline{M}_{g(B),0}(\mathcal{X})$ contains the union of a section parametrized by $M_{1}$ with $r$ general fibers of $Y \to B$.
\end{enumerate}
\end{theo}

\begin{proof}
First note that we have $g(B) \geq 1$ since any dominant family of rational sections generically parameterizes a free section.
By Corollary \ref{coro:themaximumnumber} we have $1 \leq q \leq g(B)$.

As in the statement of the theorem, let $N$ be any component of the sublocus of $M$ parametrizing curves through a fixed set of $q$ general points.  A dimension count shows that the sections parametrized by $N$ will sweep out a surface $Y$.  Let $\phi: \widetilde{Y} \to Y$ be a resolution and let $\widetilde{C}$ denote the strict transform of a general deformation of $C$ in $Y$.  By our height bound on $C$, we know there is at least a $g(B)$ dimensional family of deformations of $\widetilde{C}$ contained in $\widetilde{Y}$. Since this situation can only happen when $g(B) \geq 1$, by Bend-and-Break we deduce that the induced map $\psi: \widetilde{Y} \to B$ is generically a $\mathbb{P}^{1}$ fibration.  We also see that $h^{0}(\widetilde{C},N_{\widetilde{C}/\widetilde{Y}}) > g(B)$. By Kodaira vanishing the dimension of the space of sections of any line bundle of degree $2g-1$ on $\widetilde{C}$ is equal to $g(B)$.  Thus the normal bundle has degree $\geq 2g$, and in particular, $h^{1}(\widetilde{C},N_{\widetilde{C}/\widetilde{Y}}) = 0$.  We conclude that the family of deformations of $\widetilde{C}$ in $\widetilde{Y}$ has the expected dimension on $\widetilde{Y}$.  Since the family also has at least the expected dimension on $\widetilde{X}$, we deduce
\begin{equation*}
-K_{\widetilde{Y}/B} \cdot \widetilde{C} + (1-g(B)) \geq  - \phi^{*}K_{\mathcal{X}/B}|_{Y} \cdot \widetilde{C} + 2(1-g(B)) - 2q.
\end{equation*}

Apply Corollary \ref{coro:isolatingsection} to find a ruled surface $\mathbb{F}$ and a birational map $\rho: \widetilde{Y} \to \mathbb{F}$ which is an isomorphism on a neighborhood of $\widetilde{C}$.  As in Proposition \ref{prop:totalbreaking} (1) we can write $\widetilde{C} \equiv \widetilde{C}_{0} + kF + T$ where $F$ is a general fiber of the map to $B$ and
\begin{equation*}
k =  \frac{1}{2}(-K_{\widetilde{Y}/B} \cdot \widetilde{C}) - \frac{1}{2}\neg(\mathbb{F},-K_{\mathbb{F}/B}).
\end{equation*}
Then
\begin{align}
3g(B) -1 & \geq 2q + g(B) - 1  \nonumber \\
& \geq K_{\widetilde{Y}/B} \cdot \widetilde{C} - \phi^{*}K_{\mathcal{X}/B}|_{Y} \cdot \widetilde{C} \nonumber \\
& =  K_{\widetilde{Y}/B} \cdot \widetilde{C} - \phi^{*}K_{\mathcal{X}/B}|_{Y} \cdot \widetilde{C}_{0} + k (-\phi^{*}K_{\mathcal{X}/B}|_{Y} \cdot F) -\phi^{*}K_{\mathcal{X}/B}|_{Y} \cdot T \nonumber \\
& \label{eqn:toomuchdeformingequation2} \geq  K_{\widetilde{Y}/B} \cdot \widetilde{C} + \neg(\mathcal{X},-K_{\mathcal{X}/B}) + k (-\phi^{*}K_{\mathcal{X}/B}|_{Y} \cdot F)
\end{align}
where we have used the fact that $-K_{\mathcal{X}/B}$ is relatively nef at the last step.
Suppose for a contradiction that there is an inequality $a(Y_\eta, -K_{\mathcal X/B}|_{Y}) < a(X_\eta, -K_{\mathcal X/B})$.  This is equivalent to saying that $-\phi^{*}K_{\mathcal{X}/B}|_{Y} \cdot F \geq 3$.  Rearranging and substituting in our value for $k$, we find
\begin{align*}
\frac{1}{2}(-K_{\widetilde{Y}/B} \cdot \widetilde{C})  \leq \frac{3}{2}\neg(\mathbb{F},-K_{\mathbb{F}/B}) -\neg(\mathcal{X},-K_{\mathcal{X}/B}) + 3g(B) - 1.
\end{align*}
Using the inequality $\neg(\mathbb{F},-K_{\mathbb{F}/B}) \leq g(B)$ from Corollary \ref{coro:isolatingsection}, we see that
\begin{align*}
-K_{\mathcal{X}/B} \cdot C & \leq -K_{\widetilde{Y}/B} \cdot \widetilde{C} + 3g(B) - 1 \\
& \leq -2\neg(\mathcal{X},-K_{\mathcal{X}/B}) + 12g(B) - 3.
\end{align*}
contradicting our height bound.  Thus we see that $a(Y_\eta, -K_{\mathcal X/B}|_{Y}) \geq a(X_\eta, -K_{\mathcal X/B})$.  Since the subvarieties $Y$ with larger generic $a$-invariant lie in a closed set but $M$ defines a dominant family, we see that we must have equality of generic $a$-invariants, i.e.~$-\phi^{*}K_{\mathcal{X}/B}|_{Y} \cdot F = 2$.  Applying Equation \eqref{eqn:toomuchdeformingequation2} we see that 
\begin{equation}
\label{equation:-negupper}
-\neg(\mathbb{F},-K_{\mathbb{F}/B}) \leq -\neg(\mathcal{X},-K_{\mathcal{X}/B}) + 3g(B) - 1.
\end{equation}
Then
\begin{align}
-K_{\widetilde{Y}/B} \cdot \widetilde{C} & \geq -\phi^{*}K_{\mathcal{X}/B}|_{Y} \cdot \widetilde{C} + 1 - 3g(B)  \nonumber \\
& \geq 1 + 5g(B) + 2 \max\{4g(B)-2, 5g(B)-1 - \neg(\mathcal X, -K_{\mathcal X/B}) \} \nonumber \\
& \label{eqn:toomuchdeformingequation4} \geq 1 + 4g(B)+ \neg(\mathbb{F},-K_{\mathbb{F}/B}) + 2 \max \{4g(B) - 2,2g(B)  -\neg(\mathbb{F},-K_{\mathbb{F}/B})\}
\end{align}
where we have used the fact that $\neg(\mathbb{F},-K_{\mathbb{F}/B}) \leq g(B)$ as proved in Corollary \ref{coro:isolatingsection}.  This verifies that the sections will satisfy the height bound of Corollary \ref{coro:curvescoveringfamily} with respect to $Y$.

Proposition \ref{prop:totalbreaking} (2) guarantees the existence of a dominant family of sections $\widetilde{C}_{1}$ on $Y$ which satisfy $\widetilde{C} \equiv \widetilde{C}_{1} + (-s - 2g(B))F$ where
\begin{equation*}
s = \frac{1}{2}(K_{Y/B} \cdot C) +\frac{1}{2}\neg(\mathbb{F},-K_{\mathbb{F}/B})+  \max \{4g(B) - 2,2g(B)  -\neg(\mathbb{F},-K_{\mathbb{F}/B})\} 
\end{equation*}
and $F$ denotes a general fiber of $\pi$.  Equation \eqref{eqn:toomuchdeformingequation4} shows that $-s-2g(B) \geq 1$.  In particular, by Proposition \ref{prop:totalbreaking} (2) $\widetilde{C}$ deforms to the union of $\widetilde{C}_{1}$ with $-s-g(B)$ general fibers of $Y$.

Since $\widetilde{C}$ and $\widetilde{C}_{1}$ differ by a multiple of $F$, we have
\begin{align*}
3g(B) -1 & \geq K_{\widetilde{Y}/B} \cdot \widetilde{C} - \phi^{*}K_{\mathcal{X}/B}|_{Y} \cdot \widetilde{C} \\
& = K_{\widetilde{Y}/B} \cdot \widetilde{C}_{1} - \phi^{*}K_{\mathcal{X}/B}|_{Y} \cdot \widetilde{C}_{1}.
\end{align*}
Note that
\begin{equation*}
-K_{\widetilde{Y}/B} \cdot \widetilde{C}_{1} = 4g(B) + \neg(\mathbb{F},-K_{\mathbb{F}/B})+ 2\max \{4g(B) - 3,2g(B)  -\neg(\mathbb{F},-K_{\mathbb{F}/B}) \}.
\end{equation*}
Using the bound $-\neg(\mathbb{F},-K_{\mathbb{F}/B}) \leq -\neg(\mathcal{X},-K_{\mathcal{X}/B}) + 3g(B) - 1$ from (\ref{equation:-negupper}), we conclude
\begin{align*}
- \phi^{*}K_{\mathcal{X}/B}|_{Y} \cdot \widetilde{C}_{1} & \leq 3g(B) -1 + 4g(B) + \neg(\mathbb{F},-K_{\mathbb{F}/B}) \\
& \hspace{.3in} + 2 \max \{4g(B) - 2,5g(B) -1 - \neg(\mathcal{X},-K_{\mathcal{X}/B}) \} \\
& \leq 8g(B) - 1 +  2\max \{4g(B) - 2,5g(B) - 1 - \neg(\mathcal{X},-K_{\mathcal{X}/B}) \}
\end{align*}
Now suppose we carry out this construction for every surface $Y$ obtained above.  Let $C_{1}$ denote the image of $\widetilde{C}_{1}$ in $\mathcal{X}$.  Using a Hilbert scheme argument, we see that the corresponding families of sections $C_{1}$ constructed for general $Y$ must be members of some fixed component $M_{1} \in \Sec(\mathcal{X}/B)$. Indeed, as we vary over general sets of $q$ points in $\mathcal{X}$ the surfaces $Y$ constructed above are parametrized by an irreducible variety.
  For each surface, we have constructed some sections $C_1$ which have bounded relative anticanonical degree in $\mathcal{X}$.  Since there are only finitely many components of $\Sec(\mathcal{X}/B)$ parametrizing sections of bounded degree, for a general surface $Y$ these sections are all parametrized by the same irreducible component.  Thus we have verified that these sections satisfy properties (1), (2), and (3).
\end{proof}

\begin{coro}
\label{coro:boundedsurfacesfordominant}
Let $\pi: \mathcal{X} \to B$ be a del Pezzo fibration such that $-K_{\mathcal{X}/B}$ is relatively nef.  There is a bounded family $\mathcal{S} \to N$ of surfaces $Y \subset \mathcal{X}$ with $a(Y_{\eta},-K_{\mathcal{X}/B}|_{Y}) = a(\mathcal{X}_{\eta},-K_{\mathcal{X}/B})$ such that every component $M \subset \Sec(\mathcal{X}/B)$ parametrizing a dominant family of non-relatively free sections $C$ with
\begin{align*}
-K_{\mathcal{X}/B} \cdot C \geq \max \{& 6g(B)-2 + m(\mathcal X/B), -2\neg(\mathcal{X},-K_{\mathcal{X}/B}) + 12g(B) - 2, \\
&  2\max \{4g(B) - 2,5g(B) - 1 - \neg(\mathcal{X},-K_{\mathcal{X}/B}) \}  + 8g(B) + 1 \},
\end{align*}
will generically parametrize sections contained in these surfaces $Y$.
\end{coro}

\begin{proof}
Let $q$ denote the maximal number of general points contained in deformations of $C$.  As in Theorem \ref{theo: non-relativelyfreecurve}, the set of sections through $q$ general points sweeps out a surface $Y$.  The theorem also shows that the general such $Y$ contains a section $C_{1}$ with
\begin{equation*}
-K_{\mathcal{X}/B} \cdot C_{1} \leq 8g(B) - 1 +  2\max \{4g(B) - 2,5g(B) - 1 - \neg(\mathcal{X},-K_{\mathcal{X}/B}) \}.
\end{equation*}
We conclude by Proposition \ref{prop:boundednessequalainv} that the set of surfaces $Y$ obtained in this way forms a bounded family.
\end{proof}

Later on we will need to compare the normal bundles of the sections $C$ and $C_{1}$ in Theorem \ref{theo: non-relativelyfreecurve}.  This will be accomplished by the following lemma.

\begin{lemm}
\label{lemm:normalsheavesunchanged}
Let $\pi: \mathcal{X} \to B$ be a del Pezzo fibration.
Let $Y\subset \mathcal X$ be a surface such that if $\phi: \widetilde{Y}\to Y$ denotes the minimal resolution 
then the general fiber $F$ of $\pi\circ \phi : \widetilde{Y} \to B$ is isomorphic to $\mathbb P^1$, satisfies $-K_{\mathcal X/B} \cdot F = 2$, and its image in $\mathcal{X}$ has at worst nodal singularities.
 
Let $C_1, C_2$ be sections of $\pi\circ \phi : \widetilde{Y} \to B$ such that 
\begin{itemize}
\item both $C_i$s are contained in the smooth locus of $\mathcal X$;
\item the difference $C_1 - C_2$ is numerically equivalent to a multiple of $F$, and;
\item there exists an integer $c$ such that $-\phi^*K_{\mathcal{X}/B} \cdot C_i + K_{\widetilde{Y}/B} \cdot C_i \leq c$ and $-K_{\widetilde{Y}/B} \cdot C_i \geq \max \{2g(B), c+2g(B)-1 \}$.
\end{itemize}
Then we have
\[
h^1(C_1, N_{C_1/\mathcal X}) = h^1(C_2, N_{C_2/\mathcal X}).
\]
\end{lemm}

\begin{proof}
For the surfaces $\widetilde{Y}$ as above, define the normal sheaf $N_{\widetilde{Y}/\mathcal{X}}$ as the cokernel of the map $T_{\widetilde{Y}} \to \phi^{*}T_{\mathcal{X}}$.  We first prove some properties of the normal sheaf.  We assume above that a general fiber of the map $\pi|_{Y}: Y \to B$ is a rational conic in a fiber of $\pi: \mathcal{X} \to B$ with only nodal singularities. Thus the support of the torsion subsheaf of $N_{\widetilde{Y}/\mathcal{X}}$ does not intersect a general fiber of $Y$.
Now consider the exact sequence
\begin{equation*}
T_{\widetilde{Y}} \to \phi^{*}T_{\mathcal{X}} \to N_{\widetilde{Y}/\mathcal{X}} \to 0.
\end{equation*}
If we restrict this sequence to $C_i$ and remove the direct summands coming from $T_{C_i}$, we obtain an exact sequence
\begin{equation*}
N_{C_i/\widetilde{Y}} \to N_{C_i/\mathcal{X}} \to N_{\widetilde{Y}/\mathcal{X}}|_{C_i} \to 0
\end{equation*}
which must be exact on the left due to the fact that $N_{C_i/\widetilde{Y}}$ is an invertible sheaf.  Since $C_i$ is a general section in a family of sections on $\widetilde{Y}$ which satisfy $-K_{\widetilde{Y}/B} \cdot C_i \geq 2g(B)$, by Lemma \ref{lemm:expdimatpoints} we may ensure that $C_i$ avoids any codimension $2$ locus in $\widetilde{Y}$. 

In particular, let $(N_{\widetilde{Y}/\mathcal{X}})_{tf}$ denote the torsion-free quotient of $N_{\widetilde{Y}/\mathcal{X}}$.  Note that $(N_{\widetilde{Y}/\mathcal{X}})_{tf}$ is locally free away from a codimension $2$ subset, since it is locally free at every codimension $1$ point of $\widetilde{Y}$.  Thus we may ensure that $C_i$ is contained in the locus where $(N_{\widetilde{Y}/\mathcal{X}})_{tf}$ is locally free.  Altogether this implies that
\begin{align*}
(N_{\widetilde{Y}/\mathcal{X}}|_{C_i})_{tor} & = (N_{\widetilde{Y}/\mathcal{X}})_{tor}|_{C_i} \\
& = \mathcal{O}_{Z}|_{C_i}
\end{align*}
where $Z$ is the divisor whose components are the codimension $1$ components of the support of $(N_{\widetilde{Y}/\mathcal{X}})_{tor}$ and the multiplicities are the lengths of the torsion subsheaf along the generic point of each component.

We claim that the restriction of $N_{\widetilde{Y}/\mathcal{X}}$ to $C_1$ and to $C_{2}$ are isomorphic sheaves.  
We first prove that the restrictions of $(N_{\widetilde{Y}/\mathcal{X}})_{tor}$ agree.  Since a general fiber of $\widetilde{Y} \to B$ will not intersect the torsion subsheaf of $N_{\widetilde{Y}/\mathcal{X}}$ but $C_{1}$ and $C_{2}$ differ by a sum of general fibers, we see that the torsion subsheaves of $N_{\widetilde{Y}/\mathcal{X}}|_{C_{1}}$ and $N_{\widetilde{Y}/\mathcal{X}}|_{C_{2}}$ are isomorphic and in particular have the same length.
 
Next consider the restriction of $(N_{\widetilde{Y}/\mathcal{X}})_{tf}$.  Since a general fiber $F$ of $\widetilde{Y} \to B$ is $\mathbb P^1$ and the degree of the restriction of $(N_{\widetilde{Y}/\mathcal{X}})_{tf}$ to $F$ is zero, in the open locus of $\widetilde{Y}$ where the torsion free part is locally free the invertible sheaf $(N_{\widetilde{Y}/\mathcal{X}})_{tf}$ is linearly equivalent to a sum of vertical curves. Since $C_1$ and $C_2$ have the same intersection profile, we conclude that the restrictions of the torsion free part to $C_i$ are isomorphic.

Note that
\begin{equation*}
\deg(N_{\widetilde{Y}/\mathcal{X}}|_{C_{i}}) = -\phi^*K_{\mathcal{X}/B} \cdot C_i + K_{\widetilde{Y}/B} \cdot C_i \leq c.
\end{equation*}
On the other hand, we have
\begin{align*}
-K_{\widetilde{Y}/B} \cdot C_{i}
& \geq  c + 2g(B)-1\\
& \geq \deg(N_{\widetilde{Y}/\mathcal{X}}|_{C_{i}}) + 2g(B) - 1
\end{align*}
Let $\mathcal{K}$ be the kernel fitting in the exact sequence
\begin{equation*}
0 \to \mathcal{K} \to N_{C_i/\mathcal{X}} \to (N_{\widetilde{Y}/\mathcal{X}}|_{C_i})_{tf} \to 0.
\end{equation*}
The argument above shows that the degree of $\mathcal{K}$ is at least $2g(B)-1$ more than the degree of the rightmost term.  Thus $\Ext^{1}((N_{\widetilde{Y}/\mathcal{X}}|_{C_i})_{tf},\mathcal{K}) = 0$ and we conclude that the exact sequence splits.  By degree considerations $h^{1}(C_{i},\mathcal{K}) = 0$ and so $(N_{\widetilde{Y}/\mathcal{X}}|_{\widetilde{C}_{i}})_{tf}$ contributes all of $H^{1}(C_{i},N_{C_{i}/\mathcal{X}})$.  As explained above this implies that the dimension of this vector space is the same for both curves.
\end{proof}

\section{Stable maps through general points} \label{sect:genpoints}

Let $\pi: \mathcal{X} \to B$ be a del Pezzo fibration.  Suppose that $f: C \to \mathcal{X}$ is a genus $g(B)$ stable map obtained by deforming a section of $\pi$.  The goal of this section is to relate the number of general points of $\mathcal{X}$ contained in $f(C)$ with its anticanonical degree.  When $C$ is irreducible this is straightforward, but the situation is more subtle when $C$ is reducible.

In particular, we will show that when $f(C)$ contains (almost) as many points as possible given its degree then $C$ will have a particularly simple form.  In the next section, we use this fact to deduce that a section with semistable normal bundle will break in a particularly simple way, implying a special case of Movable Bend-and-Break.  

We first discuss irreducible curves:

\begin{lemm} \label{lemm:deformingheightbound}
Let $\pi: \mathcal{X} \to B$ be a del Pezzo fibration.  Fix $n \geq 2g(B)+1$. 
Let $M \subset \Sec(\mathcal{X}/B)$ denote a component such that for any $n$ general points of $\mathcal{X}$ there exists a member of $M$ containing those points.  Then $M$ has the expected dimension and the curves $C$ parametrized by $M$ have height $\geq 2n+2g(B)-2$.
\end{lemm}

\begin{proof}
Corollary \ref{coro:generalpointsandfreeness} shows that a general $C$ is relatively free so that $M$ has the expected dimension.  Since the points are general they will impose independent conditions on $M$, so that $\dim(M) \geq 2n$.  Thus
\begin{equation*}
2n \leq -K_{\mathcal{X}/B} \cdot C + 2(1-g(B))
\end{equation*}
which rearranges to give the desired height bound.
\end{proof}

We will use the following definition to control the failure of families of sections to have the expected dimension.

\begin{defi}
Let $\pi: \mathcal{X} \to B$ be a del Pezzo fibration.  Fix an integer $d$.  We let $\maxdef(d)$ denote the maximum dimension of any component $M \subset \Sec(\mathcal{X}/B)$ parametrizing sections of height $d$. When there is no section of height $d$, we simply set $\maxdef(d) = -\infty$. 
We also define
\begin{equation*}
\maxdef(\leq q) = \max \left\{ 0 , \max_{d \leq q} \maxdef(d) \right\}.
\end{equation*}
\end{defi}

\begin{lemm} \label{lemm:maxdefbound}
Let $\pi: \mathcal{X} \to B$ be a del Pezzo fibration such that $-K_{\mathcal{X}/B}$ is relatively ample.  Fix an integer $d$.  Choose a positive integer $n$ satisfying $n \geq \maxdef(d)$.  Suppose that we have a connected effective curve $C = C_{0} + \sum a_{i}T_{i}$ where $C_{0}$ is a section of height $d$ and each $T_{i}$ is a $\pi$-vertical rational curve.  If $C$ contains $n$ general points of $\mathcal{X}$ then $-K_{\mathcal{X}/B} \cdot C \geq  d + 3n - \frac{3}{2}\maxdef(d)$.
\end{lemm}

\begin{proof}
Let $M$ be the component of $\Sec(\mathcal{X}/B)$ containing $C_{0}$.  We split into several cases depending upon the dimension of the locus swept out by the sections in $M$. 

First suppose that $C_{0}$ is a rigid section.  In this case each general point contained in $C$ is contained in a $\pi$-vertical component.  There are only finitely many $\pi$-vertical conics through a general point of $\mathcal{X}$.  By generality, none of these conics will intersect $C_{0}$.  Thus any $\pi$-vertical curve through a general point of $\mathcal{X}$ that intersects $C_{0}$ must have anticanonical degree $\geq 3$.  Thus $-K_{\mathcal{X}/B} \cdot C \geq d + 3n$.

Next suppose that deformations of $C_{0}$ sweep out a surface $Y$.  Fix a general point $p$ of $\mathcal{X}$ and consider the finite set of $\pi$-vertical conics through this point.  Each such conic will intersect $Y$ in a finite set of points.  If a $\pi$-vertical conic through $p$ is a component of $C$, then $C_{0}$ must contain one of these intersection points with $Y$ to ensure that $C$ is connected.  Furthermore, the conditions imposed on $C_{0}$ by insisting that it contain the attachment points for conics through $p$ will be independent for different general points $p$.  Thus $C$ can have at most $\maxdef(d)$ components which are $\pi$-vertical conics through general points.  This proves that
\begin{align*}
-K_{\mathcal{X}/B} \cdot C & \geq d + 2\maxdef(d) + 3(n-\maxdef(d)) \\
& = d + 3n - \maxdef(d).
\end{align*}

Finally suppose that deformations of $C_{0}$ form a dominant family.  Suppose that $q$ of the general points are contained in $C_{0}$; this imposes $2q$ independent conditions on $M$.  In particular $2q \leq \maxdef(d)$.  Meeting a $\pi$-vertical conic through a fixed general point imposes one additional condition on $C_{0}$.  Thus $C$ can have at most $\maxdef(d) - 2q$ components which are $\pi$-vertical conics through general points.  This proves that
\begin{align*}
-K_{\mathcal{X}/B} \cdot C & \geq d + 2(\maxdef(d)-2q) + 3(n-q-(\maxdef(d)-2q)) \\
& \geq d + 3n - \maxdef(d) - q \\
& \geq d + 3n - \frac{3}{2}\maxdef(d).
\end{align*}
\end{proof}

The two following propositions describe which stable maps can pass through the maximal number of general points of $\mathcal{X}$.  The first handles the case of even height.

\begin{prop} \label{prop:evendegreebound}
Let $\pi: \mathcal{X} \to B$ be a del Pezzo fibration such that $-K_{\mathcal{X}/B}$ is relatively ample.  Fix a positive integer
\begin{equation*}
n \geq \frac{3}{2}\maxdef(\leq 6g(B) ) +2g(B)+ 2 + \max\{0, - \neg(\mathcal{X},-K_{\mathcal{X}/B})\}.
\end{equation*}
Suppose that $f: C \to \mathcal{X}$ is a genus $g(B)$ stable map with anticanonical height $2n+2g(B)-2$ such that the unique component of $C$ whose image is not $\pi$-vertical maps birationally to a section.  Then:
\begin{enumerate}
\item Suppose the image of $C$ contains $n$ general points of $\mathcal{X}$.  Then $f$ is a birational map to a relatively free section.
\item Fix a general curve $Z$ in a basepoint free linear series in a general fiber of $\pi$.  Suppose the image of $C$ contains $n-1$ general points of $\mathcal{X}$ and intersects $Z$.  Suppose also that the image of $C$ is reducible and at least one of our general points is contained in a $\pi$-vertical component of $C$.  Then $C$ has exactly two components and $f$ maps one component birationally onto a relatively free section and the other birationally onto either a conic or a cubic in a general fiber of $\pi$. Moreover the intersection of the section and the $\pi$-vertical curve is a smooth point of the $\pi$-vertical curve.
\end{enumerate}
\end{prop}

\begin{cons}
Before giving the proof, we clarify what ``general'' means in the statement of the proposition.  We let $Z$ be a fixed general member of a basepoint free linear system in a general fiber of $\pi$.  We may then ensure that the $n$ points satisfy the following conditions:
\begin{enumerate}
\item Suppose we fix any subset of our set of points of size $d$.  Then we may ensure that these points impose independent conditions on every family of height $\leq 2n+2g(B)-2$.  
In particular, when $d \geq 2g(B)+1$ then Corollary \ref{coro:generalpointsandfreeness} (3) shows that a general deformation of a section through all $d$ points will be relatively free.  Thus when choosing $d$ general points we may ensure that there are only finitely many sections of height $2d+2g(B)-2$ which contain this subset of points and that every section of this height containing the points will be relatively free.
Also we may ensure there are finitely many loci of dimension $\leq 1$ in $\Sec(\mathcal{X}/B)$ 
that parametrize sections of height $2d+2g(B)-1$ which contain this set of points and that the general section in these families is relatively free.  
\item We may ensure that every fiber of $\pi$ that contains one of the general points or $Z$ is a smooth del Pezzo surface. Moreover we may assume that each general point and $Z$ are contained in different fibers of $\pi$. Furthermore, we may ensure that a general point is not contained in any $(-1)$-curve in its fiber $F$ and that the finite set of conics in $F$ through the point have the expected dimension of intersection against the sections described in (1); namely, for any $0$-dimensional component of the parameter space the corresponding section is disjoint from each conic, and for any $1$-dimensional component of the parameter space there are only finitely many sections intersecting each conic.  We may also ensure that these finitely many sections are relatively free when they contain more than $2g(B)$ general points.
\item By choosing $Z$ general we may ensure that the sections in (1) and (2) meet $Z$ in the expected dimension: for any $0$-dimensional component of the parameter space the corresponding section is disjoint from $Z$, and for any $1$-dimensional component of the parameter space there are only finitely many sections intersecting $Z$.  We may further ensure that these finitely many sections are relatively free when they contain more than $2g(B)$ general points.
\item Consider a subset of $d$ of our general points where $d-1 \geq 2g(B)+1$.  Fix one of these general points $p$ and let $T$ be any $\pi$-vertical anticanonical conic which contains $p$.  By choosing the points general we may ensure that there are only finitely many sections of height $2d+2g(B)-2\leq 2n+2g(B)-2$ passing through the remaining $d-1$ general points and meeting with $Z$ and $T$.  Indeed, by Corollary \ref{coro:generalpointsandfreeness} (3) it suffices to consider families which generically parametrize relatively free sections.  Lemma \ref{lemm:expdimatpoints2} shows that a general member of such a family will not intersect $Z$.  
Due to the generality of the point $p$ we may ensure that no $\pi$-vertical conic through $p$ meets every section that intersects $Z$.  Finally, since general points impose codimension $2$ conditions on the remaining subfamily we obtain the desired claim. 
Furthermore, we may ensure that the finitely many sections satisfying these incidence conditions are relatively free. 
\item  Note that if we take a general point in a del Pezzo surface $S$ there are only finitely many rational anticanonical cubics containing this point and any other general point in $S$.  We ensure that all of our $n$ points satisfy this condition in their fiber.

In (1), (2), (3), and (4) we have constructed certain incidence conditions such that the families of sections satisfying these conditions are parametrized by a $0$-dimensional space.  For each choice of conditions, consider the set of $\pi$-vertical rational anticanonical cubics passing through one of our $n$ general points. These anticanonical cubics are parameterized by a $1$-dimensional family, and they sweep out the fiber $F$ of $\pi$ containing the point.  By choosing our remaining points general, we may ensure that sections satisfying the above incidence conditions go through a general point of $F$.
Using an incidence correspondence, we see that for general choices of $n$ points there will only be finitely many anticanonical cubics through one of our general points which meet these sections.  

\item We may ensure that there is no section as in (2) or (4) that intersects two different conics through two different points in our set.  Indeed, since each such conic imposes independent conditions on the parameter space of sections, the set of sections through two general conics cannot meet the $1$-dimensional locus parametrizing sections through the other general points and $Z$.  By a similar argument we may ensure that any section as in (1)-(4) that is parametrized by a $0$-dimensional set does not meet with any line intersecting $Z$.
\item Since $n \geq \maxdef(\leq 6g(B))$, we may ensure that for any $d<6g(B)+1$ and any subset of our set of points of size $\geq \maxdef(d)$ the points are general in the sense of Lemma \ref{lemm:maxdefbound}.
\end{enumerate}
\end{cons}

\begin{proof}[Proof of Proposition \ref{prop:evendegreebound}:]
Note that our hypotheses imply that the $\pi$-vertical connected components of $C \backslash C_{0}$ are trees of rational curves.

(1).  It suffices to show that the image of $f$ is irreducible.  Suppose otherwise, so that $f(C) = C_{0} + \sum_{i \in I} a_{i}T_{i}$ for some $\pi$-vertical curves $T_{i}$.   Let $d$ denote the height of $C_{0}$ and set $t_{i} = -K_{\mathcal{X}/B} \cdot T_{i}$, so that
\begin{equation*}
2n+2g(B)-2 = d + \sum_{i \in I} a_{i}t_{i}.
\end{equation*}
Note that $C_{0}$ can contain at most $\max \{ 0, \lfloor \frac{d}{2} \rfloor + 1 - g(B) \}$ general points of $\mathcal{X}$ if it is relatively free.  By Corollary \ref{coro:generalpointsandfreeness} $C_{0}$ can contain at most $2g(B)$ general points of $\mathcal{X}$ if it is not relatively free.  An irreducible vertical curve $T_{i}$ can contain at most $1$ general point and if it does then $t_{i} \geq 2$.   Let $I' \subset I$ denote the set of vertical curves that contain one of the general points.  

We now break the argument into several cases.

\noindent \textbf{Case 1: $d \geq 6g(B)$.}  
Then the number of general points contained in $f(C)$ is bounded above by the number of general points contained in $C_{0}$ and in the $T_{i}$.  Thus:
\begin{align*}
\max \left\{ 0, \left\lfloor \frac{d}{2} \right\rfloor  + 1 - g(B), 2g(B) \right\} + |I'| & \geq n \\
& = \frac{d}{2} + 1 - g(B) + \sum_{i \in I} a_{i}t_{i}/2.
\end{align*}
Since $d \geq 6g(B)$ the RHS is an upper bound for the LHS.  Thus the inequality above must be an equality.  This means that $d$ is even, that $C_{0}$ goes through $\left\lfloor \frac{d}{2} \right\rfloor  + 1 - g(B)$ general points, that $C_{0}$ is relatively free (since it contains at least $2g(B)+1$ general points), that each $a_{i}=1$, and that each component of $T$ is a free vertical curve through one of the general points with $t_{i} = 2$.  In particular, the set of $\frac{d}{2}+1-g(B)$ general points determines a finite number of possibilities for $C_{0}$, and each vertical curve is also determined by a general point up to a finite set of possibilities.  If there are any vertical components, then for general choices $f(C)$ will not be connected, an impossibility.  Thus $f(C)$ is irreducible.

\noindent \textbf{Case 2: $d < 6g(B)$.}
Due to our lower bound on $n$ we may apply Lemma \ref{lemm:maxdefbound}.  It shows that $C$ must have height 
\begin{align*}
-K_{\mathcal{X}/B} \cdot C & \geq \neg(\mathcal{X},-K_{\mathcal{X}/B}) + 3n - \frac{3}{2}\maxdef(\leq 6g(B)-1)
\end{align*}
Applying our assumption on the lower bound for $n$, we see that
\begin{align*}
-K_{\mathcal{X}/B} \cdot C > 2n + 2g(B)-2
\end{align*}
contradicting our assumption that $C$ have height $2n+2g(B)-2$.

\bigskip

(2).  Write $f(C) = C_{0} + \sum_{i \in I} a_{i}T_{i}$ for some $\pi$-vertical curves $T_{i}$.   Let $d$ denote the height of $C_{0}$ and set $t_{i} = -K_{\mathcal{X}/B} \cdot T_{i}$, so that
\begin{equation*}
2n + 2g(B) - 2 = d + \sum_{i \in I} a_{i}t_{i}.
\end{equation*}
 Let $I' \subset I$ denote the set of vertical curves that contain one of the general points.  Again we separate into cases:

\noindent \textbf{Case 1: $d \geq 6g(B) +1$ and $C_{0}$ intersects $Z$.} 
Meeting $Z$ will impose one condition on a relatively free section.  Since $C_{0}$ intersects $Z$ it can contain at most $\max \{ 0, \lfloor \frac{d+1}{2} \rfloor - g(B), 2g(B) \}$ general points of $\mathcal{X}$.   The number of general points contained in $f(C)$ is bounded above by the number of general points contained in $C_{0}$ and in the $T_{i}$.  Thus:
\begin{align*}
\max \left\{ 0, \left\lfloor \frac{d+1}{2}  \right\rfloor - g(B), 2g(B) - 1 \right\} + |I'| & \geq n-1 \\
& \geq \frac{d}{2} - g(B) + \sum_{i \in I} a_{i}t_{i}/2.
\end{align*}
If $d$ is even, then we must have equality everywhere.  This means that $C_{0}$ contains the maximal number of general points, the set of deformations of $C_{0}$ which meet $Z$ and go through the maximal number of points is at most $1$-dimensional, that every $a_{i}=1$, and that each component of $T$ is a free vertical curve through one of the general points with $t_{i} = 2$.  However, for $C_{0}$ to meet a vertical conic through a general point is a codimension $1$ condition, so by generality there can be at most one vertical component and we obtain the desired expression. Moreover the intersection of the section and the vertical curve must be a smooth point of the vertical curve due to generality.  Since $C_0$ contains at least $2g(B)+1$ general points, it must be relatively free by generality.

If $d$ is odd, then the equation simplifies to $1/2 + |I'| \geq \sum_{i \in I} a_{i}t_{i}/2$.  There are only three options for the vertical components:
\begin{enumerate}
\item every component of $T$ has anticanonical degree $2$ and contains a general point, or
\item every component of $T$ but one has anticanonical degree $2$ and contains a general point, and the last one has anticanonical degree $3$ and contains a general point, or,
\item every component of $T$ but one has anticanonical degree $2$ and contains a general point, and the last component has anticanonical degree $1$.
\end{enumerate}
Note that there are only finitely many deformations of $C_{0}$ which meet $Z$ and go through $\frac{d+1}{2}-g(B)$ points.  Thus such $C_0$ is relatively free by generality as it contains at least $2g(B)+1$ general points. Recall that by assumption there is a vertical component of $f(C)$ through a general point.  However, since there are only finitely many vertical conics through a general point, by generality no such conic can intersect $C_{0}$.  This rules out the first and third situations, showing that $C$ must be the union of a free section and a cubic in a fiber. Moreover there are only finitely many cubics containing a general point and meeting with $C_0$ thus by generality the intersection of the section and the cubic is a smooth point of the cubic.

\noindent \textbf{Case 2: $d \geq 6g(B) +1$ and $C_{0}$ does not intersect $Z$.}
Just as before, the number of general points contained in $f(C)$ is bounded above by the number of general points contained in $C_{0}$ and in the $T_{i}$.  Thus:
\begin{align*}
\max \left\{ 0, \left\lfloor \frac{d}{2} \right\rfloor + 1 - g(B),2g(B)-1 \right\} + |I'| & \geq n - 1 \\
& \geq \frac{d}{2} -g(B) + \sum_{i \in I} a_{i}t_{i}/2.
\end{align*}
Since $C_{0}$ does not intersect $Z$, there must be a vertical curve that does intersect $Z$ but does not contain any general points, so that
\begin{equation*}
\sum_{i \in I} a_{i}t_{i}/2 - |I'| \geq \sum_{i \in I} a_{i}t_{i}/2 - \sum_{i \in I'} a_{i}t_{i}/2 \geq 1/2.
\end{equation*}
Thus when $d$ is even, $C_{0}$ must contain the maximal number of points $\frac{d}{2}+1 - g(B)$ and must be relatively free by generality.  In this situation the difference between $|I'|$ and $\sum_{i \in I'} a_{i}t_{i}/2$ is at most $1/2$.  There are four options for the vertical components:
\begin{enumerate}
\item every component of $T$ but one has anticanonical degree $2$ and contains a general point, and the last has anticanonical degree $1$ and meets $Z$.
\item every component of $T$ but two has anticanonical degree $2$ and contains a general point, one has anticanonical degree $1$ and meets $Z$, and the last one has anticanonical degree $3$ and contains a general point,
\item every component of $T$ but two has anticanonical degree $2$ and contains a general point, and the last two components have anticanonical degree $1$, one of which meets $Z$,
\item every component of $T$ has anticanonical degree $2$, all but one contain a general point, and the last component meets $Z$.
\end{enumerate}
Since by generality there are only finitely many deformations of $C_{0}$ through the required number of points, such $C_{0}$ can not intersect a vertical line meeting $Z$ or a vertical conic through a general point.  This rules out the first three cases immediately, and the fourth is also ruled out since by assumption there exists at least one vertical curve which contains some general point.

When $d$ is odd then $C_{0}$ must contain the maximal number of points $\lfloor \frac{d}{2} \rfloor+1-g(B)$.  In this case there is only one option: every component of $T$ but one has anticanonical degree $2$ and contains a general point, and the last has anticanonical degree $1$ and meets $Z$.  However, since by generality $C_{0}$ can only deform in a one-parameter family while containing the maximal number of points, it is impossible for $C_{0}$ to meet both a line intersecting $Z$ and a conic through a general point.  Since by assumption $f(C)$ contains a vertical component through a general point, this case is also ruled out.

\noindent \textbf{Case 3: $d < 6g(B) +1$.}
In this case Lemma \ref{lemm:maxdefbound} shows that $C$ must have height
\begin{align*}
-K_{\mathcal{X}/B} \cdot C & \geq \neg(\mathcal{X},-K_{\mathcal{X}/B}) + 3(n-1) - \frac{3}{2}\maxdef(\leq 6g(B)) \\
& > 2n + 2g(B)-2
\end{align*}
proving the impossibility of this case.
\end{proof}

The next proposition is the analogue of Proposition \ref{prop:evendegreebound} for sections of odd height.

\begin{prop} \label{prop:odddegreebound} 
Let $\pi: \mathcal{X} \to B$ be a del Pezzo fibration such that $-K_{\mathcal{X}/B}$ is relatively ample.  Fix a positive integer
\begin{equation*}
n \geq \frac{3}{2}\maxdef(\leq 6g(B))+ 2g(B) + \max\{0, - \neg(\mathcal{X},-K_{\mathcal{X}/B})\}.
\end{equation*}
Suppose that $f: C \to \mathcal{X}$ is a genus $g(B)$ stable map with anticanonical height $2n + 2g(B)-1$ such that the unique component of $C$ whose image is not $\pi$-vertical maps birationally to a section.  Then:
\begin{enumerate}
\item Fix a general curve $Z$ in a basepoint free linear series in a general fiber of $\pi$. Suppose the image of $C$ contains $n$ general points of $\mathcal{X}$ and intersects $Z$. Then $f$ is a birational map to a relatively free section.
\item Suppose the image of $C$ contains $n$ general points of $\mathcal{X}$.  Suppose also that the image of $C$ is reducible and at least one general point is contained in a $\pi$-vertical component of $C$.  Then $C$ has exactly two components and $f$ maps one component birationally onto a relatively free section and the other birationally onto a conic or a cubic in a general fiber of $\pi$. Moreover the intersection of the section and the $\pi$-vertical curve is a smooth point of the $\pi$-vertical curve.
\end{enumerate}
\end{prop}

The proof is essentially the same as the proof of Proposition \ref{prop:evendegreebound}, but slightly easier.

\section{Movable Bend-and-Break for del Pezzo fibrations}

We next establish Movable Bend-and-Break for sections of del Pezzo fibrations such that $-K_{\mathcal{X}/B}$ is relatively ample.

\begin{theo} \label{theo:mbbmain}
Let $\pi: \mathcal{X} \to B$ be a del Pezzo fibration such that $-K_{\mathcal{X}/B}$ is relatively ample.  Define
\begin{align*}
Q(\mathcal{X}) = \max \{& 10g(B)+3, 3\maxdef(\leq 6g(B) ) +6g(B)+ 3 + 2\max\{0, - \neg(\mathcal{X},-K_{\mathcal{X}/B})\},\\
&  -2\neg(\mathcal{X},-K_{\mathcal{X}/B})  + 12g(B) +5,  -2\neg(\mathcal{X},-K_{\mathcal{X}/B}) + 21g(B) - 3,\\
& 3 \maxdef(\leq 8g(B) - \neg(\mathcal{X},-K_{\mathcal{X}/B})) - \neg(\mathcal{X},-K_{\mathcal{X}/B}) + 6g(B) + 1\}.
\end{align*}
Suppose that $M \subset \Sec(\mathcal{X}/B)$ is a component that parametrizes a family of relatively free sections $C$ satisfying $-K_{\mathcal{X}/B} \cdot C \geq Q(\mathcal{X})$.  Then the closure of $M$ in $\overline{M}_{g(B),0}(\mathcal{X})$ contains a point representing a stable map whose domain has exactly two components, one which maps birationally onto a relatively free section and one which maps birationally onto a free $\pi$-vertical curve. Moreover the section and the vertical curve meet at a smooth point of the vertical curve.
\end{theo}

It is helpful to introduce the following definition:

\begin{defi}
Let $\pi: \mathcal{X} \to B$ be a del Pezzo fibration such that $-K_{\mathcal{X}/B}$ is relatively ample.  We define
$\MBB(\mathcal{X})$ to be the smallest non-negative integer $r$ such that any component of $\overline{M}_{g(B),0}(\mathcal{X})$ which generically parametrizes  relatively free sections of height at least $r$  will also parametrize 
a stable map whose domain has exactly two components, one which maps birationally onto a relatively free section and one which maps birationally onto a free $\pi$-vertical curve, and such that the section intersects a smooth point of the vertical curve.
\end{defi}

Thus Theorem \ref{theo:mbbmain} establishes that $\MBB(\mathcal{X})$ exists and gives an explicit upper bound $\MBB(\mathcal{X}) \leq Q(\mathcal{X})$.

We will prove Theorem \ref{theo:mbbmain} by splitting into several cases.  The first situation (Cases (1) and (2)) is when $N_{C/\mathcal{X}}$ is ``almost'' semistable.  In this case we can appeal to Proposition \ref{prop:evendegreebound} and Proposition \ref{prop:odddegreebound} to obtain a suitable breaking of $C$.  The second situation (Cases (3) and (4)) is when $N_{C/\mathcal{X}}$ is unstable.  By taking all deformations of $C$ through a suitable number of general points we obtain a surface $\Sigma$ as in Section \ref{subsec:Shen}.  Since $C$ deforms a lot in $\Sigma$, we can break $C$ inside $\Sigma$ and then show that this breaking satisfies the desired properties.  Note that in this case the argument depends upon whether the degree of the destabilizing quotient of $N_{C/\mathcal{X}}$ is large or small.

\begin{proof}[Proof of Theorem \ref{theo:mbbmain}:]
Suppose that $M$ is a component of $\Sec(\mathcal{X}/B)$ parametrizing a relatively free family of sections that satisfy
\begin{align*}
-K_{\mathcal{X}/B} \cdot C \geq Q(\mathcal{X}).
\end{align*}
Since each section is smooth we can consider its normal bundle, and we separate into several cases based on the normal bundle of the general curve $C$ parametrized by $M$.

\textbf{Case 1:} $N_{C/\mathcal{X}}$ is semistable.
First assume that the height of $C$ is even.  Choose $n$ so that $-K_{\mathcal{X}/B} \cdot C = 2n + 2g(B)-2$. Our height bound implies that 
\[
n \geq \frac{3}{2}\maxdef(\leq 6g(B) ) +2g(B)+ \frac{5}{2} + \max\{0, - \neg(\mathcal{X},-K_{\mathcal{X}/B})\} \geq \frac{5}{2},
\]
which implies that $n \geq 3$.  Since $C$ is relatively free Riemann-Roch implies that
\[
h^0(C, N_{C/\mathcal X}) = 2n.
\]
Fix $n-1$ general points and a basepoint free curve $Z$ in a general fiber of $\pi$.
Proposition~\ref{prop:pointsandnormalbundle} shows there is a one parameter family of deformations of $C$ passing through $n-1$ general points and meeting with $Z$. By Bend-and-Break as in Lemma~\ref{lemm:strongerbandb} $C$ degenerates to a stable map whose image is the union of a section with $\pi$-vertical rational curves where at least one $\pi$-vertical component contains a general point. Now our assertion follows from Proposition~\ref{prop:evendegreebound}.  

When the height is odd, choose $n$ so that $-K_{\mathcal X/B} \cdot C = 2n + 2g(B)-1$.  Our height bound implies that
\[
n \geq \frac{3}{2}\maxdef(\leq 6g(B))+ 2g(B) +2  + \max\{0, - \neg(\mathcal{X},-K_{\mathcal{X}/B})\} \geq 2.
\]
Since $C$ is relatively free Riemann-Roch implies that
\[
h^0(C, N_{C/\mathcal X}) = 2n+1.
\]
Fix $n$ general points..
Proposition~\ref{prop:pointsandnormalbundle} shows there is a one parameter family of deformations of $C$ passing through $n$ general points. By Bend-and-Break as in Lemma~\ref{lemm:strongerbandb} $C$ degenerates to a stable map whose image is the union of a section with $\pi$-vertical rational curves where at least one $\pi$-vertical component contains a general point. Now our assertion follows from Proposition~\ref{prop:odddegreebound}.

\textbf{Case 2:} $N_{C/\mathcal{X}}$ is unstable with the maximal destabilizing exact sequence
\[
0 \to \mathcal L_2 \to N_{C/\mathcal X} \to \mathcal L_1 \to 0
\]
such that $\deg(\mathcal{L}_{2}) - \deg(\mathcal{L}_{1}) \leq 1$.  
Due to the relative freeness of the curve and the bound on the degrees of the components we have
\begin{equation*}
h^{0}(C,\mathcal{L}_{1}) = \left\lfloor \frac{h^{0}(C,N_{C/\mathcal{X}})}{2} \right\rfloor.
\end{equation*}
Thus the same proof for Case 1 works without any modification.

\textbf{Case 3:} $N_{C/\mathcal{X}}$ is unstable with the maximal destabilizing exact sequence
\[
0 \to \mathcal L_2 \to N_{C/\mathcal X} \to \mathcal L_1 \to 0
\]
such that $3g(B) \leq \deg(\mathcal{L}_{1}) \leq \deg(\mathcal{L}_{2}) -2$. 

We write $a = H^{0}(C,\mathcal{L}_{1}) = \deg(\mathcal{L}_{1}) + 1 - g(B)$ and $b = H^{0}(C,\mathcal{L}_{2}) = \deg(\mathcal{L}_{2}) + 1 - g(B)$. Our assumption implies that $a \geq 2g(B)+1$.

By Proposition \ref{prop:pointsandnormalbundle} deformations of $C$ can contain $a$ general points of $\mathcal{X}$.  Consider a component $\Sigma$ of the locus swept out by the curves $C$ parametrized by $M$ through $a$ general points $\{x_{1},\ldots,x_{a}\}$; we may assume that the general curve sweeping out $\Sigma$ has generic normal bundle and $\{x_{1},\ldots,x_{a}\}$ are general on $C$.  We know that $\Sigma \subsetneq \mathcal{X}$ and a dimension count shows that it is a surface.  In fact, since each general point imposes $2$ independent conditions on the moduli space, there is a $(b-a)$-dimensional family of curves through these general points which sweep out $\Sigma$.

Let $\nu: \Sigma' \to \Sigma$ be the normalization and let $\phi: \widetilde{\Sigma} \to \Sigma$ be the minimal resolution.  
The results in Section~\ref{subsec:Shen} show that $\widetilde{\Sigma} \to \Sigma$ is finite on an open neighborhood of a general section $C$ sweeping out $\Sigma$ and the strict transforms $\widetilde{C}$ will define a dominant family of sections of $\widetilde{\Sigma}$ which have normal bundle $\mathcal{L}_{2}$. In particular the induced map $\psi: \widetilde{\Sigma} \to B$ has connected fibers.  Since $\deg(\mathcal{L}_{2}) \geq g(B)+2$ the usual Bend-and-Break theorem shows that the general fiber of $\psi$ is $\mathbb{P}^{1}$.

By Corollary \ref{coro:isolatingsection} we find a ruled surface $\mathbb{F}$ and a birational map $\rho: \widetilde{\Sigma} \to \mathbb{F}$ which is an isomorphism on a neighborhood of $\widetilde{C}$.  Let $C_{0}$ be a section of $\mathbb{F} \to B$ of minimal height and write $\rho_{*}\widetilde{C} \equiv C_{0} + kF$ where $F$ is some fixed general fiber. 

\begin{clai}  \label{clai:deformationofpushforward} 
The section $\rho_{*}\widetilde{C}$ deforms on $\mathbb{F}$ to the sum of a relatively free section on $\mathbb{F}$ and a general fiber of $\psi$.
\end{clai}

\begin{proof}[Proof of claim:] Assume for a contradiction that it does not.  By Theorem \ref{theo:bandbforsurfaces} we see that
\begin{equation*}
\max\{ 2, 2g(B)+1, 4g(B) + 1 -\neg(\mathbb{F},-K_{\mathbb{F}/B})\} > -K_{\mathbb{F}/B} \cdot \widetilde{C} = \deg(\mathcal{L}_{2}).
\end{equation*}
By Lemma \ref{lemm:indecomrk2} we have $-\neg(\mathbb{F},-K_{\mathbb{F}/B}) \geq -g(B)$, so that the supremum is achieved by $4g(B)+1 - \neg(\mathbb{F},-K_{\mathbb{F}/B})$ unless $g(B) = 0$ where the supremum is $2$.  Thus (even in the genus $0$ case) we have
\begin{equation} \label{eq:boundfrom4.7}
-\neg(\mathbb{F},-K_{\mathbb{F}/B}) \geq \deg(\mathcal{L}_{2}) - 4g(B) - 2.
\end{equation}
On the other hand, note that
\begin{align*}
\deg(\mathcal{L}_{1}) + \deg(\mathcal{L}_{2}) & = -K_{\mathcal{X}/B} \cdot C \\
& = -\phi^{*}K_{\mathcal{X}/B} \cdot \rho^{*}C_{0} - k \phi^{*}K_{\mathcal{X}/B} \cdot F \\
& \geq -\phi^{*}K_{\mathcal{X}/B} \cdot \widetilde{C}_{0} - k \phi^{*}K_{\mathcal{X}/B} \cdot F
\end{align*}
where the last inequality follows from the fact that $\rho^{*}C_{0} - \widetilde{C}_{0}$ is effective and $-\phi^{*}K_{\mathcal{X}/B}$ is $\psi$-relatively nef.  Also, by arguing on $\mathbb{F}$ we see that
\begin{equation}  \label{eqn:mbbargumentequation2}
\deg(\mathcal{L}_{2}) = -K_{\mathbb{F}/B} \cdot \rho_{*}\widetilde{C} = -K_{\mathbb{F}/B} \cdot C_{0} + 2k.
\end{equation}
Substituting and simplifying, we see that 
\begin{equation} \label{eqn:mbbargumentequation}
\deg(\mathcal{L}_{1}) \geq \neg(\mathcal{X},-K_{\mathcal{X}/B})  + k (- \phi^{*}K_{\mathcal{X}/B} \cdot F - 2) - \neg(\mathbb{F},-K_{\mathbb{F}/B}).
\end{equation}
We now split the argument into two cases.  First, suppose that $-K_{\mathcal{X}/B} \cdot F \geq 3$.   Then we obtain
\begin{align*}
\deg(\mathcal{L}_{1}) & \geq \neg(\mathcal{X},-K_{\mathcal{X}/B})  + k  - \neg(\mathbb{F},-K_{\mathbb{F}/B}) \\
& = \neg(\mathcal{X},-K_{\mathcal{X}/B})  + \frac{1}{2}\deg(\mathcal{L}_{2})  - \frac{3}{2}\neg(\mathbb{F},-K_{\mathbb{F}/B})
\end{align*}
Using the lower bound on $-\neg(\mathbb{F},-K_{\mathbb{F}/B})$ from Equation \eqref{eq:boundfrom4.7}, we find altogether 
\begin{equation*}
\deg(\mathcal{L}_{2}) \geq \deg(\mathcal{L}_{1}) \geq \neg(\mathcal{X},-K_{\mathcal{X}/B})  +  2 \deg(\mathcal{L}_{2}) - 6g(B) - 3.
\end{equation*}
Thus
\begin{equation*}
-K_{\mathcal{X}/B} \cdot C \leq 2 \deg(\mathcal{L}_{2})-2 \leq -2\neg(\mathcal{X},-K_{\mathcal{X}/B})  + 12g(B) +4
\end{equation*}
contradicting our degree bounds.

Second, suppose that $-K_{\mathcal{X}/B} \cdot F = 2$.  Then the inequality  \eqref{eqn:mbbargumentequation} and \eqref{eq:boundfrom4.7} yield 
\begin{align*}
\deg(\mathcal{L}_{1}) & \geq \neg(\mathcal{X},-K_{\mathcal{X}/B}) + \deg(\mathcal{L}_{2}) - 4g(B) - 2.
\end{align*}
In particular this means that
\begin{equation} \label{eqn:mmbargumentequation3}
-K_{\mathcal{X}/B} \cdot C = \deg(\mathcal{L}_{1}) + \deg(\mathcal{L}_{2}) \leq 2\deg(\mathcal{L}_{1}) - \neg(\mathcal{X},-K_{\mathcal{X}/B}) + 4g(B) + 2.
\end{equation}
Then we have
\begin{align*}
-K_{\mathcal{X}/B} \cdot C & = -\rho_{*}\phi^{*}K_{\mathcal{X}/B} \cdot \rho_{*}\widetilde{C} \\
& \geq -\phi^{*}K_{\mathcal{X}/B} \cdot \widetilde{C}_{0} + 2k
\end{align*}
where the last inequality follows from the relative nefness of $-K_{\mathcal{X}/B}$.  Combining with Equation \eqref{eqn:mbbargumentequation2} and \eqref{eq:boundfrom4.7}, we see that \begin{equation*}
2\deg(\mathcal{L}_{1}) - 2\deg(\mathcal{L}_{2}) + 8g(B) + 4 - \neg(\mathcal{X},-K_{\mathcal{X}/B}) > -\phi^{*}K_{\mathcal{X}/B} \cdot \widetilde{C}_{0}.
\end{equation*}
In particular, the image $C_{0}$ of $\widetilde{C}_{0}$ in $\mathcal{X}$ has height at most $8g(B) - \neg(\mathcal{X},-K_{\mathcal{X}/B})$.  Note that any $a$ general points of $\widetilde{\Sigma}$ will also be general in $\mathcal{X}$.  Thus by adding on $a$ general fibers to $\widetilde{C}_{0}$ and taking the image in $\mathcal{X}$, we obtain a reducible curve $C_{0} + S$ where $C_{0}$ is a section, $S$ is $\pi$-vertical, and the entire curve goes through $a$ general points of $\mathcal{X}$.  By Lemma \ref{lemm:maxdefbound} either
\begin{equation*}
-K_{\mathcal{X}/B} \cdot C_{0} + 2a \geq -K_{\mathcal{X}/B} \cdot C_{0} + 3a - \frac{3}{2} \maxdef(\leq 8g(B) - \neg(\mathcal{X},-K_{\mathcal{X}/B}))
\end{equation*}
or $a < \maxdef(\leq 8g(B) - \neg(\mathcal{X},-K_{\mathcal{X}/B}))$.  Together these show that $a \leq \frac{3}{2} \maxdef(\leq 8g(B) - \neg(\mathcal{X},-K_{\mathcal{X}/B}))$.  Then by Equation \eqref{eqn:mmbargumentequation3} and the fact that $\deg(\mathcal{L}_{1}) = a + g(B) - 1$ we have
\begin{equation*}
-K_{\mathcal{X}/B} \cdot C \leq 3 \maxdef(\leq 8g(B) - \neg(\mathcal{X},-K_{\mathcal{X}/B})) - \neg(\mathcal{X},-K_{\mathcal{X}/B}) + 6g(B) 
\end{equation*}
contradicting our degree bound.
\end{proof}

Having proved the claim, we now return to the main argument.  Since a general deformation of a relatively free section on $\mathbb F$ will avoid any codimension $2$ locus, we may assume that the curves $\widetilde{C}_{1}$ and $F$ constructed on $\mathbb{F}$ in Claim \ref{clai:deformationofpushforward} avoid the exceptional centers of $\widetilde{\Sigma} \to \mathbb F$. Thus by taking the strict transforms of these curves, we have shown that on $\widetilde{\Sigma}$ the curves $\widetilde{C}$ deform to curves of the form $\widetilde{C}_{1} + F$ where $F$ is a general fiber of $\pi$ and $\widetilde{C}_{1}$ is relatively free in $\widetilde{\Sigma}$ that is general in its deformation class.  Moreover since $\widetilde{C}$ avoids the preimage of the singularities in $\mathcal X$ and $\Sigma'$, we may ensure that $\widetilde{C}_{1}$ avoids the preimage of the singularities of $\mathcal X$ and $\Sigma'$.  Indeed since a relatively free section avoids any codimension $2$ loci, $\widetilde{C}_1$ avoids $0$-dimensional preimages of singularities. For $1$-dimensional preimages of singularities since $\widetilde{C}$ avoids such loci, the intersection number of $\widetilde{C}$ to such loci is $0$. Then since $\widetilde{C}_{1}$ shares the same intersection property with $\widetilde{C}$, we conclude that it avoids any preimage of singularities.  Moreover we may assume that the image of the intersection of $\widetilde{C}_1$ and $F$ is a smooth point of the image of $F$ in $\mathcal X$.

Our plan is to show that $\widetilde{C}_{1}$ is relatively free in $\mathcal{X}$ and that $F$ is free in $\mathcal{X}$, yielding the desired deformation of $C$.  We separate the argument into two cases depending on the stability of the normal bundle of $\widetilde{C}_{1}$.

If $N_{\widetilde{C}_1/\mathcal X}$ is semistable, then $\widetilde{C}_1$ is relatively free by Proposition~\ref{prop:pointsandnormalbundle} and Corollary~\ref{coro:generalpointsandfreeness} (2).
Indeed, it follows from Proposition~\ref{prop:normalbundleinsurface} that the normal bundle $N_{\widetilde{C}/\widetilde{\Sigma}}$ is given by $\mathcal L_2$. Then since $\deg (N_{\widetilde{C}_1/\widetilde{\Sigma}}) = \deg (\mathcal L_2) -2 \geq \deg(\mathcal L_1)$ and since $\widetilde{\Sigma}$ contains the $a$ general points $\{x_1, \cdots, x_a\}$ that we have fixed,
we conclude that a deformation of $\widetilde{C}_1$ also contains these $a$ general points so that $h^{0}(\widetilde{C}_{1},N_{\widetilde{C}_{1}/\mathcal{X}}) \geq 4g(B) + 2$ by Proposition~\ref{prop:pointsandnormalbundle} (2). Thus $\widetilde{C}_1$ is relatively free by Corollary~\ref{coro:generalpointsandfreeness} (2). Since $F$ is free on $\widetilde{\Sigma}$ and $\Sigma$ contains a general point of $\mathcal{X}$, the image of $F$ is free on $\mathcal{X}$, finishing the argument in this case. 

Suppose that $N_{\widetilde{C}_1/\mathcal X}$ is unstable with the maximal destabilizing exact sequence:
\[
0 \to \mathcal L_2' \to N_{\widetilde{C}_1/\mathcal X} \to \mathcal L_1' \to 0.
\] 
Then after reselecting $a$ general points $q_1, \cdots, q_a$ on $\widetilde{C}_1$ we have $$h^0(\widetilde{C}_1, N_{\widetilde{C}_1/\mathcal X}) = h^0(\widetilde{C}_1, N_{\widetilde{C}_1/\mathcal X}(-q_1 - \cdots -q_a)) + 2a.$$ Note that our surface $\Sigma$ does not depend on the choice of $q_i$'s because of Lemma~\ref{lemm:Shen_independent}.  This means that 
\[
h^0(\widetilde{C}_1, \mathcal L_1') = h^0(\widetilde{C}_1, \mathcal L_1'(-q_1 - \cdots -q_a)) + a \geq a.
\]
Corollary~\ref{coro:generalpointsandfreeness} (1) shows that $\widetilde{C}_{1}$ is relatively free on $\mathcal{X}$. 
Similarly, since $F$ is free on $\widetilde{\Sigma}$, it is also free on $\mathcal{X}$.

\textbf{Case 4:}  $N_{C/\mathcal{X}}$ is unstable with the maximal destabilizing exact sequence
\[
0 \to \mathcal L_2 \to N_{C/\mathcal X} \to \mathcal L_1 \to 0
\]
such that $\deg(\mathcal{L}_{1}) < 3g(B)  < 7g(B) + 3 \leq \deg(\mathcal{L}_{2})$.   
As before we write $a = H^{0}(C,\mathcal{L}_{1}) = \deg(\mathcal{L}_{1}) + 1 - g(B) \leq 2g(B)$ and $b = H^{0}(C,\mathcal{L}_{2}) = \deg(\mathcal{L}_{2}) + 1 - g(B)$. 
By Proposition \ref{prop:pointsandnormalbundle} deformations of $C$ can contain $a$ general points of $\mathcal{X}$.  Consider a component $\Sigma$ of the locus swept out by the curves in $M$ through $a$ general points $\{x_{1},\ldots,x_{a}\}$; we may assume that the normal bundle of the general curve sweeping out $\Sigma$ has a Harder-Narasimhan filtration of generic type and $\{x_{1},\ldots,x_{a}\}$ are general on $C$.  We know that $\Sigma \subsetneq \mathcal{X}$ and a dimension count shows that it is a surface.  In fact, since each general point imposes $2$ independent conditions on the moduli space, there is a $(b-a)$-dimensional family of curves through these general points which sweep out $\Sigma$.

Let $\phi: \widetilde{\Sigma} \to \Sigma$ be the minimal resolution. The results in Section~\ref{subsec:Shen} show that $\widetilde{\Sigma} \to \Sigma$ is finite on an open neighborhood of a general section $C$ sweeping out $\Sigma$ and the strict transforms $\widetilde{C}$ will define a dominant family of sections of $\widetilde{\Sigma}$ which have normal bundle $\mathcal{L}_{2}$. In particular the induced map $\psi: \widetilde{\Sigma} \to B$ has connected fibers.  Since $\deg(\mathcal{L}_{2}) \geq g(B)+2$ the usual Bend-and-Break theorem shows that the general fiber of $\psi$ is $\mathbb{P}^{1}$.

Apply Corollary \ref{coro:isolatingsection} to find a ruled surface $\mathbb{F}$ and a birational map $\rho: \widetilde{\Sigma} \to \mathbb{F}$ which is an isomorphism on a neighborhood of $\widetilde{C}$.  As in Proposition \ref{prop:totalbreaking} (1) we can write $\widetilde{C} \equiv \widetilde{C}_{0} + kF + T$ where $F$ is a general fiber of the map to $B$ and
\begin{equation*}
k =  \frac{1}{2}(-K_{\widetilde{\Sigma}/B} \cdot \widetilde{C}) - \frac{1}{2}\neg(\mathbb{F},-K_{\mathbb{F}/B}).
\end{equation*}
Since $N_{\widetilde{C}/\widetilde{\Sigma}} = \mathcal L_2$ has vanishing $H^1$, the component parametrizing deformations of $\widetilde{C}$ in $\Sigma$ has the expected dimension. Thus we have
\[
-K_{\widetilde{\Sigma}/B} \cdot \widetilde{C} + 1-g(B) \geq -K_{\mathcal X/B} \cdot C +2(1-g(B)) -2a
\]
Then
\begin{align}
5g(B) -1 & \geq 2a + g(B) - 1  \nonumber \\
& \geq K_{\widetilde{\Sigma}/B} \cdot \widetilde{C} - \phi^{*}K_{\mathcal{X}/B}|_{\Sigma} \cdot \widetilde{C} \nonumber \\
& =  K_{\widetilde{\Sigma}/B} \cdot \widetilde{C} - \phi^{*}K_{\mathcal{X}/B}|_{\Sigma} \cdot \widetilde{C}_{0} + k (-\phi^{*}K_{\mathcal{X}/B}|_{\Sigma} \cdot F) -\phi^{*}K_{\mathcal{X}/B}|_{\Sigma} \cdot T \nonumber \\
& \label{eqn:toomuchdeformingequation3} \geq  K_{\widetilde{\Sigma}/B} \cdot \widetilde{C} + \neg(\mathcal{X},-K_{\mathcal{X}/B}) + k (-\phi^{*}K_{\mathcal{X}/B}|_{\Sigma} \cdot F)
\end{align}
where we have used the fact that $-K_{\mathcal{X}/B}$ is relatively nef at the last step.
Suppose for a contradiction that there is an inequality $a(\widetilde{\Sigma}_\eta, -\phi^{*}K_{\mathcal X/B}|_{\Sigma}) < a(X_\eta, -K_{\mathcal X/B})$.  This is equivalent to saying that $-\phi^{*}K_{\mathcal{X}/B}|_{\Sigma} \cdot F \geq 3$.  Rearranging and substituting in our value for $k$, we find
\begin{align*}
\frac{1}{2}(-K_{\widetilde{\Sigma}/B} \cdot \widetilde{C})  - \frac{3}{2}\neg(\mathbb{F},-K_{\mathbb{F}/B}) \leq -\neg(\mathcal{X},-K_{\mathcal{X}/B}) + 5g(B) - 1.
\end{align*}
Using the inequality $\neg(\mathbb{F},-K_{\mathbb{F}/B}) \leq g(B)$ from Corollary \ref{coro:isolatingsection}, we see that
\begin{align*}
-K_{\mathcal{X}/B} \cdot C & \leq -K_{\widetilde{\Sigma}/B} \cdot \widetilde{C} + 5g(B) - 1 \\
& \leq -2\neg(\mathcal{X},-K_{\mathcal{X}/B}) + 18g(B) - 3.
\end{align*}
contradicting our height bound.  Thus we see that $a(\Sigma_\eta, -K_{\mathcal X/B}|_{\Sigma}) \geq a(X_\eta, -K_{\mathcal X/B})$.  Since the subvarieties $\Sigma$ with larger generic $a$-invariant lie in a closed set but $C$ is a member of a dominant family, we see that we must have equality of generic $a$-invariants, i.e.~$-\phi^{*}K_{\mathcal{X}/B}|_{\Sigma} \cdot F = 2$. Using the arguments in Case 3 we conclude that $\rho_{*}\widetilde{C}$ deforms on $\mathbb{F}$ to the sum of a section $\widetilde{C}_1$ that is relatively free in $\widetilde{\Sigma}$ and a general fiber $F$ of $\psi$. 
When $\Sigma$ contains some singular points of $\mathcal X$, we may assume that $\widetilde{C}$ does not meet with the preimage of these singular points on $\widetilde{\Sigma}$. This implies that we can assume that $\widetilde{C}_1$ also avoids this preimage.  Note that since $\Sigma$ contains a general point of $\mathcal{X}$, the general fiber of $\Sigma \to B$ will be a nodal rational curve of anticanonical degree $2$ by Lemma~\ref{lemm:conicsondP}. 
Applying Lemma~\ref{lemm:normalsheavesunchanged} to $\widetilde{C}$ and $\widetilde{C}_{1}$ with the surface $\Sigma$ and the constant $c = 5g(B)-1$, we see that
\[
h^1(\widetilde{C}_1, N_{\widetilde{C}_1/\mathcal X}) = h^1(C, N_{C/\mathcal X})=0.
\]
Hence $\widetilde{C}_1$ is relatively free in $\mathcal X$.
\end{proof}

By applying Movable Bend-and-Break repeatedly, one can prove the following version.

\begin{coro} \label{coro:maxbreaking}
Let $\pi: \mathcal{X} \to B$ be a del Pezzo fibration such that $-K_{\mathcal{X}/B}$ is relatively ample.  Let $M$ be a component of $\Sec(\mathcal{X}/B)$ parametrizing a family of relatively free sections of height $\geq \MBB(\mathcal{X})$.  Then the closure $\overline{M}$ in $\overline{M}_{g(B),0}(\mathcal{X})$ contains a point representing the union of a relatively free section of height $< \MBB(\mathcal{X})$ with a chain of $\pi$-vertical free curves of anticanonical degree $\leq 3$ contained in a general fiber of $\pi$.  Furthermore the section meets the chain of $\pi$-vertical curves at a smooth point of the chain.
\end{coro}

The argument is very similar to the arguments in \cite[Section 7]{LTdPI}, and so we will only sketch the proof.

\begin{proof}
First, by applying Movable Bend-and-Break inductively we can find a stable map in $\overline{M}$ whose domain is a comb -- that is, the union of a curve $C_{0}$ of genus $g(B)$ with several rational curves meeting $C_{0}$ transversally at distinct points -- which maps $C_{0}$ onto a relatively free section of height $< \MBB(\mathcal{X})$ and each rational curve onto a free $\pi$-vertical curve.  By deforming and gluing the $\pi$-vertical components as in \cite[Theorem 7.7]{LTdPI} using Lemma~\ref{lemm: gluingandsmoothing}, one can then find a stable map in $\overline{M}$ whose domain is the union of a genus $g(B)$ curve and a single rational curve such that the genus $g(B)$ curve is mapped to a relatively free section of height $< \MBB(\mathcal{X})$ and the rational curve is mapped to a free $\pi$-vertical curve in a general fiber $F$ of $\pi$. Moreover we may assume that the intersection of the section and the vertical rational curve is a smooth point of the vertical curve. Finally, applying \cite[Lemma 2.12]{LTdPI} we can deform this $\pi$-vertical curve in $F$ to a chain of anticanonical conics and cubics while fixing the attachment point to $C_{0}$.
\end{proof}

By regluing all but one of the vertical components to the section, we obtain the following version:

\begin{coro} \label{coro:minbreaking}
Let $\pi: \mathcal{X} \to B$ be a del Pezzo fibration such that $-K_{\mathcal{X}/B}$ is relatively ample.  Let $M$ be a component of $\Sec(\mathcal{X}/B)$ parametrizing a family of sections of height $\geq \MBB(\mathcal{X})$.  Then the closure $\overline{M}$ in $\overline{M}_{g(B),0}(\mathcal{X})$ contains a point representing the union of a relatively free section with a $\pi$-vertical free curve of anticanonical degree $\leq 3$.  Furthermore the section meets the $\pi$-vertical curve at a smooth point.
\end{coro}

We also obtain a statement about the fibers of the evaluation map for the universal family of sections.

\begin{theo} \label{theo:connectedfibers}
Let $\pi: \mathcal{X} \to B$ be a del Pezzo fibration with $-K_{\mathcal{X}/B}$ relatively ample.  Suppose that $M$ is a component of $\Sec(\mathcal{X}/B)_{d}$ where $d \geq \MBB(\mathcal{X}) + 2$ and that $M$ generically parametrizes relatively free sections.  Then any resolution of the evaluation map of the universal family over $M$ has connected fibers.
\end{theo}

The proof is similar to \cite[Corollary 7.10]{LTdPI}: let $Y \to \mathcal{X}$ denote the Stein factorization of the evaluation map.  We first use Movable Bend-and-Break to break off some $\pi$-vertical free rational curves on $\mathcal{X}$.  By construction these must be the images of some $\pi$-vertical free rational curves on $Y$.  Then by applying \cite[Theorem 6.2]{LT16} to the generic fiber of $\pi$ we see that $Y \to \mathcal{X}$ must be birational.  We omit the details here, referring the reader to \cite[Corollary 7.10]{LTdPI}.

\section{Batyrev's heuristic}

Suppose that $\pi: \mathcal{X} \to B$ is a del Pezzo fibration. 
Batyrev's heuristic for Manin's Conjecture indicates that the number of components of $\Sec(\mathcal{X}/B)_{d}$ should be bounded above by a polynomial in $d$.  In this section we prove this statement for  Fano fibrations of relative dimension $1$ and $2$.

\begin{prop} \label{prop:batyrevconjdim1}
Let $\pi: \mathcal{X} \to B$ be a Fano fibration of relative dimension $1$.  Fix a numerical class $\alpha \in N_{1}(\mathcal{X})$.  When the anticanonical degree of $\alpha$ is sufficiently large, there is at most one component of $\Sec(\mathcal{X}/B)$ representing $\alpha$.  Equivalently, for $d$ sufficiently large there are exactly $|\Lambda|$ components of $\Sec(\mathcal{X}/B)_{d}$ where $\Lambda$ denotes the number of intersection profiles.
\end{prop}

In particular this implies that there is a polynomial $P(d)$ such that the number of components of $\Sec(\mathcal{X}/B)_{d}$ is bounded above by $P(d)$.

\begin{proof}
For each intersection profile $\lambda_{i}$ Lemma \ref{lemm:intersectionprofilebirationalmap} yields a birational map $\phi_{i}: \mathcal{X} \to \mathbb{F}_{i}$ to a ruled surface $\mathbb{F}_{i}$.  The set of components of $\Sec(\mathcal{X}/B)_{d}$ with intersection profile $\lambda$ is in bijection with the set of components of $\Sec(\mathbb{F}_{i}/B)_{d}$.  Thus it suffices to prove the statement for ruled surfaces.

Let $\mathcal{E}$ denote a rank $2$ bundle on $B$ defining the $\mathbb{P}^{1}$-bundle $\mathbb{F}_{i}$.  Suppose that the height $d$ is sufficiently large so that every twist of $\mathcal{E}$ by a divisor of degree $d$ is globally generated and has vanishing $H^{1}$.  Then there is a unique component of $\Sec(\mathbb{F}_{i}/B)$ of degree $d$.  This proves the first statement.

To see the second statement, note that in each intersection profile there is at most one numerical class of anticanonical degree $d$.  By combining with the argument above we obtain the second statement.
\end{proof}

We next turn to the del Pezzo fibration case.

\begin{theo} \label{theo:batyrevconjdim2}
Let $\pi: \mathcal{X} \to B$ be a del Pezzo fibration such that $-K_{\mathcal{X}/B}$ is relatively ample.  There is some polynomial $P(d)$ such that the number of components of $\Sec(\mathcal{X}/B)_{d}$ is bounded above by $P(d)$.
\end{theo}

The hardest case of the theorem is when the component $M$ parametrizes a dominant family of sections that is not relatively free so we explain this case separately.

\begin{prop} \label{coro:countingnonrelativelyfree}
Let $\pi: \mathcal{X} \to B$ be a del Pezzo fibration such that $-K_{\mathcal{X}/B}$ is relatively nef. There is an upper bound $R = R(\mathcal{X})$ on the number of components of $\Sec(\mathcal{X}/B)$ which parametrize a dominant family of non-relatively free curves representing a fixed numerical class.
\end{prop}

\begin{proof}
Suppose that $M$ parametrizes a dominant family of non-relatively-free sections on $\mathcal{X}$ of sufficiently large height and let $q$ be the maximal number of general points contained in the sections parametrized by $M$.  Fix the bounded family of surfaces $\mathcal{S} \to N$ defined by Corollary \ref{coro:boundedsurfacesfordominant}.  
Let $M^{[q]}$ be the space of sections parametrized by $M$ with $q$ marked points.  This is an irreducible variety; we denote its normalization by $\widetilde{M}^{[q]}$.  Consider the evaluation map $\widetilde{M}^{[q]} \to X^{q}$.  Theorem \ref{theo: non-relativelyfreecurve} 
shows that for any component of a general fiber of this map 
 the sections parametrized by this component will sweep out an irreducible
surface parametrized by $\mathcal{S}$.  Thus the Stein factorization of $\widetilde{M}^{[q]} \to X^{q}$ admits a rational map to the parameter space $\mathcal{N}$.
This means that if we let $\{ N_{i} \}$ denote the stratification of $N$ as in Lemma \ref{lemm:ruledsurfacedeform} there is a unique stratum $N_{i}$ which contains the general surface swept out by sections in $M$ through $q$ general points.
Corollary \ref{coro:curvescoveringfamily} shows that in fact for every surface parametrized by $N_{i}$ there will be a family of sections on a resolution of the surface which will be parametrized by some sublocus of $M$.  Thus the sublocus of $M$ swept out by families of sections on the surfaces parametrized by the stratum $N_{i}$ will contain a dense open subset of $M$.

In particular, fix a surface $\widetilde{Y}$ parametrized by $N_{i}$.  The component $M$ is uniquely determined by the anticanonical degree of the sections in $\widetilde{Y}$ and the intersection profile of any family of sections in $\widetilde{Y}$ which are contained in $M$.  Note that if we fix the numerical class of the sections in $\mathcal{X}$ and fix an intersection profile in $\widetilde{Y}$ then there is a unique family on $\widetilde{Y}$ which pushes forward to the given numerical class on $\mathcal{X}$ because our height is sufficiently large.
Thus we see that there is a universal upper bound on the number of families of curves that sweep out surfaces in the stratum $N_{i}$ and which push forward to a given numerical class on $\mathcal{X}$.  Since there are only finitely many strata $N_{i}$, we obtain the desired upper bound.  
\end{proof}

\begin{proof}[Proof of Theorem \ref{theo:batyrevconjdim2}:]
We consider separately non-dominant families, dominant but not relatively free families, and relatively free families.
Note that by Lemma \ref{lemm:northcott} in each case it suffices to prove the statement only for those sections whose anticanonical degree is sufficiently large.

Every non-dominant family of sections of sufficiently large height will sweep out a surface $Y$ in $\mathcal{X}$.  By Theorem \ref{theo:toomuchdeforming2} when the height is sufficiently large this surface $Y$ will be a ruled surface and by Corollary \ref{coro:nondomfambounded} only finitely many such surfaces can be obtained in this way.  It suffices to prove the statement for each surface $Y$ separately.  This follows from Proposition \ref{prop:batyrevconjdim1}.  (Note that the height function used here may be different from the height function in Proposition \ref{prop:batyrevconjdim1}.  Nevertheless, it is easy to show that a polynomial bound with respect to one height function will yield a polynomial bound for the other.)  

Suppose that $M \subset \Sec(\mathcal{X}/B)_{d}$ is a component parametrizing a dominant family of non-relatively free sections.  Proposition \ref{coro:countingnonrelativelyfree} proves that there is an upper bound $R$ on the number of such components representing any fixed numerical class, and this immediately implies the desired polynomial bound.

Finally, suppose that $M \subset \Sec(\mathcal{X}/B)_{d}$ is a component parametrizing a dominant family of sections which are generically relatively free.  By Corollary \ref{coro:maxbreaking} $M$ can be obtained by gluing a relatively free section $C$ of height $< \MBB(\mathcal{X})$ to a chain of $\pi$-vertical free curves each which has anticanonical degree $2$ or $3$.  By \cite[Theorem 5.13]{LT17} there is a polynomial $Q(r)$ which bounds the number of components of $\overline{M}_{0,0}(\mathcal{X})$ of anticanonical degree $r$ that can be obtained by smoothing chains of free curves of these types.  Let $c$ be the maximal degree of the Stein factorization of the evaluation map for any family of relatively free sections of height $< \MBB(\mathcal{X})$ and set $r = d + K_{\mathcal{X}/B} \cdot C$. Then by gluing $C$ to a smoothing of a chain of $\pi$-vertical curves of total anticanonical degree $r$ we can get at most $cQ(r)$ components of $\overline{M}_{g(B),0}(\mathcal{X})$.  By summing up these contributions over the finite set of families of relatively free sections of height $< \MBB(\mathcal{X})$, we obtain the desired bound.
\end{proof}

\section{Geometric Manin's Conjecture} \label{sect:maninconj}

In this section we recall the definition of the counting function in Geometric Manin's Conjecture and clarify its relationship to the structural theorems proved earlier.

Suppose that $\pi: \mathcal{X} \to B$ is a del Pezzo fibration such that $\mathcal X$ is smooth. Recall from Section~\ref{sec:preliminaries} that $N_1(\mathcal X)$ is defined to be the space of real $1$-cycles modulo numerical equivalence and $N_1(\mathcal X)_\bZ$ is the lattice of integral $1$-cycles.  Similarly, $N^1(\mathcal X)$ is the space of $\mathbb R$-divisors modulo numerical equivalence and $N^1(\mathcal X)_\bZ$ is the lattice of integral Cartier divisors.  Dual to the restriction map $N^{1}(\mathcal{X}) \to N^{1}(\mathcal{X}_{\eta})$ we have a pushforward map $N_{1}(\mathcal{X}_{\eta}) \to N_{1}(\mathcal{X})$.  Henceforth we will identify $N_{1}(\mathcal{X}_{\eta})$ with this subspace of $N_{1}(\mathcal{X})$.  Note however that the lattice $N_{1}(\mathcal{X}_{\eta})_{\mathbb{Z}}$ may be strictly contained in the lattice $N_{1}(\mathcal{X}_{\eta}) \cap N_{1}(\mathcal{X})_{\mathbb{Z}}$.

\subsection{Counting components}

Let $\pi: \mathcal{X} \to B$ be a Fano fibration.  A key piece of Geometric Manin's Conjecture is a precise bound on the number of components of $\Sec(\mathcal{X}/B)$ representing a fixed numerical class.  We expect the following principle to hold:

\begin{prin} \label{prin:uniquenessforalgebraicclasses}
Any ``sufficiently positive'' algebraic equivalence class of sections on $\mathcal{X}$ is represented by at most one family of relatively free sections.
\end{prin}

Here ``sufficiently positive'' is supposed to be taken with respect to the cone $\Nef_{1}(\mathcal{X}_{\eta})$ in $N_{1}(\mathcal{X})$.   We will give a more precise formulation of this statement for del Pezzo fibrations below.  The best evidence for Principle \ref{prin:uniquenessforalgebraicclasses} comes from fibrations of the form $X \times \mathbb{P}^{1}$ where $X$ is a Fano variety: there are many classes of Fano variety for which the components of $\Mor(\mathbb{P}^{1},X)$ have been classified and all examples that we know of are compatible with this principle (see Section~\ref{subsec:history}). 

In order to translate this statement to numerical equivalence, one must know how many algebraic equivalence classes there are representing a fixed numerical class.  Let $\mathcal{X}$ be a smooth threefold admitting a morphism $\pi$ to a curve $B$ whose fibers are connected and rationally connected.  For any section $C$ of $\pi$ the pushforward map $\CH_{0}(C) \to \CH_{0}(\mathcal{X})$ is surjective.  Thus by \cite[Theorem 10.17]{Voi03} we have $H^{2,0}(\mathcal{X}) = 0$ so that 
\begin{equation*}
\Br(\mathcal{X}) \cong H^{3}(\mathcal{X}, \mathbb{Z})_{\mathrm{tors}}.
\end{equation*}
According to the universal coefficient theorem, we can equally well think of $\Br(\mathcal{X})$ as the torsion classes of $H_{2}(\mathcal{X},\mathbb{Z})$.  Let $Q_{1}(\mathcal{X})$ denote the set of algebraic equivalence classes of curves of $\mathcal{X}$.  \cite[Theorem 1]{BS83} shows that algebraic and homological equivalence coincide for curve classes on $\mathcal{X}$ and \cite[Theorem 2]{Voisin06} proves the integral Hodge conjecture for $\mathcal{X}$.  Together these show:

\begin{theo}[\cite{BS83}, \cite{Voisin06}] \label{theo:algequiv}
Let $\mathcal{X}$ be a smooth threefold admitting a morphism to a curve $B$ whose fibers are connected and rationally connected.  Then $|\Br(\mathcal{X})|$ is the size of the kernel of the quotient map $q: Q_{1}(\mathcal{X}) \to N_{1}(\mathcal{X})_{\mathbb{Z}}$.
\end{theo}

Thus, Principle \ref{prin:uniquenessforalgebraicclasses} can be interpreted as follows:  

\begin{conj} \label{conj:numberofcomponents}
Let $\pi: \mathcal{X} \to B$ be a del Pezzo fibration with $\mathcal{X}$ smooth.  Fix an intersection profile $\lambda$ and let $N_{\lambda}$ denote the affine subset of $N_{1}(\mathcal{X})$ consisting of curve classes whose intersection numbers with $\pi$-vertical divisors are described by $\lambda$.  Define $\Nef_{\lambda} = \Nef_{1}(\mathcal{X}) \cap N_{\lambda}$.  There is some translate $\mathcal{T}$ of $\Nef_{\lambda}$ in $N_{\lambda}$ such that every class in $\mathcal{T}_{\mathbb{Z}}$ is represented by exactly $|\Br(\mathcal{X})|$ different families of relatively free sections.
\end{conj}

\subsection{Formulating Geometric Manin's Conjecture}
Throughout this section $\pi: \mathcal{X} \to B$ denotes a del Pezzo fibration such that $\mathcal{X}$ is smooth and $-K_{\mathcal{X}/B}$ is relatively ample.  For simplicity, we will assume that the general fibers are not isomorphic to $\mathbb{P}^{2}$ or $\mathbb{P}^{1} \times \mathbb{P}^{1}$.  (In these two cases one must adjust the counting function slightly to reflect the fact that the general fiber does not contain any curve class with anticanonical degree $1$.)

The most general version of Geometric Manin's Conjecture is formulated using the Fujita invariant to construct the exceptional set as in \cite{LST18}.  For del Pezzo fibrations, the general construction of the exceptional set in \cite{LST18} has the following simple description.

\begin{defi}[\cite{LTRes} Theorem 9.1] \label{defi:accumulatingcomponent}
We say that a component $M \subset \Sec(\mathcal{X}/B)$ is an accumulating component if either:
\begin{enumerate}
\item the sections sweep out a surface $Y$ whose intersection with a general fiber is an anticanonical line, or
\item $\rho(\mathcal X_\eta)  = 1$ and the sections sweep out a surface $Y$ whose intersection with a general fiber is a singular anticanonical conic.
\end{enumerate}
Any component $M$ which is not an accumulating component is called a Manin component.
We let $\Manin_{i}$ denote the set of Manin components that parametrize sections $C$ satisfying $-K_{\mathcal{X}/B} \cdot C = i$.
\end{defi}

\begin{rema}
Definition \ref{defi:accumulatingcomponent} (2) is very restrictive.  Since the condition $\rho(\mathcal{X}_{\eta})=1$ implies that the Picard group is generated by a multiple of $-K_{\mathcal{X}_{\eta}}$ the existence of an anticanonical conic means that $\mathcal{X}_{\eta}$ must be a del Pezzo surface of degree $\leq 2$.
\end{rema}

\begin{rema}
The relationship between Definition \ref{defi:accumulatingcomponent} and families of non-relatively free sections is a little subtle.  We will make this connection precise in the following section.
\end{rema}

The counting function in Geometric Manin's Conjecture encodes the number and dimension of Manin components representing sections of height at most $d$.

\begin{defi} \label{defi:countingfunction}
Fix a real number $q>1$.  For any positive integer $d$ define
\begin{equation*}
N(\mathcal{X},-K_{\mathcal{X}/B},q,d) := \sum_{i = 1}^{d} \sum_{M \in \Manin_{i}} q^{\dim M}.
\end{equation*}
\end{defi}

Geometric Manin's Conjecture predicts the asymptotic growth rate of this counting function as we let $d$ go to $\infty$.  The expected growth rate is
\begin{equation*}
N(\mathcal{X},-K_{\mathcal{X}/B},q,d) \sim_{d \to \infty} Cq^{d}d^{\rho(\mathcal{X}_{\eta}) - 1}
\end{equation*}
for some constant $C$.

\begin{rema} \label{rema:enoughsections}
Given a section $C$ of a del Pezzo fibration $\pi: \mathcal{X} \to B$, we can always construct a relatively free section as follows.  Choose a resolution $\phi: \widetilde{\mathcal{X}} \to \mathcal{X}$ and let $\widetilde{C}$ denote the strict transform of $C$.  By gluing $\widetilde{C}$ with suitably chosen free curves in the fibers of $\phi \circ \pi$ we can ensure that after smoothing we obtain a section $\widetilde{T}$ through $2g(B)+1$ general points of $\widetilde{\mathcal{X}}$.  By Corollary \ref{coro:generalpointsandfreeness} a general deformation of the pushforward of $\widetilde{T}$ to $\mathcal{X}$ will be relatively free.

By varying the choice of curves we glue on, it is not hard to show that the resulting relatively free sections generate a translate of a full-dimensional subcone of $\Nef_{1}(\mathcal{X}_\eta)$.  In particular, this shows that the asymptotic growth rate of $N(\mathcal{X},-K_{\mathcal{X}/B},q,d)$ is bounded below by $Cq^{d}d^{\rho(\mathcal{X}_{\eta})-1}$ for some constant $C$.
\end{rema}

\subsection{Upper bounds on the counting function} 
As in the previous section $\pi: \mathcal{X} \to B$ denotes a del Pezzo fibration such that $\mathcal{X}$ is smooth and $-K_{\mathcal{X}/B}$ is relatively ample and (for simplicity) we assume that the general fiber of $\pi$ is not isomorphic to $\mathbb{P}^{2}$ or $\mathbb{P}^{1} \times \mathbb{P}^{1}$.  

In order to verify Geometric Manin's Conjecture, we first must explicitly identify the accumulating components of $\Sec(\mathcal{X}/B)$.  We also must bound the contributions of the Manin components to the counting function.  The following claims summarize the roles of each type of component of $\Sec(\mathcal{X}/B)$ in Geometric Manin's Conjecture:
\begin{enumerate}
\item Most non-dominant families of sections will be accumulating components; the ones which are not will give a negligible contribution to the asymptotic growth of the counting function $N(\mathcal{X},-K_{\mathcal{X}/B},q,d)$.
\item Most dominant families of sections which are not relatively free will give a negligible contribution to the asymptotic growth of $N(\mathcal{X},-K_{\mathcal{X}/B},q,d)$; the ones which do not will be accumulating components.
\item Dominant families of sections which are generically relatively free will give a contribution to the asymptotic growth that can be computed using Movable Bend-and-Break.
\end{enumerate}
We will verify these claims in the proof of the following theorem.

\begin{theo}
Let $\pi: \mathcal{X} \to B$ be a del Pezzo fibration such that $\mathcal X$ is smooth, $-K_{\mathcal{X}/B}$ is relatively ample, and the general fiber is a del Pezzo surface that is not $\mathbb{P}^{2}$ or $\mathbb{P}^{1} \times \mathbb{P}^{1}$.  Then there is some non-negative integer $r$ such that
\begin{equation*}
N(\mathcal{X},-K_{\mathcal{X}/B},q,d) = O(q^{d} d^{r}).
\end{equation*}
\end{theo}

\begin{proof} 
Suppose we fix an intersection profile $\lambda$.  The set of sections with intersection profile $\lambda$ will be contained in an affine translate $\mathcal{N}$ of the subspace $N_{1}(\mathcal{X}_{\eta})$.  We will further restrict our attention to a fixed coset $\Xi$ of $N_{1}(\mathcal{X}_{\eta})_{\mathbb{Z}}$ inside of $\mathcal{N}_{\mathbb{Z}}$.  Since our goal is to prove an upper bound and since there are only finitely many cosets in all the intersection profiles, it suffices to prove an upper bound of the desired type for each coset individually.

Fix a coset $\Xi$ and choose a translation $\psi: \mathcal{N} \to N_{1}(\mathcal{X}_{\eta})$ which identifies the coset $\Xi$ with the lattice $N_{1}(\mathcal{X}_{\eta})_{\mathbb{Z}}$.  We now analyze separately the three different types of components of $\Sec(\mathcal{X}/B)$ which lie in $\Xi$.  Since we only care about the asymptotic behavior, it suffices to restrict our attention to the families of sections with sufficiently large height.  

\textbf{Case 1:} Consider a component $M \subset \Sec(\mathcal{X}/B)$ parametrizing a non-dominant family of sections of sufficiently large height.  By Theorem \ref{theo:toomuchdeforming2} there is a finite set of surfaces $Y$ satisfying $a(Y_{\eta},-K_{\mathcal{X}/B}) \geq a(\mathcal{X}_{\eta},-K_{\mathcal{X}/B})$ such that the sections parametrized by $M$ will sweep out one of these surfaces $Y$.  In the case when the inequality of a-invariants is strict the component $M$ will be accumulating.  If instead we have an equality of $a$-invariants, there are two cases to consider.  First, suppose that $\rho(\mathcal{X}_{\eta}) > 1$.   According to Proposition \ref{prop:batyrevconjdim1} the surface $Y$ can only contain finitely many components of relative anticanonical degree $d$.  On the other hand, due to the Picard rank assumption the number of components of $\Sec(\mathcal{X}/B)_{d}$ will grow at least linearly.  So in this case the contributions of $Y$ to the counting function are negligible.  Second, suppose that $\rho(\mathcal{X}_{\eta}) = 1$.  Note that in this situation the generic fiber $Y_{\eta}$ must be geometrically irreducible, yielding an anticanonical conic over the ground field in $\mathcal{X}_{\eta}$.  Due to the restriction on the Picard rank, the only option is that $\mathcal{X}_{\eta}$ is a del Pezzo surface of degree $\leq 2$ and Picard rank $1$ and that this conic is either:
\begin{enumerate}
\item a rational curve in $|-K_{X}|$ if $\mathcal{X}_{\eta}$ has degree $2$, or
\item a rational curve in $|-2K_{X}|$ if $\mathcal{X}_{\eta}$ has degree $1$.
\end{enumerate}
Note that each curve of these types is singular and thus is included in the exceptional set.  So every component of this type is an accumulating component.

\textbf{Case 2:} Consider the dominant families of non-relatively free sections with sufficiently large height.  Let $T_{1},\ldots,T_{r}$ denote all the families of rational anticanonical conics on $\mathcal{X}_{\eta}$ which are defined over the ground field and let $\beta_{1},\ldots,\beta_{r} \in N_{1}(\mathcal{X}_{\eta})$ denote their numerical classes (which may coincide).  We claim that there is some bounded subset $\mathcal{S} \subset N_{1}(\mathcal{X}_{\eta})$ such that the $\psi$-image of a class $\alpha$ of any dominant but not relatively free family is contained in $\mathcal{S} + m\beta_{j}$ for some index $j$ and some non-negative integer $m$.  This is an immediate consequence of Theorem \ref{theo: non-relativelyfreecurve} (3) which shows that any such class $\alpha$ is the sum of the class of a section of bounded height with some multiple of the class of a $\pi$-vertical conic.

We are now ready to analyze the contributions of such families to the counting function.  First, suppose that $\rho(\mathcal{X}_{\eta}) > 1$.  As demonstrated by Corollary \ref{coro:countingnonrelativelyfree} there is an upper bound $R$ on the number of components of dominant families of sections that are not relatively free which represent any fixed numerical class.  Furthermore, Lemma \ref{lemm:dominantexpecteddim} proves there is a constant $m(\mathcal{X}/B)$ such that $\dim(M) \leq -K_{\mathcal{X}/B} \cdot C + 2 - 2g(B) + m(\mathcal{X}/B)$ for all such components $M$.  Thus, the contribution of all the dominant non-relatively-free families contained in $\Xi$ to the counting function is bounded above by
\begin{equation*}
\sum_{i=1}^{d} R \cdot |\mathcal{S}| \cdot r \cdot q^{i + s}
\end{equation*}
for some constant $s$.  Since we are assuming that $\rho(\mathcal{X}_{\eta}) > 1$, this contribution to the counting function is asymptotically negligible compared to the growth rate as in Remark \ref{rema:enoughsections}.  Second, suppose that $\rho(\mathcal{X}_{\eta}) = 1$.  As explained in Case 1, in this situation the rational curves on surfaces swept out by fiberwise conics must be accumulating components.

\textbf{Case 3:} Consider the dominant families of relatively free sections with sufficiently large height.  In order to count the contributions of such components, our plan is to translate into a sum over lattice points in $N(\mathcal{X}_{\eta})$.

We first claim there is some class $\beta \in N_{1}(\mathcal{X}_{\eta})$ such that every class $\alpha$ of a dominant family of sections in our coset satisfies $\psi(\alpha) \in \beta + \Nef_{1}(\mathcal{X}_{\eta})$.  Indeed, by Theorem \ref{theo:mbbmain} we see that every dominant family of sections of sufficiently large height is numerically equivalent to a sum of a $\pi$-vertical free rational curve and a section of bounded height.  Furthermore, if we take a $\pi$-vertical free rational curve in a general fiber $F$ then its image under the map $N_{1}(F) \to N_{1}(\mathcal{X})$ is contained in $\Nef_{1}(\mathcal{X}_{\eta})$.  Thus we deduce that the $\psi$-image of all dominant families of sections lies in a finite union of translates of $\Nef_{1}(\mathcal{X}_{\eta})$, proving the claim.

We now use $\psi$ to translate the sum over to $N_{1}(\mathcal{X}_{\eta})_{\mathbb{Z}}$.  Let $b$ be the relative anticanonical degree of the class $\beta$ and let $c$ denote the constant which is the difference in relative anticanonical degrees between any numerical class in $\Xi$ and its image under $\psi$.  By definition every component $M$ of $\Sec(\mathcal{X}/B)$ that generically parametrizes relatively free sections has the expected dimension.  Furthermore, Theorem \ref{theo:batyrevconjdim2} gives a polynomial upper bound $P(d)$ on the number of components representing $\Sec(\mathcal{X}/B)_{d}$.  Altogether we see that the contribution of the relatively free sections to the counting function is bounded above by
\begin{equation*}
\sum_{\substack{\alpha \in \Nef_{1}(\mathcal{X}_{\eta})_{\mathbb{Z}} \\ -K_{\mathcal{X}/B} \cdot \alpha \leq d-b-c}} P( -K_{\mathcal{X}/B} \cdot \alpha + b + c) q^{-K_{\mathcal{X}/B} \cdot \alpha + b + c + 2 - 2g(B)}.
\end{equation*}
Using standard lattice counting techniques this shows the desired result.
\end{proof}

If we assume Conjecture \ref{conj:numberofcomponents} then we can compute explicitly the contributions of relatively free families.  By summing up over all cosets as in the proof of \cite[Theorem 9.10]{LTdPI}, we obtain the following:

\begin{theo} \label{theo:asymptoticformula}
Let $\pi: \mathcal{X} \to B$ be a del Pezzo fibration such that $\mathcal X$ is smooth, $-K_{\mathcal{X}/B}$ is relatively ample, and the general fiber is a del Pezzo surface of degree $\geq 2$ that is not $\mathbb{P}^{2}$ or $\mathbb{P}^{1} \times \mathbb{P}^{1}$.
Assume that Conjecture \ref{conj:numberofcomponents} holds. Then 
\begin{equation*}
N(\mathcal{X},-K_{\mathcal{X}/B},q,d) \mathrel{\mathop{\sim}_{d \to \infty}} \left( \tau_X \cdot \alpha (\mathcal{X}_\eta, -K_{\mathcal{X}/B}) \cdot |\Br(\mathcal{X})| \cdot \frac{q}{q-1} \right) q^{d} d^{\rho(\mathcal{X}_{\eta})-1}.
\end{equation*}
\end{theo}

Here the $\alpha$-constant of $\mathcal X_\eta$ is defined by 
\[
\alpha(\mathcal{X}_\eta, -K_{\mathcal{X}/B}) := \dim N_1(\mathcal X_\eta) \cdot \mu(\Nef_1(\mathcal X_\eta)\cap \{ \gamma \in N_1(\mathcal X_\eta) \, | -K_{\mathcal X_\eta} \cdot \gamma \leq 1\})
\]
where the volume is computed with respect to the lattice structure $N_{1}(\mathcal{X}_{\eta})_{\mathbb{Z}}$ and the $\tau$-constant of $\mathcal{X}$ is defined by
\begin{equation*}
\tau_{\mathcal{X}} = |\Lambda| \cdot [N_{1}(\mathcal{X})_{\mathbb{Z}} \cap N_{1}(\mathcal{X}_{\eta}) : N_1(\mathcal{X}_\eta)_{\mathbb Z}]
\end{equation*}
where $\Lambda$ denotes the set of allowable intersection profiles for $\mathcal{X}$.

\section{Geometric Manin's Conjecture for certain families}

In this section we prove Geometric Manin's Conjecture for certain special types of del Pezzo fibrations. 

\subsection{Singularities of fibers}
Suppose that $S$ is a (normal) del Pezzo surface with canonical singularities over an algebraically closed field of characteristic $0$.  \cite{HW81} shows that the minimal resolution of $S$ is a weak del Pezzo surface.

\begin{defi} \label{defi:systemoflines}
Let $S$ be a del Pezzo surface with canonical singularities.  We say that a line bundle $L$ on $S$ defines a system of lines if it satisfies $-K_{S} \cdot L = 3$ and $L^{2} = 1$.
\end{defi}

The terminology is motivated by the following description of such $L$.

\begin{lemm} \label{lemm:systemoflinestop2}
Let $S$ be a del Pezzo surface with canonical singularities and let $L$ define a system of lines on $S$.  Then $L$ is the pullback of $\mathcal{O}(1)$ under a birational morphism to $\mathbb{P}^{2}$.
\end{lemm}

\begin{proof}
Let $\phi: S' \to S$ be a minimal resolution so that $S'$ is a weak del Pezzo surface. We have $h^{0}(S,L) = h^{0}(S',\phi^{*}L) = 3$ by Riemann-Roch.  By \cite[Proposition 5.2.2.4]{ADHL} we see that $|\phi^{*}L|$ is basepoint free.  Thus the linear series defines a morphism to $\mathbb{P}^{2}$ and since $L^{2}=1$ this morphism must be birational.  We claim that $|L|$ is also basepoint free.  Indeed, let $s$ be any point in $S$ and let $s'$ denote a preimage in $S'$.  We can find an irreducible rational curve in $|\phi^{*}L|$ that avoids $s'$.  This divisor has vanishing intersection with every $(-2)$-curve, so it cannot intersect any such curve.  Thus this divisor avoids the entire fiber over $s$.  In this way we see that $L$ defines a birational morphism to $\mathbb{P}^{2}$.
\end{proof}

\begin{coro} \label{coro:goodmembersofsystem}
Let $S$ be a del Pezzo surface with canonical singularities and let $L$ denote a system of lines on $S$.  Fix a general point $s$ and let $|L|_{s} \subset |L|$ denote the sublinear series of divisors through $s$.  Then the general member of $|L|_{s}$ is a smooth rational curve which avoids $\Sing(S)$.  Any divisor parametrized by $|L|_{s}$ which avoids $\Sing(S)$ will be one of the following:
\begin{enumerate}
\item a smooth anticanonical cubic,
\item a $(-1)$-curve and a smooth anticanonical conic meeting transversally, or
\item a chain of three $(-1)$-curves meeting transversally.
\end{enumerate}
\end{coro}

\begin{proof}
The first statement follows from the fact that $|L|$ defines a birational morphism $\psi: S \to \mathbb{P}^{2}$.
We still must classify the possible types of divisors $D \in |L|_{s}$ which avoid $\Sing(S)$.  Note that every component will be a $(-1)$-curve, an anticanonical conic, or an anticanonical cubic.  Since $D$ avoids $\Sing(S)$, the image of $D$ in $\mathbb{P}^{2}$ is a line such that there are at most $2$ $\psi$-exceptional centers along $D$. The three cases in the theorem occur when $D$ meets $0$, $1$, or $2$ of these points.
\end{proof}

\subsection{Global systems of lines}

We now globalize our discussion from the previous section.

\begin{defi} \label{defi:globalsystem}
Let $\pi: \mathcal{X} \to B$ be a del Pezzo fibration.  We say that $\pi$ admits a global system of lines if:
\begin{enumerate}
\item Every fiber of $\pi$ is a normal del Pezzo surface with canonical singularities.
\item There is a line bundle $\mathcal{L}$ on $\mathcal{X}$ such that for every fiber $F$ of $\pi$ the restriction $\mathcal{L}|_{F}$ is a system of lines.
\end{enumerate}
\end{defi}

\cite{Ishii82} constructs a coarse moduli space parametrizing families of smooth del Pezzo surfaces which admit a global system of lines.  By \cite[Theorem 3]{Ishii82} this space can be compactified by allowing the surfaces to acquire $A_{1}$-singularities on the boundary.  Thus there are many del Pezzo fibrations satisfying the conditions of Definition \ref{defi:globalsystem}. 

Suppose that $\pi: \mathcal{X} \to B$ is a del Pezzo fibration such that every fiber is a del Pezzo surface with canonical singularities.  If $\mathcal{X}$ admits a line bundle $\mathcal{L}$ as in Definition \ref{defi:globalsystem} (2) then the restriction $\mathcal{L}|_{\mathcal{X}_{\eta}}$ defines a system of lines on the generic fiber that is defined over the ground field.  Conversely, suppose that $\mathcal{X}_{\eta}$ admits a birational morphism over $k(B)$ to projective space.  By taking the closure of an element of this linear system in $\mathcal{X}$, we obtain a Weil divisor $\mathcal{D}$ on $\mathcal{X}$ whose restriction to a general fiber is a system of lines.  This is not quite enough to verify Definition \ref{defi:globalsystem}, since $\mathcal{D}$ might not be Cartier.  If $\mathcal{X}$ is $\mathbb{Q}$-factorial, then since it is terminal Gorenstein \cite[Lemma 5.1]{Kawamata88} implies every divisor on $\mathcal{X}$ is Cartier.  In summary:

\begin{lemm} \label{lemm:linesexist}
Let $\pi: \mathcal{X} \to B$ be a del Pezzo fibration such that every fiber is a del Pezzo surface with canonical singularities.  Suppose that $\mathcal{X}$ is $\mathbb{Q}$-factorial and that $\mathcal{X}_{\eta}$ admits a birational morphism to $\mathbb{P}^{2}_{k(B)}$ over $k(B)$.  Then $\pi$ admits a global system of lines.
\end{lemm}

The key idea in this section is the following.  Suppose that $C_{1},C_{2}$ are numerically equivalent sections on $\mathcal{X}$.  We would like to show that after adding some $\pi$-vertical free curves to $C_{1}$ and $C_{2}$ we can obtain stable maps which lie in the same component of $\overline{M}_{g(B),0}(\mathcal{X})$.  Our strategy is to first find a well-behaved surface $Y$ containing $C_{1}$ and $C_{2}$ and then to construct the desired stable maps inside of $Y$.

The first step is to identify this well-behaved surface $Y$.

\begin{prop} \label{prop:wellbehavedsurface}
Let $\pi: \mathcal{X} \to B$ be a del Pezzo fibration which admits a global system of lines $\mathcal{L}$.  Suppose that $M_{1},M_{2}$ are components of $\Sec(\mathcal{X}/B)$ which generically parametrize relatively free sections which contain $\geq 2g(B)+1$ general points of $\mathcal{X}$. Then for general sections $C_{1}$ parametrized by $M_{1}$ and $C_{2}$ parametrized by $M_{2}$ there is a surface $Y \subset \mathcal{X}$ satisfying the following properties.
\begin{enumerate}
\item $Y$ contains $C_{1}$ and $C_{2}$.
\item $Y$ is contained in the smooth locus of $\mathcal{X}$ and its intersection with any fiber $F$ is contained in the smooth locus of $F$.
\item The restriction of $Y$ to any fiber $F$ of $\pi$ is a member of a system of lines on $F$.
\end{enumerate}
\end{prop}

\begin{proof}
Let $Z$ denote the union of the singular loci of all singular fibers of $\pi$ (so in particular $Z$ contains the singular locus of $\mathcal{X}$).  Note that $Z$ has codimension $3$ in $\mathcal{X}$.  Thus by Lemma \ref{lemm:expdimatpoints2} general choices of $C_{1}$ and $C_{2}$ will avoid $Z$.

Suppose we fix a singular fiber $F_{0}$ and fix a general section $C_{2}$.  In particular by Lemma \ref{lemm:expdimatpoints2} we may ensure that there are only finitely many divisors in the system of lines $|\mathcal{L}|_{F_{0}}|$ which contain $F_{0} \cap C_{2}$ and meet $Z$.  As we take the union of these divisors as we vary over all fibers $F_{0}$, we obtain a codimension $2$ subset $W$ of $\mathcal{X}$.  In particular, a general $C_{1}$ will not intersect $W$.

Consider the rational points on $\mathcal{X}_{\eta}$ corresponding to the general $C_{1},C_{2}$ as described in the previous paragraph.  There is a unique line in $\mathcal{X}_{\eta}$ in the system of lines $\mathcal{L}|_{\mathcal{X}_{\eta}}$ which contains these two points, and this line will be defined over the ground field.  Let $Y$ denote the corresponding surface in $\mathcal{X}$.  Then the restriction of $Y$ to any fiber $F$ will be the unique element in the system of lines $\mathcal{L}|_{F}$ connecting these two points.  In particular this $Y$ satisfies all the desired properties.
\end{proof}

The next lemma enables us to ``connect'' two sections of a fibration $\pi: Y \to B$ of relative dimension $1$.

\begin{lemm} \label{lemm:identifyingonsurface}
Let $Y$ be a smooth projective surface equipped with a morphism $\pi: Y \to B$ whose general fiber is $\mathbb{P}^{1}$.  Suppose that $C_{1}$ and $C_{2}$ are two sections of $\pi$.  Then there are stable maps $f_{1}: Z_{1} \to Y$ and $f_{2}: Z_{2} \to Y$ such that:
\begin{enumerate}
\item Both $Z_{1}$ and $Z_{2}$ consist of one genus $g(B)$ curve attached to trees of rational curves.
\item The map $f_{i}$ maps the genus $g(B)$ curve in $Z_{i}$ isomorphically to $C_{i}$ and the trees of rational curves to $\pi$-vertical curves.
\item Fix any point $b \in B$ and let $Z_{1,b}$ denote the part of $Z_{1}$ whose support maps to the fiber $F_{b}$ over $b$.  If $F_{b}$ is reducible then $f_{1*}Z_{1,b} < F_{b}$.  We also have the analogous statement for $f_{2}$.
\item $f_{1}$ and $f_{2}$ lie in the same component of $\overline{M}_{g(B),0}(Y)$.
\end{enumerate}
\end{lemm}

\begin{proof}
Choose any birational contraction $\rho: Y \to \mathbb{F}$ where $\mathbb{F}$ is a ruled surface.  There are $\pi$-vertical effective curves $T_{1},T_{2}$ such that $\rho^{*}\rho_{*}C_{1} = C_{1} + T_{1}$ and $\rho^{*}\rho_{*}C_{2} = C_{2} + T_{2}$.  An easy induction on the relative Picard rank of $\rho$ shows that $T_{1}$ is supported on reducible fibers and that the component of $T_{1}$ supported on any fiber is less effective than that fiber.  An analogous statement holds for $T_{2}$.

If we add on sufficiently many general fibers to $\rho_*C_{1}$ and $\rho_{*}C_{2}$ we can guarantee that the resulting curves are contained in the closure of the same component of $\Sec(\mathbb{F}/B)$ (see Proposition \ref{prop:totalbreaking}).  Note that a general section in this component will miss the $\rho$-exceptional centers, so that the strict transform will be the same as the pullback.  Thus we see that $\rho^{*}\rho_{*}C_{1}$ and $\rho^{*}\rho_{*}C_{2}$ are algebraically equivalent and both lie in the closure of a single component of $\Sec(Y/B)$.  By construction the portion of these $1$-cycles supported on any reducible fiber $F_{b}$ is less effective than the entire fiber.

Recall that $\Sec(Y/B)$ admits an embedding into $\overline{M}_{g(B),0}(Y)$.  By taking limits as the image curve approaches $\rho^{*}\rho_{*}C_{1}$ and $\rho^{*}\rho_{*}C_{2}$, we obtain stable maps $f_1 : Z_1 \to Y$ and $f_2 : Z_2 \to Y$ whose corresponding cycles are $\rho^{*}\rho_{*}C_{1}$ and $\rho^{*}\rho_{*}C_{2}$.  
These stable maps satisfy all the desired properties. 
\end{proof}

Finally, we prove that we can ``connect'' two sections in a del Pezzo fibration with a global system of lines.

\begin{lemm} \label{lemm:eventuallyequal}
Let $\pi: \mathcal{X} \to B$ be a del Pezzo fibration that admits a global system of lines.  Let $M_{1},M_{2}$ be components of $\Sec(\mathcal{X}/B)$ which generically parametrize relatively free sections which can contain $\geq 2g(B) + 1$ general points of $\mathcal{X}$.  Let $C_{1},C_{2}$ be general sections parametrized by $M_{1},M_{2}$ respectively.  There are stable maps $f_{1}: Z_{1} \to \mathcal{X}$ and $f_{2}: Z_{2} \to \mathcal{X}$ such that
\begin{enumerate}
\item the image of $f_{1}$ and the image of $f_{2}$ do not intersect the singular locus of $\mathcal{X}$ or the singular locus of any fiber $F$,
\item both $f_{1}$ and $f_{2}$ are smooth points of $\overline{M}_{g(B),0}(\mathcal{X})$,
\item the domain of $f_{1}$ is a comb which maps the handle to $C_{1}$ and the teeth to $\pi$-vertical free curves, and similarly for $f_{2}$, and
\item $f_{1}$ and $f_{2}$ lie in the same component of $\overline{M}_{g(B),0}(\mathcal{X})$.
\end{enumerate}
\end{lemm}

\begin{proof}
Using Proposition \ref{prop:wellbehavedsurface} we obtain a surface $Y$ containing both $C_{1}$ and $C_{2}$.  Note that $Y$ is contained in smooth locus of $\mathcal{X}$ and does not intersect the singular locus of any fiber of $\pi$.  Furthermore, the restriction of $Y$ to any fiber yields a system of lines.
Let $Y'$ denote a minimal resolution of $Y$.  By Corollary \ref{coro:goodmembersofsystem} the singular fibers of the map $Y \to B$ consist of a chain of two or three rational curves meeting transversally.  Thus $Y$ has only $A_{n}$ singularities.  This means that every fiber of $Y' \to B$ is a chain of rational curves.  Furthermore, the two curves on the end of this chain are not contracted by the birational map to $Y$.

Apply Lemma \ref{lemm:identifyingonsurface} to the strict transforms of $C_{1}$ and $C_{2}$ on $Y'$ to find stable maps $f'_{1}$ and $f'_{2}$ and let $\widetilde{f}_{1},\widetilde{f}_{2}$ be the corresponding stable maps to $\mathcal{X}$.  We claim that $\widetilde{f}_{1}$ and $\widetilde{f}_{2}$ represent smooth points of $\overline{M}_{g(B),0}(\mathcal{X})$.  Indeed, by combining Lemma \ref{lemm:identifyingonsurface} (3) with the explicit description of the fibers of $Y' \to B$ as chains of rational curves we see that the part of the image of $\widetilde{f}_{1}$ supported on any fiber of $Y \to B$ is less effective than the fiber itself, and the same is true of $\widetilde{f}_{2}$.  The classification of Corollary \ref{coro:goodmembersofsystem} yields the following possibilities for the vertical components of $\widetilde{f}_{1},\widetilde{f}_{2}$ contained in a fiber: a smooth anticanonical cubic, a smooth anticanonical conic, a $(-1)$-curve, or a union of two $(-1)$-curves meeting transversally.  Furthermore, since the locus of intersection points of $(-1)$-curves in fibers has codimension $2$ in $\mathcal{X}$ a general section in our families will avoid it.  
Thus in each of these cases $\widetilde{f}_{1},\widetilde{f}_{2}$ are local immersions near the neighborhood of any node and \cite[Lemma 2.6]{GHS03} applies to show the smoothness of the stable map.

The last step is to improve the properties of $\widetilde{f}_{1}$ and $\widetilde{f}_{2}$.  Consider the part of the image of $\widetilde{f}_{1}$ that is supported in a given fiber $F_{0}$.  We know that it is contained in the smooth locus of $F_{0}$.  In fact, since $F_{0}$ admits a system of lines there are many very free curves contained in the smooth locus of $F_{0}$.  By successively gluing on general very free curves in the smooth locus of $F_{0}$ to the various components of $\widetilde{f}_{1}$, we obtain another stable map $\widehat{f}_{1}$ such that the components in $F_{0}$ can be smoothed to a free curve while keeping the intersection point with the section fixed (see \cite[II.7.9 Theorem]{Kollar}).  We simultaneously glue members of these families of vertical free curves to general points of $\widetilde{f}_{2}$ to get $\widehat{f}_{2}$.   Since the original stable maps were smooth points of $\overline{M}_{g(B),0}(\mathcal{X})$ which lie in the same component, we can ensure that $\widehat{f}_{1}$ and $\widehat{f}_{2}$ will again be smooth points of $\overline{M}_{g(B),0}(\mathcal{X})$ which lie in the same component.  We then replace $\widetilde{f}_{2}$ by $\widehat{f}_{2}$ and $\widetilde{f}_{1}$ by a deformation of $\widehat{f}_{1}$ which smooths the part supported in $F_{0}$.  Repeating this process several times we obtain the desired statement.
\end{proof}

\subsection{Geometric Manin's Conjecture for del Pezzo fibrations with global systems of lines}

Let $\pi: \mathcal{X} \to B$ be a del Pezzo fibration.  We start by discussing a monoid action on the set of relatively free sections which comes from Movable Bend-and-Break.  

Fix a general fiber $F$ of $\pi$ and let $\mathcal{R}$ denote the set of components of $\overline{M}_{0,0}(F)$ which generically parametrize birational maps onto free rational curves.  If $F$ is a del Pezzo surface of degree $\geq 2$ then \cite[Theorem 5.1]{Testa09} combined with \cite[Lemma 2.13]{LTdPI} shows that every class in $\Nef_{1}(F)_{\mathbb{Z}}$ is represented by at most one element of $\mathcal{R}$.  Furthermore, two free curves on $F$ are guaranteed to intersect unless they are both fibers of the same map to $\mathbb{P}^{1}$.  With this one exception, if we choose two components of $\mathcal{R}$ then by gluing and smoothing general curves in this component we obtain a unique component of $\mathcal{R}$.  To handle the exception, for each family $R$ of anticanonical conics with square $0$ and for each integer $k \geq 2$ we formally add the element $kR$ to $\mathcal{R}$ -- conceptually, we think of these classes as ``gluing curves from $R$ $k$ times''.  With this addition and by formally adding in a $0$ class which acts trivially the gluing operation gives $\mathcal{R}$ the structure of a commutative monoid.

Let $\mathcal{M}$ denote the set of components of $\Sec(\mathcal{X}/B)$ which generically parametrize relatively free sections with height $\geq \MBB(\mathcal{X}) + 2$.  By Theorem \ref{theo:connectedfibers} the universal map over any component of $\mathcal{M}$ has connected fibers. Thus if we choose any component in $\mathcal{M}$ and any component of $\mathcal{R}$ then by gluing and smoothing we obtain a unique new component of $\mathcal{M}$.  Furthermore by \cite[Lemma 5.11]{LT17} if we perform this operation twice with two different families in $\mathcal{R}$ then the resulting component does not depend on which order we glue.  In this way $\mathcal{M}$ is given the structure of an $\mathcal{R}$-module.

\begin{lemm} \label{lemm:gmchelper}
There are finitely many components $M_{1},\ldots,M_{s} \in \mathcal{M}$ such that every component of $\mathcal{M}$ can be obtained by gluing $\pi$-vertical free curves to some $M_{i}$ and smoothing.
\end{lemm}

\begin{proof}
Let $M$ be any component of $\mathcal{M}$.  By Corollary \ref{coro:minbreaking} the closure of $M$ in $\overline{M}_{g(B),0}(\mathcal{X})$ contains a point representing the union of a relatively free section with a $\pi$-vertical curve of anticanonical degree $2$ or $3$.  By repeatedly breaking off such components, we see that we may define the $M_{i}$ to be the finite set of components of $\mathcal{M}$ which have anticanonical degree at least $\MBB(\mathcal{X}) + 2$ and no more than $\MBB(\mathcal{X}) + 5$.  
\end{proof}

It will be useful to reinterpret this lemma using the monoid action.  Let $\mathcal{S}$ denote the disjoint union of $s$ copies of $\mathcal{R}$.  Define the map $\xi: \mathcal{S} \to \mathcal{M}$ which sends a component $R$ in the $j$th copy of $\mathcal{R}$ to $R \cdot M_{j}$.  Then the lemma shows that this map is surjective.

We are now prepared to prove Geometric Manin's Conjecture for certain types of families.  Suppose that $\pi: \mathcal{X} \to B$ is a del Pezzo fibration which admits a global system of lines.  Then the relative Picard rank of $\pi$ is $\rho(\mathcal{X}_{\eta})$.  As in the discussion after Definition \ref{defi:countingfunction} there is a translate of $N_{1}(\mathcal{X}_{\eta}) \subset N_{1}(\mathcal{X})$ which contains all the classes of sections of $\pi$; we denote this translate by $N_{F}$.

\begin{theo}
Let $\pi: \mathcal{X} \to B$ be a del Pezzo fibration such that:
\begin{enumerate}
\item $\mathcal{X}_{\eta}$ is a del Pezzo surface of degree $\geq 2$ such that $\rho(\mathcal{X}_{\eta})$ coincides with the geometric Picard rank of $\mathcal{X}_{\eta}$, and
\item $\pi$ admits a global system of lines.
\end{enumerate}
There is a numerical class $\alpha \in N_{F}$ such that for any numerical class $\beta \in \alpha + \Nef_{1}(\mathcal{X}_{\eta})$ there is exactly one component of $\Sec(\mathcal{X}/B)$ which both represents $\beta$ and generically parametrizes relatively free sections.
\end{theo}

Combined with Lemma \ref{lemm:linesexist} this proves Theorem \ref{theo:maintheorem4}.

\begin{proof}
Let $i: F \to \mathcal{X}$ denote the inclusion of a general fiber.  By our assumption on the monodromy action we see that $\mathcal{X}$ has Picard rank $\rho(F) + 1$. Thus as demonstrated earlier every class in $\Nef_{1}(\mathcal{X}_{\eta}) \subset N_{1}(\mathcal{X})$ is represented by a unique family of free rational curves.  In particular, by gluing such curves onto a relatively free section and smoothing we obtain the existence of components of $\Sec(\mathcal{X}/B)$ representing numerical classes in a translation of $\Nef_1(\mathcal X_\eta)$, and it only remains to prove the uniqueness. 

Let $\mathcal{M}$ denote the set of components of $\Sec(\mathcal{X}/B)$ which generically parametrize relatively free sections with height $\geq \MBB(\mathcal{X}) + 2$.  Consider the map $\mu: \mathcal{M} \to N_{1}(\mathcal{X})_{\mathbb{Z}}$ which sends a component to its numerical class.  We would like to show that for some class $\alpha$ as in the statement of the theorem the fibers of $\mu$ over $\alpha + \Nef_{1}(\mathcal{X}_{\eta})$ are singletons.

Let $M_{1},\ldots,M_{s}$ be a finite set of components in $\mathcal{M}$ as in Lemma \ref{lemm:gmchelper} and let $\mathcal{S}$ denote the corresponding $\mathcal{R}$-module as defined above.  We have a surjective map $\xi: \mathcal{S} \to \mathcal{M}$.  Our strategy for understanding the fibers of $\mu$ is to study the fibers of $\xi$ and $\mu \circ \xi$ and to prove that they are often the same.

The first step is to better understand the fibers of $\xi$.  Suppose we fix our attention on the $i$th copy of $\mathcal{R}$ in $\mathcal{S}$, so that the restriction of $\xi$ sends $R \mapsto R \cdot M_{i}$.  Recall that there is a unique component of $\mathcal{R}$ representing any class of $\Nef_{1}(F)$.  Thus, the restrictions of both $\xi$ and $\mu \circ \xi$ to this  subset of $\mathcal{S}$ are injective.  Since $\mathcal{S}$ consists of $s$ copies of $\mathcal{R}$, we deduce that the fibers of $\xi$ and of $\mu \circ \xi$ have size at most $s$.

We will use the following claim to show that these fibers often are as large as possible.

\begin{clai}
Suppose that $M,M' \in \mathcal{M}$.  Then there are elements $R, R' \in \mathcal{R}$ such that
\begin{equation*}
R \cdot M = R' \cdot M'.
\end{equation*}
\end{clai}

\begin{proof}[Proof of claim:]
First, by adding on suitable elements of $\mathcal{R}$ and smoothing we may assume that $M, M'$ generically parametrize relatively free curves which can go through $\geq 2g(B) + 1$ general points of $\mathcal{X}$. (Note that if we prove our statement for these new components, the claim for the original components follows.)  We then apply Lemma \ref{lemm:eventuallyequal}, which immediately implies the desired claim.
\end{proof}

We now apply the claim to all possible pairs of elements from $M_{1},\ldots,M_{s}$.  In this way we obtain $R_{1},\ldots,R_{s}$ in $\mathcal{R}$ such that
\begin{equation*}
R_{1} \cdot M_{1} = R_{2} \cdot M_{2} = \ldots = R_{s} \cdot M_{s}.
\end{equation*}
Let $M'$ be this common family in $\mathcal{M}$ and let $\alpha$ be its numerical class.  Then for any $R \in \mathcal{R}$ the fiber of $\xi$ over $R \cdot M'$ is as large as possible: it contains $s$ elements, one in each component of $\mathcal{S}$.  (Precisely, the $i$th copy of $\mathcal{R}$ in $\mathcal{S}$ contributes the element $R \cdot R_{i}$ to this fiber.)  This of course means that the fibers of $\mu \circ \xi$ over $\alpha + R$ also have the maximal size $s$.  These two statements together show that for any $R \in \mathcal{R}$ the map $\mu$ is injective over the class $\alpha + R$.  Since the numerical class map from $\mathcal{R}$ to $\Nef_{1}(F)_{\mathbb{Z}}$ is surjective, we have proved the desired statement.
\end{proof}

\bibliographystyle{alpha}
\bibliography{Sections2}

\end{document}